\documentclass{article}

\usepackage{setspace}
\usepackage{ulem}
\usepackage{amssymb,amsmath,amsthm,eucal,todonotes}
\usepackage{bm}
\usepackage[T1]{fontenc}
\usepackage[utf8]{inputenc}
\usepackage[letterpaper,top=2cm,bottom=2cm,left=1.5cm,right=1.5cm,marginparwidth=1.75cm]{geometry}

\usepackage{graphicx}
\usepackage{amsfonts}
\usepackage{float}
\usepackage{subcaption}
\usepackage{algorithm}
\usepackage{algorithmic}
\usepackage{xcolor}
\usepackage{color}
\usepackage{lineno}
\usepackage{todonotes}
\definecolor{darkGreen}{RGB}{0,100,0}

\usepackage{mathtools}

\def \x {{\bf x}}
\def \b {{\bf b}}
\def \r {{\bf r}}

\def\RR{\mathbb{R}}

\def \pmatrix{ \left( \begin{array} }
\def \endpmatrix{ \end{array} \right) }
\def \u {{\bf u}}
\def \y {{\bf y}}
\def \v {{\bf v}}

\def \argmin {\mbox{argmin}}

\newtheorem{df}{Definition}[section]

\newtheorem{rmk}[df]{Remark}

\title{Piecewise DMD for oscillatory and Turing spatio-temporal dynamics}

\author{Alessandro Alla$^{1}$, Angela Monti$^{2,3}$, Ivonne Sgura$^{2}$\\
\small $^{1}$ Università Ca' Foscari Venezia, Dipartimento di Scienze Molecolari e Nanosistemi, Venezia, Italy,\\
\small e-mail: alessandro.alla@unive.it\\
\small $^{2}$ Università del Salento, Dipartimento di Matematica e Fisica ‘‘E. De Giorgi’’, Lecce, Italy,\\ 
\small e-mail: \{angela.monti, ivonne.sgura\}@unisalento.it \\
\small $^{3}$Istituto per le Applicazioni del Calcolo ``M. Picone'',
CNR, Bari, Italy\\
\small e-mail: a.monti@ba.iac.cnr.it}

\begin{document}
\maketitle

\begin{abstract}

Dynamic Mode Decomposition (DMD) is an equation-free method that aims at reconstructing the best linear fit from temporal datasets. In this paper, we show that DMD does not provide accurate approximation for datasets describing oscillatory dynamics, like spiral waves and relaxation oscillations, or spatio-temporal Turing instability.
Inspired from the classical "divide  and conquer" approach, we propose a piecewise version of DMD (pDMD) to overcome this problem.
The main idea is to split the original dataset in $N$ submatrices and then apply the exact (randomized) DMD method in each subset of the obtained partition.
We describe the pDMD algorithm in detail and we introduce some error indicators to evaluate its performance when $N$ is increased. 
Numerical experiments show that very accurate reconstructions are obtained by pDMD for datasets arising from time snapshots of some reaction-diffusion PDE systems, like the FitzHugh-Nagumo model, the $\lambda$-$\omega$ system and the DIB morpho-chemical system for battery modeling.

\end{abstract}


\section{Introduction}
The large amount of temporal datasets has increased in the last decade and with that the study of hidden structures. Specifically, mathematical models can describe rigorously datasets where no information is provided. A data-driven model can help to understand physical phenomena and forecast future development. 
Recently, the use of machine learning techniques has further improved the capability to discovery mathematical models and has significantly enlarged this area of research. The literature on this topic is rather spread. Here we, first, recall two methods based on different strategies: sparse optimization and neural networks. Finally, we will revoke the Dynamic Mode Decomposition (DMD) which is the building block for the method used in this manuscript.

One technique goes back to 2016 where the authors in \cite{BPK16} have used sparse optimization methods, e.g. Lasso algorithm, to recover Ordinary Differential Equations (ODEs). The method relies on a large library including elements that may appear in the model and through a sparse optimization algorithm it is possible to discover the weights of the terms in the library. Thus, if a coefficient is zero, the corresponding element in the library does not appear in the model. This also justifies the use of sparse optimization methods, because only few terms are required in the searched model. 
Later, this method has been extended in \cite{RBPK17} to PDEs with constant coefficients. The library in this case was built from derivatives of the dataset to include further terms in the model. Finally, extension to non constant coefficients for PDEs was presented in \cite{RABK19}.

Another method to construct data-driven modeling is built on Deep Neural Networks (DNNs). Specifically, a new class of DNNs, namely Physics-Informed Neural Networks (PINNs), is trained to solve supervised learning using information from the hidden physical law one wants to discover. The physical law may be described by an ODE or a PDE model. The great novelty in PINNs is the use of the physical laws together with the mean square error for dataset in the minimization of loss function. Therefore, the output of this method provides the coefficients needed to discover the model. In this way, the method is forced to converge to the model and to consider the physics behind the dataset. PINNs was introduced in \cite{RPK19, KKLPWY21,YLK22}.

Another powerful technique for data-driven modeling is the Dynamic Mode Decomposition (DMD).
DMD was firstly introduced in \cite{S10} and its algorithm finds the best linear fit model without explicit knowledge of the dynamics hidden in the data.
Specifically, the DMD algorithm determines eigenvalues and eigenvectors of an approximate linear model.  Later, in \cite{TRDLBK14} the authors introduced the exact DMD based on a low rank approximation of the original method. 
Afterwards other algorithms have been proposed to further improve the above DMD method; we refer to e.g. the optimal DMD \cite{HH17} and the higher order DMD (HODMD, \cite{LV17, LV17b}) which is based on the fitting of multilinear models. Despite the tremendous effort to improve the method and its computational efficiency, there are (at least) two classes of datasets where DMD still does not work properly, to the best of authors' knowledge.

Specifically, periodic datasets and Turing instability may lead to a wrong (in the first case) or a not accurate (in the second case) DMD reconstruction, as shown in \cite{VBGRC22, BMS21} and later in Section \ref{sec_fail_dmd} of this paper.

In \cite{VBGRC22}, the authors show that the DMD algorithm for systems is more accurate when a unique dataset is considered instead of a dataset for each variable. They refer to this approach as coupled or uncoupled DMD, respectively. Furthermore, in \cite{VBGRC22}, it is mentioned that DMD may fail in the case of reaction diffusion systems.
In \cite{BMS21}, the uncoupled DMD implementation was applied to Turing instability dynamics leading to stationary pattern solutions. The authors have shown that DMD cannot reconstruct accurately that datasets even if large ranks are used in the algorithm. 

The main goal of this paper is to propose an alternative DMD algorithm to face with these drawbacks. Stemming on the ``divide and conquer'' principle, well known in the numerical analysis framework, we introduce a {\it piecewise} version of the exact DMD, that we will define as pDMD.
In the first part of the paper, we will present in details which drawbacks arise when the exact DMD is applied to reconstruct the spatio--temporal dynamics of a selection of Reaction-Diffusion PDE (RD--PDE) systems in two variables. Following the results in \cite{VBGRC22}, in this manuscript, we always use a coupled DMD approach. Furthermore, a randomized DMD algorithm based on the $QB$ decomposition (see \cite{EMKB19}) is applied to reduce the computational costs.

In particular, we consider the following RD--PDE systems whose solutions exhibit an oscillatory dynamics: i) the FitzHugh-Nagumo model \cite{CS10} 1D in space and with relaxation oscillations in time leading to a limit cycle in the phase plane; ii) the $\lambda$-$\omega$ system from \cite{MV06, BG2006} with spiral waves solutions. 
  In the first case, also for the full rank approximation, the oscillatory solutions and the corresponding limit cycle cannot be recovered at all. In the second case, for a certain range of ranks, the final spiral wave is approximated with low accuracy, but the amplitude and phase of the oscillating behaviour are not preserved. Moreover, for larger ranks, where a better approximation is expected, ill-conditioning of the fitting procedure behind DMD emerges and the time dynamics is lost also at the final time. These results are reported in Section \ref{sec:periodic}.
  
To deepen the discussion started in \cite{BMS21} about DMD defects, we apply here the ``coupled'' exact randomized DMD to follow the Turing dynamics of the morphochemical RD--PDE model, introduced in \cite{DIB13} and known as DIB model. In Section \ref{sec:turing}, we deal with the Turing instability, where a transient unstable regime is present (said reactivity zone) before reaching the spatially inhomogeneous Turing pattern at the steady state (stabilizing zone). Finally, in Section \ref{ex_dib_hopf}, we also consider an example of Turing--Hopf instability \cite{DIB15}, that is, there exists an interplay between Turing and Hopf instabilities, where, after the initial unstable behaviour, the solutions are patterns that oscillate both in space and time.

In the first case (Section \ref{sec:turing}), we show that after a certain rank the DMD error dramatically increases and indeed blows up, because also here the ill--coditioning of the fitting procedure appears. Nevertheless for small ranks it is possible to reconstruct with few accuracy the time dynamics and the final pattern (see e.g. Figure \ref{dib_mean_classic}). For the more complicated Turing-Hopf dynamics (Section \ref{ex_dib_hopf}), that, at best of author's knowledge, has not been studied so far by DMD, very inaccurate approximations are obtained until the full rank choice and the best case exhibits large errors both for the final pattern and all the time history, as documented by the limit cycle approximation (see e.g Figures \ref{dib_hopf_phase}, \ref{dib_hopf_dmd_classic}). 

Therefore, our proposed pDMD algorithm will work as follows. Given a dataset $S$, corresponding e.g. to periodic or Turing dynamics in the time interval $[0,T]$, and a tolerance $\overline{tol} >0$, we start with DMD on the whole dataset. If the obtained error, that is the maximum of the worst-approximation in time between the dataset and its DMD reconstruction, is above the threshold $\overline{tol}$ we split $S=S_1 \cup S_2$ into $N=2$ parts and compute DMD in each submatrix $S_i$. This error indicator is checked at each iteration and, if the $i-$th dataset is not reconstructed accurately, we directly increase the number of subdivisions $N$. 

We iterate this splitting in $N$ parts till we reach the desired accuracy on all subsets $S_i, i=1, \dots, N$. This procedure will identify the first acceptable partition size, say $N^*$. 
In a second step, we decide to increase the partition size for $N > N^*$ in order to look for a sort of convergence along the whole time dynamics, by controlling another error indicator in the Frobenius norm accounting for the whole time DMD reconstruction. 
This further piecewise iteration will tend to improve an initial good reconstruction, as we will show in our numerical experiments. More details will be discussed in the formalized algorithms reported in Section \ref{sec:pdmd} and in Section \ref{sec:testpdmd}.
It is worth remarking that our pDMD can work with several datasets of small sizes, if many subintervals are required, and this allows to better follow "locally" the dynamics of the problem and then to better capture its behaviour with a linear regression method. 

An extensive numerical study of our algorithm is deeply discussed by considering the snapshot matrices obtained by the numerical solutions of the RD--PDE models described above. We will show that the pDMD increases the accuracy of the approximation removing the drawbacks discussed.

The paper is organized as follows. Section \ref{sec:dmd} recalls the exact DMD method and its variant based on the randomized QB decomposition. In Section \ref{sec:fom}, we briefly introduce the general reaction-diffusion PDE system in exam and the IMEX Euler method in matrix form for its approximation, because the obtained numerical solutions are used to build our datasets. As discussed above, in Section \ref{sec_fail_dmd}, we show how the exact (randomized) DMD fails on different datasets for periodic data or Turing dynamics. The pDMD is introduced in Section \ref{sec:pdmd} together with a complete description of our algorithm. Finally, in Section \ref{sec:testpdmd}, we present our numerical results. Conclusions are drawn in Section \ref{sec:end}. All numerical simulations have been performed in MATLAB (ver. 2019a) on a computer DELL, i7 Intel Core processor 2.8 GHz and 16Gb RAM.
%


\section{Dynamic Mode Decomposition}\label{sec:dmd}
The Dynamic Mode Decomposition (DMD) technique aims at reconstructing the best linear dynamical system hidden in a given temporal dataset $S = [\x_0, \x_1, \dots, \x_m]\in\mathbb{R}^{n \times (m+1)}$ where the $i-$th column of the matrix $S$ corresponds to the data at time $t_{i+1}$ and $t_{i+1}>t_i, i = 0, \dots, m-1.$ DMD fits the following linear model on these data:
\begin{align}\label{hid:dyn}
\begin{aligned}
\dot{\y}(t) &= A\y(t)\quad t\in [0,T],\\
\y(0) &= \y_0\in\RR^{n},
\end{aligned}
\end{align}
where $A\in\mathbb{R}^{n \times n}$ is an unknown operator, $\y(t):[0,T]\rightarrow \mathbb{R}^{n}$, the initial condition $\y_0$ coincides with the first element of the dataset $\x_0$ and the data are such that $\x_i \approx \y(t_i), i=1, \dots, m.$ 
To discover the matrix $A,$ DMD starts by splitting the snapshot matrix $S$ into two matrices
$$ S_L= \begin{bmatrix} | & | & \dots & | \\
\x_0 & \x_1 & \dots & \x_{m-1} \\
| & | & \dots & |  
\end{bmatrix} \in \mathbb{R}^{n \times m}, \quad S_R = \begin{bmatrix} | & | & \dots & | \\
\x_1 & \x_2 & \dots & \x_m \\
| & | & \dots & |  
\end{bmatrix} \in \mathbb{R}^{n \times m}$$
and then, assuming that there exists a linear relation between $\x_{k+1}$ and $\x_k$ for $k=0, \ldots, m-1,$ tries to find the best fitting matrix $A\in\mathbb{R}^{n \times n}$ such that 
\begin{equation}
S_R = AS_L, \quad   \Longleftrightarrow \quad  \x_{k+1}=A \x_k, \ k = 0, \dots, m-1.
\label{DMDfit}
\end{equation}
%
Therefore, for this purpose, the following least squares optimization problem can be formulated
\begin{equation}\label{opt_dmd}
A := \underset{\mathcal{A}\in\mathbb{R}^{n \times n}}{\argmin} \|S_R - \mathcal{A}S_L\|_F
\end{equation}
where $\| \cdot \|_F$ is the Frobenius norm.
It is well known that \eqref{opt_dmd} can be solved by computing the Moore-Penrose pseudo-inverse of $S_L$, such that the best fit solution is given by $A = S_R S_L^{ \dagger}.$ However, the dimension $n$ of the problem may be large and, from a computational point of view, it is not convenient to calculate $A$ directly by the above product.
For this reason, the exact DMD algorithm proposed in \cite{TRDLBK14} computes a low rank approximation adding a rank constraint to the optimization problem \eqref{opt_dmd}, as follows:

\begin{equation}\label{A_dmdr}
\begin{aligned}
A:=& \underset{\mathcal{A}\in\mathbb{R}^{n \times n},\, \text{rank}(\mathcal{A})=r}{\argmin}
& &  \|S_R - \mathcal{A}S_L\|_F \\
\end{aligned}
\end{equation}
%
where usually $r\ll n.$ 
%
Instead of solving directly the rank constrained optimization problem \eqref{A_dmdr} which is hard to tackle, DMD starts by computing a reduced SVD of the matrix $S_L \approx \Psi_r \Sigma_r V_{r}^T$, where $\Sigma_r \in \mathbb{R}^{r \times r}$ is a diagonal matrix whose entries $\sigma_i \geq 0$ are the first $r$ singular values of $S_L$ sorted in a decreasing order, $\Psi_r \in \mathbb{R}^{n \times r}$ and $V_r \in \mathbb{R}^{m \times r}$ are orthogonal matrices. Therefore, an approximation of the full matrix $A$ can be obtained by computing the pseudoinverse of the rank reduced $S_L,$ that is $S_L ^{\dagger} \approx V_r \Sigma_r^{-1} \Psi_r^T,$ such that
$$ A \approx S_R V_r \Sigma_r^{-1} \Psi_r^T $$
and then by projecting it onto the POD modes (the $r$ leading left singular vectors $\Psi_r$) it is possible to compute the reduced matrix  
\begin{equation}\label{eq:atilde}
 \tilde{A} = \Psi_{r}^{T} S_R  V_r \Sigma_{r}^{-1} \in \mathbb{R}^{r \times r}.
 \end{equation}
We observe that $\tilde{A}$ has the same leading $r$ eigenvalues of $A$. 
Thus, we compute the spectral decomposition of $\tilde{A}$, $ \tilde{A} W = W \Lambda $
where the columns of $W$ are the eigenvectors of $\tilde{A}$ and $\Lambda= diag(\lambda_1, \dots, \lambda_r)$ is a diagonal matrix containing the corresponding leading $r$ eigenvalues of the full matrix $A$. 
Finally, we reconstruct the high-dimensional DMD modes of $A$ by
$$\Phi = S_R V_r \Sigma_{r}^{-1} W \in \RR^{n \times r}.$$
We observe that these DMD modes are eigenvectors of the matrix $A$, with corresponding eigenvalues $\Lambda$:
\begin{equation*}
A \Phi = (S_R V_r \Sigma_r^{-1} \Psi_r^T)(S_R V_r \Sigma_r^{-1} W) = S_R V_r \Sigma_r^{-1} \tilde{A} W = S_R V_r \Sigma_r^{-1} W \Lambda  = \Phi \Lambda
\end{equation*}
and reconstruct the state variable as
\begin{equation}
\label{DMD_sol}
\x_k \approx \widetilde \x_k:=\sum_{i=1}^{r} \bm{\phi}_i \lambda_i^{k} b_i  = \Phi \Lambda^k \mathbf{b}
\end{equation}
where the DMD modes $\boldsymbol{\phi}_i$ are the columns of $\Phi$ (eigenvectors of $A$), $\lambda_i$ are the corresponding eigenvalues, while $b_i$ can be obtained by solving in the least squares sense the overdetermined system $\Phi  \mathbf{b} = \sum_{i=1}^{r} \bm{\phi}_i b_i = \x_0, $ where $\x_0$ is the first snapshot. The notation $\widetilde\x_k$ will be the DMD reconstruction of $\x_k$. \\
We briefly summarize these steps in Algorithm \ref{DMD}, defined as the exact DMD in \cite{TRDLBK14}.
\begin{algorithm}[h!]
\caption{(Exact) DMD}
\label{DMD}
\begin{algorithmic}[1]
\STATE {\bf{INPUT}} Snapshots $\{\x_0, \x_1, \dots, \x_m\}$, rank $r$ 
\STATE {\bf{OUTPUT}} DMD modes $\{\boldsymbol{\phi}_1, \dots, \boldsymbol{\phi}_r\}$, eigenvalues $\{\lambda_1,\ldots, \lambda_r\}$, reduced solution $\{\widetilde \x_k\}_{k=0}^m$
\STATE Set $S_L = [\x_0, \dots, \x_{m-1}], S_R = [\x_1, \dots, \x_m]$
\STATE Compute the truncated SVD of $S_L$, $S_L \approx \Psi_r \Sigma_r V_r^T$
\STATE Compute $\tilde{A} = \Psi_{r}^{T} S_R  V_{r} \Sigma_{r}^{-1}$
\STATE Compute the spectral decomposition of $\tilde{A}$, $\tilde{A} W = W \Lambda$
\STATE Calculate the DMD modes as $\Phi = S_R V_r \Sigma_{r}^{-1} W$
\STATE  Set $\b=\Phi^\dagger \x_0$
\STATE Set $\widetilde\x_k = \sum_{i=1}^{r} \bm{\phi}_i \lambda_i^{k} b_i , \quad k=0,\ldots, m$ 
\STATE $\widetilde S(r)=[\widetilde \x_0, \widetilde \x_1, \dots, \widetilde \x_m]$ 
\end{algorithmic}
\end{algorithm}
\subsection{Randomized DMD}
%
The exact DMD method as presented in Algorithm \ref{DMD} can be still computationally very expensive if $n$ is very large. Therefore, in this subsection we recall the method introduced in \cite{EMKB19} based on the randomized $QB$ decomposition. The aim is to write the dataset $S\approx QB$ where $Q\in\RR^{n \times r}$ is an orthogonal matrix and  $B\in\RR^{r\times (m+1)}$.
Once the QB decomposition is obtained for the snapshot matrix $S$, one will directly apply the DMD algorithm to the matrix $B=[\b_0,\ldots, \b_m]\in\RR^{r\times (m+1)}$ which is clearly much smaller than the original matrix $S$. In this way, after splitting the matrix into $B_L=[\b_0,\ldots, \b_{m-1}]$ and $B_R=[\b_1,\ldots, \b_{m}]$ we can solve the optimization problem

\begin{equation}\label{opt:qb}
\widetilde{A}:=\underset{\mathcal{\widetilde A}\in\mathbb{R}^{r\times r}}{\argmin} \|B_R - \mathcal{\widetilde A} B_L\|_F.
\end{equation}
It turns out that the solution of \eqref{opt:qb} is $\widetilde A=B_R B_L^\dagger$.  We note that $\widetilde A$ is a low rank approximation of $A$ and differs from \eqref{eq:atilde}. This computation is now doable since the dimensions of the matrices $B_L$ and $B_R$ are (eventually) very small. Then, we can compute the eigenvalue decomposition of the matrix $\tilde A W = W\Lambda$ and set the DMD modes $\Phi = Q B_R V\Sigma^{-1}W$. Therefore, the reduced solution $\x_k, \forall k$ can be obtained as in Algorithm \ref{DMD}. \\
In the remainder of this section we recall how to obtain the randomized QB decomposition. \\ First of all, one has to choose the so called {\it target rank} $r$ and the number of oversampling $p$ usually $5\leq p \leq 10.$ Then, we generate a random test matrix $\Omega \in\RR^{(m+1) \times \ell}$ with $\ell = r + p$ drawn from the normal Gaussian distribution. The oversampling needs to guarantee the target rank $r$, in fact it is common to build slightly larger test matrix to obtain improved basis.
The sampling matrix can be computed as $Y = S\Omega\in \RR^{n \times \ell} $ or by using the power iteration method as $Y = ((SS^T)^q S)\Omega, q\in\mathbb{N}.$
The matrix $Q$ is then obtained from the $QR$ decomposition of the sampling matrix $Y$. It is shown in \cite{EMKB19} that the power iteration improves the quality of the approximated basis matrix $Q$ using just one or two iterations, i.e. $q=\{1,2\}$.
Finally, the low rank matrix $B$ will be such that $B = Q^T S$. The QB algorithm is summarized in Algorithm \ref{alg:qb} using Matlab notations.
\begin{algorithm}[h!]
\caption{Randomized QB decomposition}
\label{alg:qb}
\begin{algorithmic}[1]
\STATE {\bf{INPUT}} Snapshots $S$, target rank $r$, oversampling $p$, number of power iterations $q$
\STATE {\bf{OUTPUT}} $Q\in\RR^{n \times r}, B\in\RR^{r\times (m+1)}$
\STATE $\ell = r + p$
\STATE $\Omega = {\tt{rand}}(m+1, \ell)$
\STATE $Y = S\Omega$
\FOR{j=1,\ldots,q}
\STATE $[Q, \cdot] = {\tt qr}(Y)$
\STATE $[Z, \cdot] = {\tt qr}(S^T Q)$
\STATE $Y = SZ$
\ENDFOR
\STATE $[Q, \cdot]= {\tt qr}(Y)$
\STATE $B = Q^TS$
\end{algorithmic}
\end{algorithm}
\section{Full model and its numerical approximation} 
\label{sec:fom}
In this paper, our aim is to apply DMD to reconstruct in time both oscillatory dynamics, leading for example to relaxation oscillations and spiral waves, and Turing pattern formation dynamics, that presents a transient unstable regime, known as reactivity, before reaching a structured spatially inhomogeneous pattern as stationary solution at the steady state. A common feature is that all these time behaviours can characterize the solutions of a RD--PDE system for different choices of the involved parameters. More details can be found e.g. in \cite{Murray03book}.

For this reason, here we consider the following general RD--PDE system:
\begin{equation}
\label{RDPDE}
\begin{cases}
u_t = d_u \Delta u + f(u,v) , \quad (x,y) \in \Omega \subset 
\mathbb{R}^d , \quad t \in (0,T], \\
v_t = d_v \Delta v + g(u,v), \\
(\mathbf{n} \nabla u)_{| \partial \Omega} = b_u(t), \quad (\mathbf{n} \nabla v)_{| \partial 
\Omega}  = b_v(t), \\
u(x,y,0) = u_0(x,y) , \quad v(x,y,0) = v_0(x,y).
\end{cases}
\end{equation}
where $d_u,d_v\in\RR^+$ are the diffusion coefficients, $T>0$ the final time of integration,  $d \in \{1,2\}$ is the space dimension of \eqref{RDPDE}.
The nonlinear reaction terms $f, g : \RR^2 \rightarrow \RR$ account for biological, chemical and other kind of phenomena.  We will consider Neumann boundary conditions, where $\mathbf{n}$ denotes the exterior normal to the boundary $\partial \Omega$ and $b_u(t), b_v(t): [0,T] \rightarrow \RR $ are scalar functions, identically zero in the case of homogeneous Neumann BCs.
Our aim is to solve numerically \eqref{RDPDE} to generate a dataset $S$ of our interest and then apply the DMD directly to this dataset without using any extra information coming from the (known) PDE. 

The model \eqref{RDPDE} depends on various parameters that will be chosen {\it ad hoc} in order to study different kinds of dynamics, as follows.
In fact, in the next sections we will consider: 1) the FitzHugh-Nagumo model \cite{CS10} with  relaxation oscillations and related limit cycle and 2) a $\lambda$-$\omega$ system \cite{MV06, BG2006} with spiral waves. Concerning Turing pattern formation, we consider 
the DIB morphochemical model \cite{DIB13, DIB15, DIB17sphere, SLB19}. Furthermore, we study this model also in presence of a combination of oscillatory and Turing behaviours, arising from the so-called  Turing-Hopf patterns, that are spatial inhomogeneous Turing patterns oscillating both in space and time (\cite{DIB15, SS2016}). 

For all the above models, the construction of the snapshot matrix $S$, to feed the DMD method, follows from the numerical approximation of \eqref{RDPDE}. 
Hence, for the spatial semi-discretization we apply standard finite differences with a total number $n$ of meshpoints inside $\Omega$. For the approximation in time, we apply the IMEX Euler scheme (i.e. we treat implicitly the diffusion part and explicitly the nonlinear reaction terms) on the meshgrid $\tau_{i+1}=\tau_i + h_t$ $i=0, \dots, n_T$ with timestep $h_t=T/n_T$. To simulate oscillatory solutions and Turing patterns we require both fine spatial meshes and integration for long times ($T \gg 1$) to attain the {\it standing} asymptotic oscillations (i.e. the limit cycle in the phase space) or the stationary Turing pattern. For this reason the computational load of the usual {\it vector approach} solving a large sparse linear system at each time step for the IMEX Euler method, can be very expensive. Then, here we apply the recent {\it matrix-oriented} approach and in particular the {\it rEuler} method proposed in \cite{DSS20}, solving at each timestep a Sylvester matrix equation in the reduced spectral space. 

To further reduce the computational cost, we will store only some of the computed snapshots on a temporal sub-grid $t_{i+1} = t_i + \kappa h_t$, where $\kappa\in\mathbb{N}$ allows to select equidistributed snapshots from the original grid $\{\tau_i\}_{i=0}^{n_T}$. Note that, if $\kappa=1$ we consider all the snapshots, if $\kappa=4$ we store snapshots every 4 time steps from the original grid. Hence, for the simulations presented in the next sections, we build the snapshot matrix $S \in \RR^{2n \times (m+1)}$, where $m+1 =n_T/\kappa$ and, for $i=0,\ldots,m$, the $(i+1)$-th column is the extended vector $\x_{i+1}:=[\u_{i+1}; \v_{i+1}]  \in \RR^{2n}$ given by the concatenation of the numerical solutions for both unknowns $u$ and $v$, i.e. $ \u_{i+1}\approx u(t_i),  \v_{i+1}\approx v(t_i)$.
 

%
\section{Drawbacks of the DMD method}
\label{sec_fail_dmd}
The examples discussed in this section have in common that the exact randomized DMD method, recalled in Section \ref{sec:dmd}, does not approximate the dataset accurately. In the first two examples for periodic datasets, we will see that DMD completely fails even with a full rank approximation. In the case of Turing instability dynamics, we will show that DMD exhibits poor approximation with also an error behaviour dramatically increasing with the rank $r$ in case of stationary patterns.

In all tests presented, to measure the quality of the DMD approximation of rank $r$, we use the following relative error in the Frobenius norm between the dataset $S$ and its DMD reconstruction $\widetilde{S}$:
\begin{equation}
\label{ef_data}
\mathcal{E}(\widetilde{S},r) =  \frac{\| S-\widetilde{S}\|_F}{\|S \|_F}.
\end{equation}
%

Furthermore, for each kind of dynamics considered, we are interested in comparing the behaviour in time of the full dataset with that approximated by the DMD for a given rank $r$. For this reason, we will compare the time behaviour of the spatial mean of the full and reduced solutions,  defined for $u$ by
\begin{equation}
\label{mean}
\langle u(t) \rangle := \frac{1}{|\Omega|} \int_\Omega u(x,y,t) dx dy \approx 
\text{mean}(\u_k), \quad  \quad k =0, \dots, m 
\end{equation}
and $\text{mean}(\widetilde\u_k)$, $k=0, \dots, m$, respectively. Similar computations are done for the variable $v$. In all examples shown in this section, DMD is performed with the QB decomposition from Algorithm \ref{alg:qb} and the results presented are referred to the value $r$ which minimizes $\mathcal{E}(\widetilde S,r)$, as can be extracted by the corresponding reported figures.

\subsection{Examples with periodic datasets}
\label{sec:periodic}
%
The first example will focus on the approximation of the limit cycle generated by the FitzHugh-Nagumo model \cite{CS10}. The second example concerns the reconstruction of the spiral wave solution and dynamics of the $\lambda$-$\omega$ system in \cite{MV06, BG2006}. As already outlined, we are interested in the reconstruction of the whole spatio-temporal history.
\subsubsection{FitzHugh-Nagumo model: limit cycle}
\label{ex_FHN}
\begin{figure}[htbp]
\centering
 \includegraphics[scale=0.4]{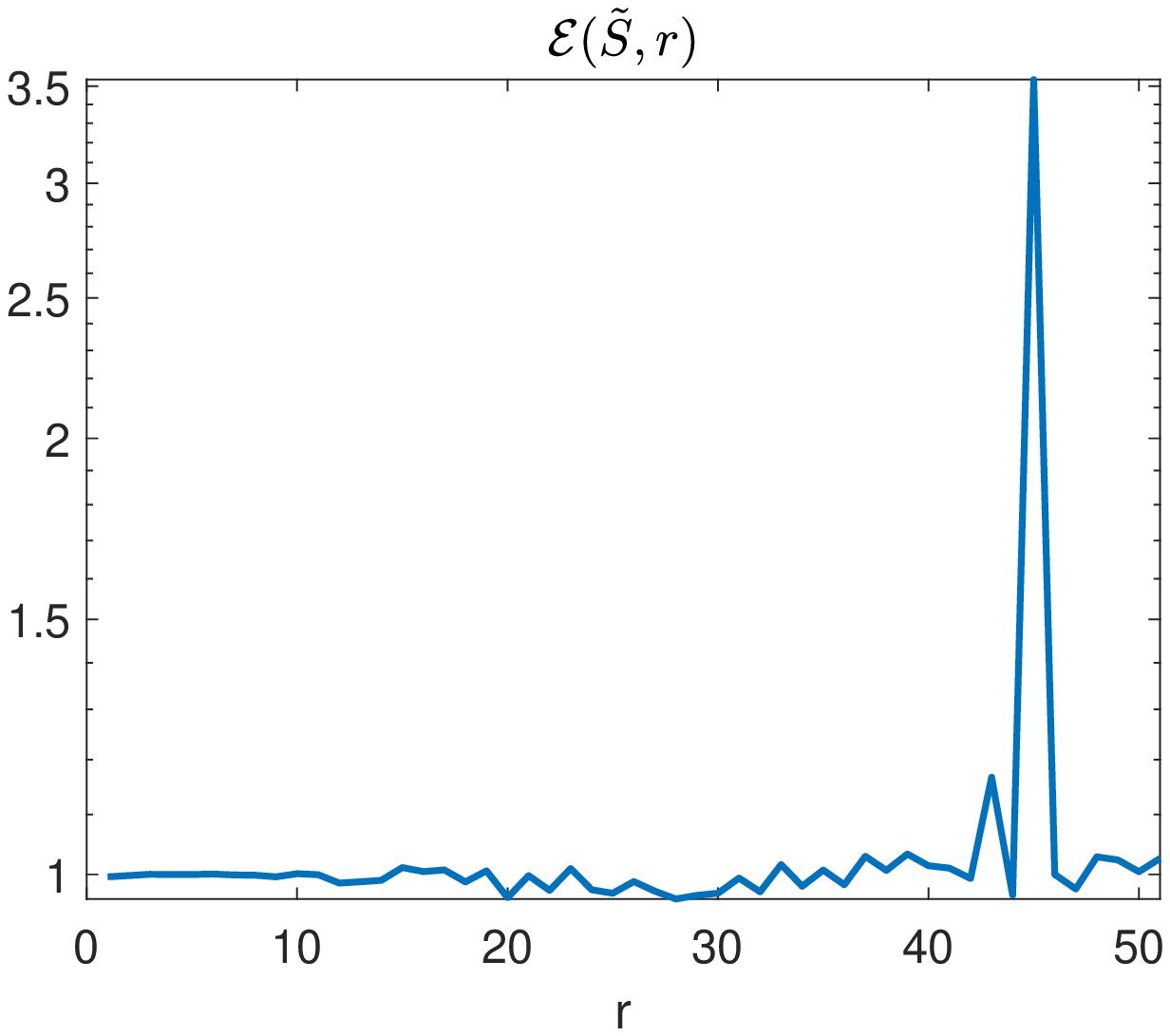}
\includegraphics[scale=0.4]{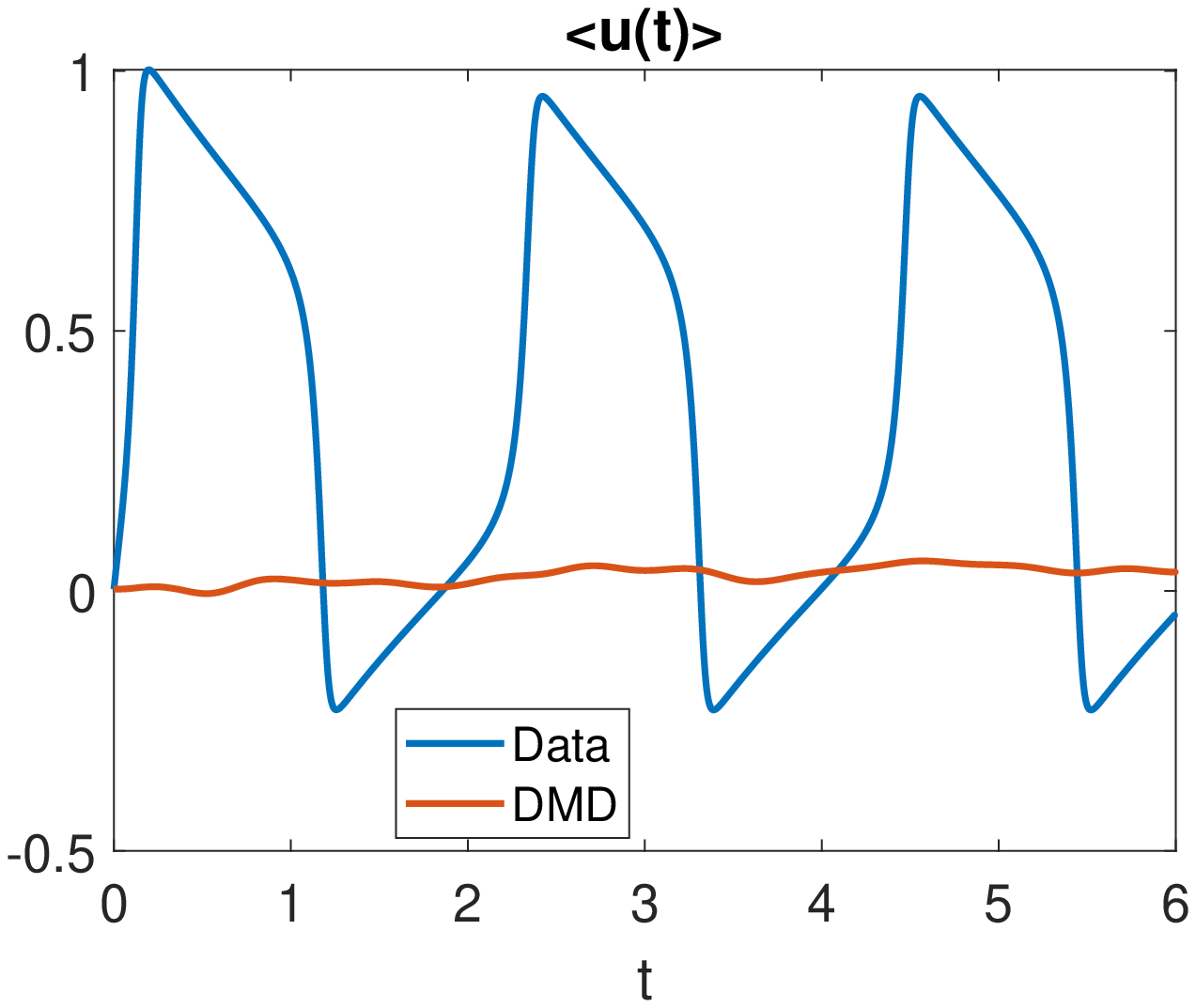}
\includegraphics[scale=0.4]{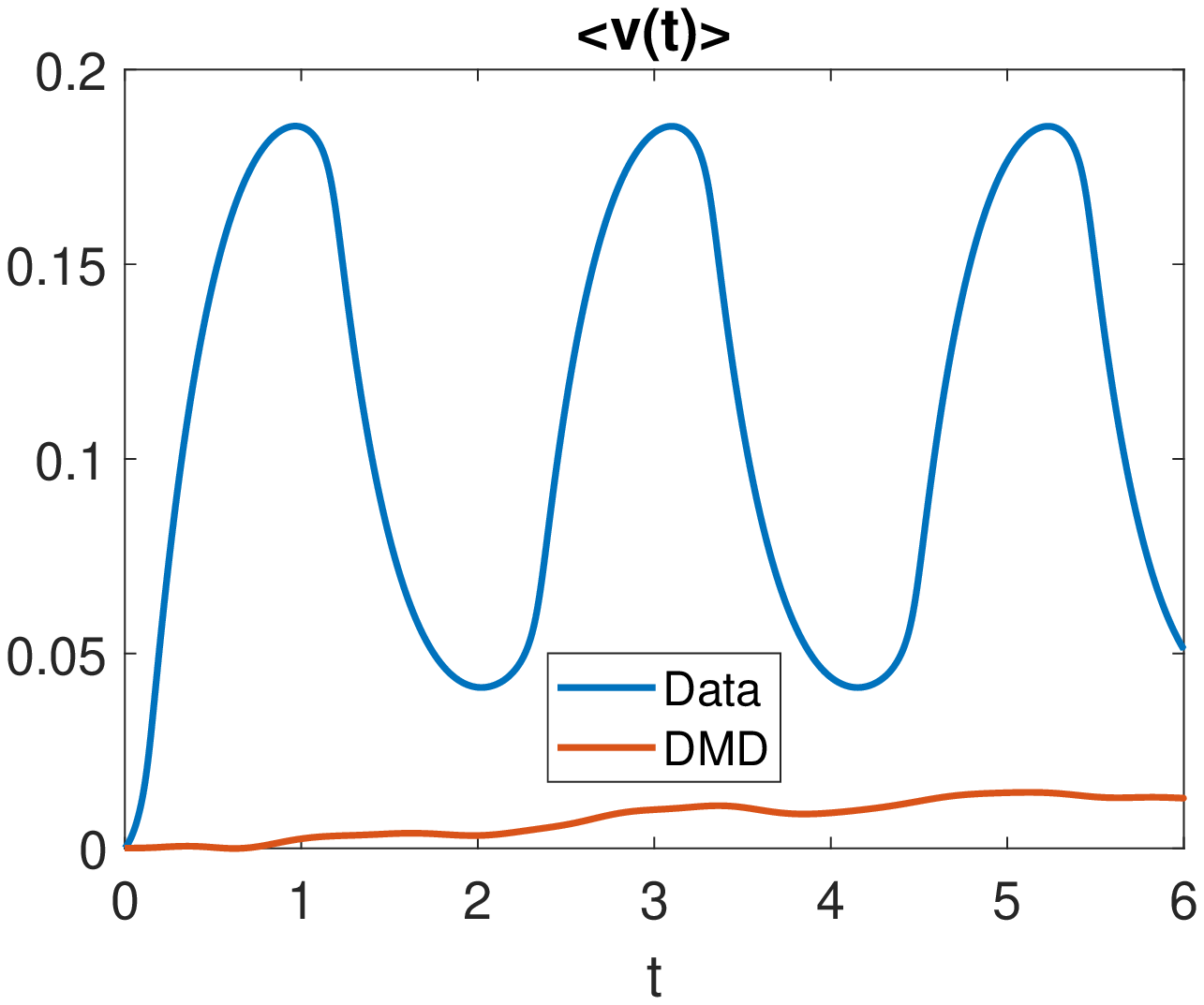}
\captionsetup{justification=justified}
\caption{FHN model, relaxation oscillations. Relative error \eqref{ef_data} of the DMD reconstruction $\widetilde S$ with respect to the dataset $S$ (left). Comparison of the spatial mean for the variables $u$ (center) and $v$ (right) for DMD of rank $r = 28$ corresponding to the minimum of the error. }
\label{fhn_err}
\end{figure}

The one dimensional FitzHugh-Nagumo (FHN) model describes the activation and deactivation dynamics of a spiking neuron and it is a simplified version of the more famous Hodgkin-Huxley model \cite{CS10}. The nonlinear reaction terms in \eqref{RDPDE} are given by
\begin{equation}
\label{par_fhn}
\begin{cases}
f(u,v) = \displaystyle \frac{u(u-0.1)(1-u)}{d_u} -\frac{v}{d_u} +\frac{c}{d_u}\\
g(u,v) = b u - \gamma v + c
\end{cases}
\end{equation}
to build the dataset $S$, here we solve the FHN system on the 1D domain $\Omega = [0,1]$, for $t \in [0,T], T = 6$, and $d_u = 0.015, d_v = 0, b = 0.5, \gamma = 2, c = 0.05$.
The initial and boundary conditions are given by
\begin{equation}\label{par:modCS}
\begin{aligned}
& u_0(x) = 0, \quad v_0(x) = 0, \quad x \in \Omega, \\
& u_x(0,t) = b_u(t)= -(5 \cdot 10^4 t^3 e^{(-15t)}), \quad u_x(1,t) = b_v(t)= 0, \quad t \in [0,T]. 
\end{aligned}
\end{equation}
For the spatial meshgrid we consider $n = 1024$ points. We integrate in time by applying the IMEX Euler scheme with $h_t = 10^{-3}$, thus we have $m + 1 = 6000$. 
The parameters in \eqref{par_fhn} and \eqref{par:modCS} are taken from \cite{CS10}.\\
%
In the left panel of Figure \ref{fhn_err}, we show the relative error $\mathcal{E}(\widetilde{S},r),$ for $r = 1, \dots, R$, where $R=51$ is the rank of the snapshot matrix $S \in \RR^{2048 \times 6000}$. 
This error is very high and erratic, indicating that the DMD reconstruction for both variables is completely wrong. Even worse, the error increases when we consider higher values for the rank, which is something not expected a priori. The minimum value is reached for $r=28$, i.e. $\mathcal{E}(\widetilde{S},28) = 0.9618$.
For this $r$ value, in the middle and right panels of Figure \ref{fhn_err} we compare the spatial mean of the data \eqref{mean} with that of the DMD reconstruction, for both $u$ and $v$. It is clear that the DMD reconstruction does not capture neither the periodic dynamics nor the amplitude of the relaxation oscillations of the FHN model.

\begin{figure}[htbp]
\centering
\includegraphics[scale=0.4]{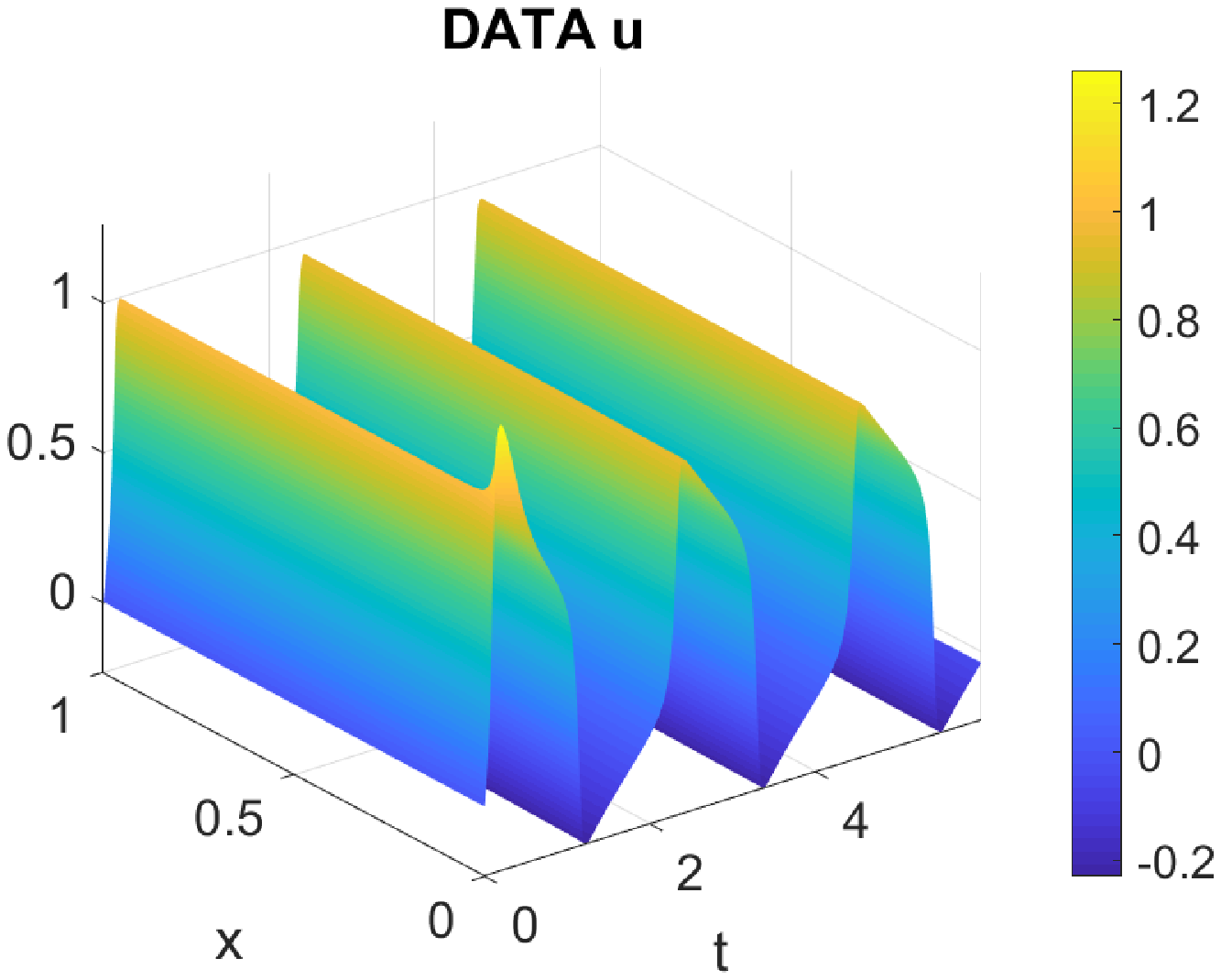}
\includegraphics[scale=0.4]{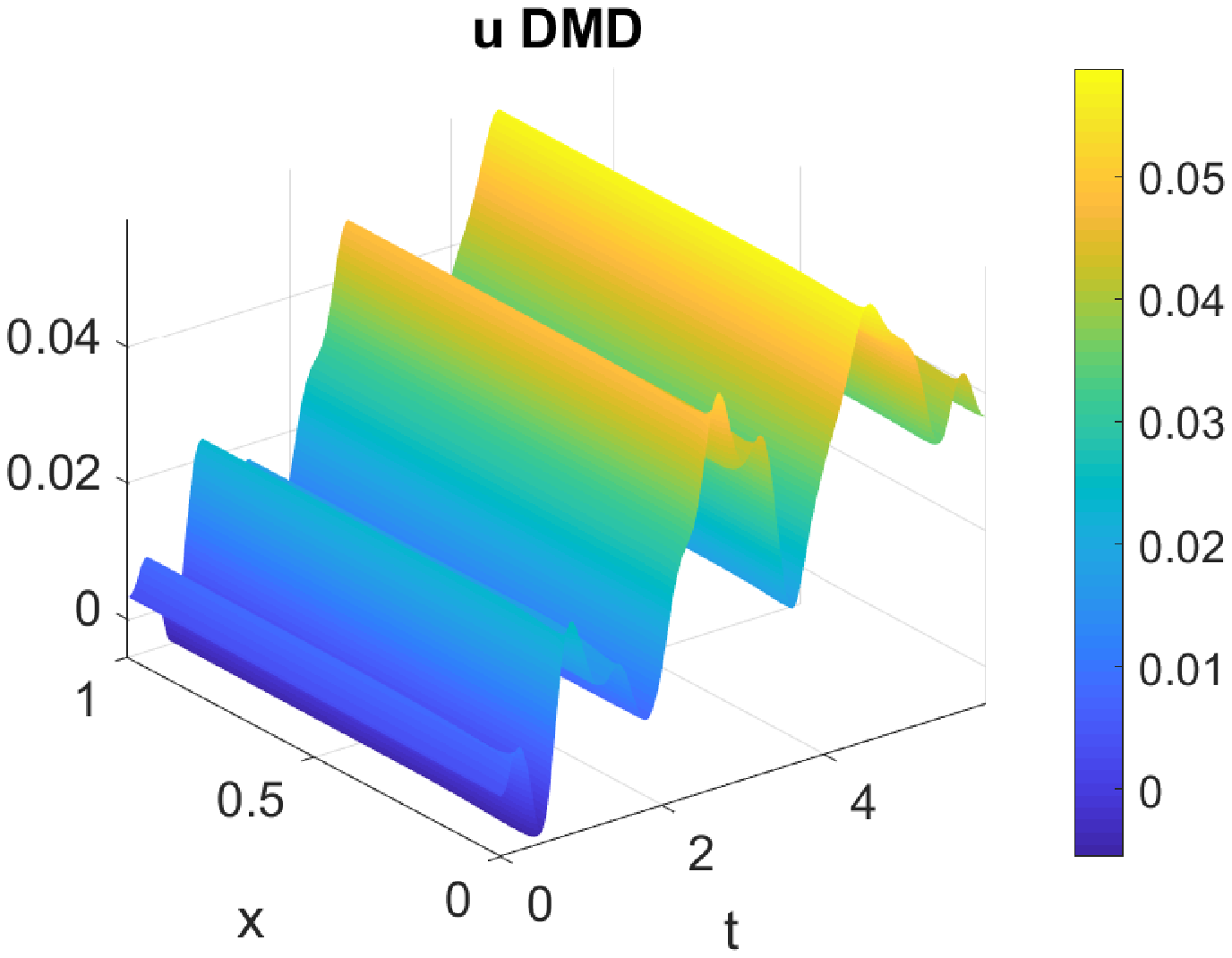}\\
\includegraphics[scale=0.4]{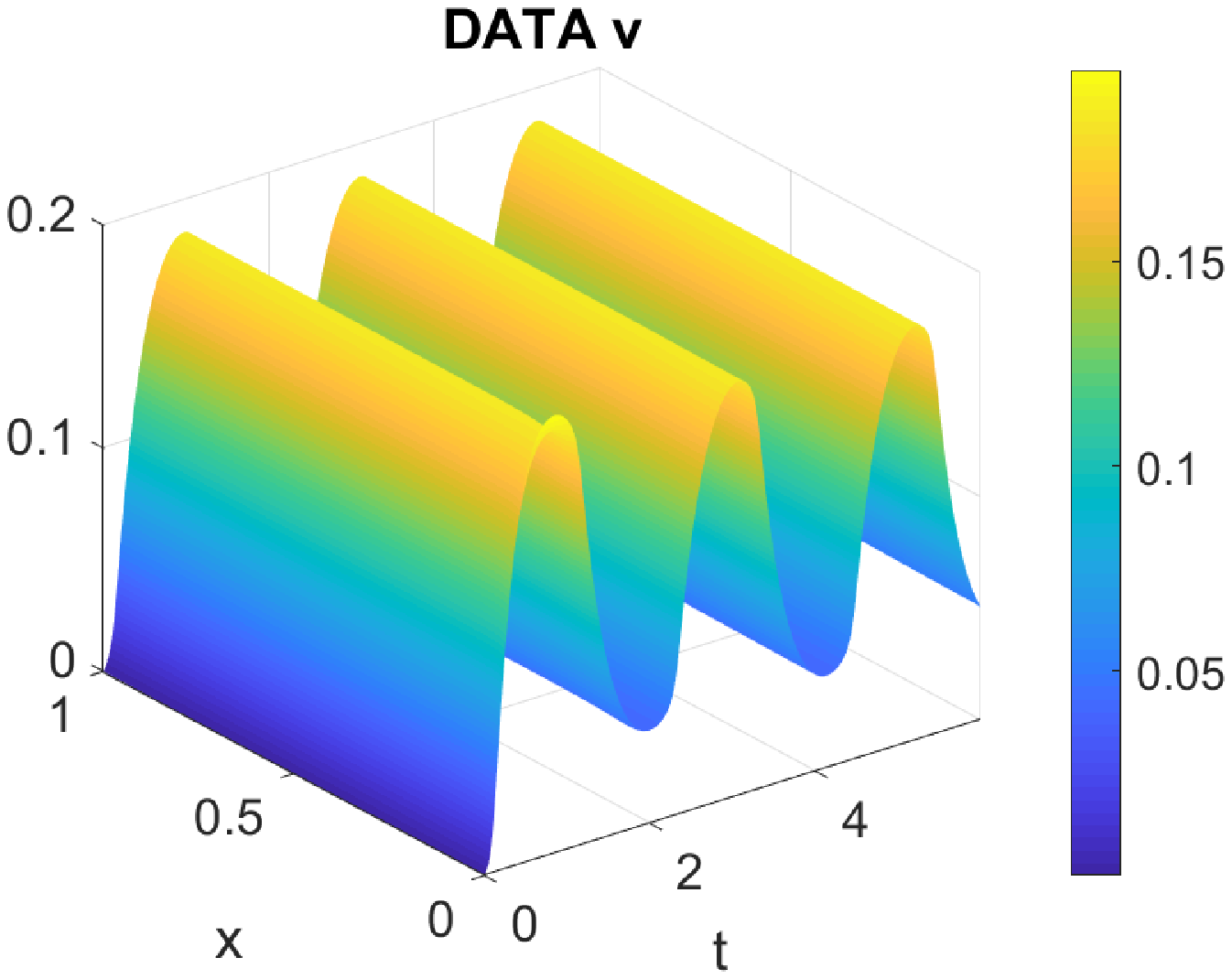}
\includegraphics[scale=0.4]{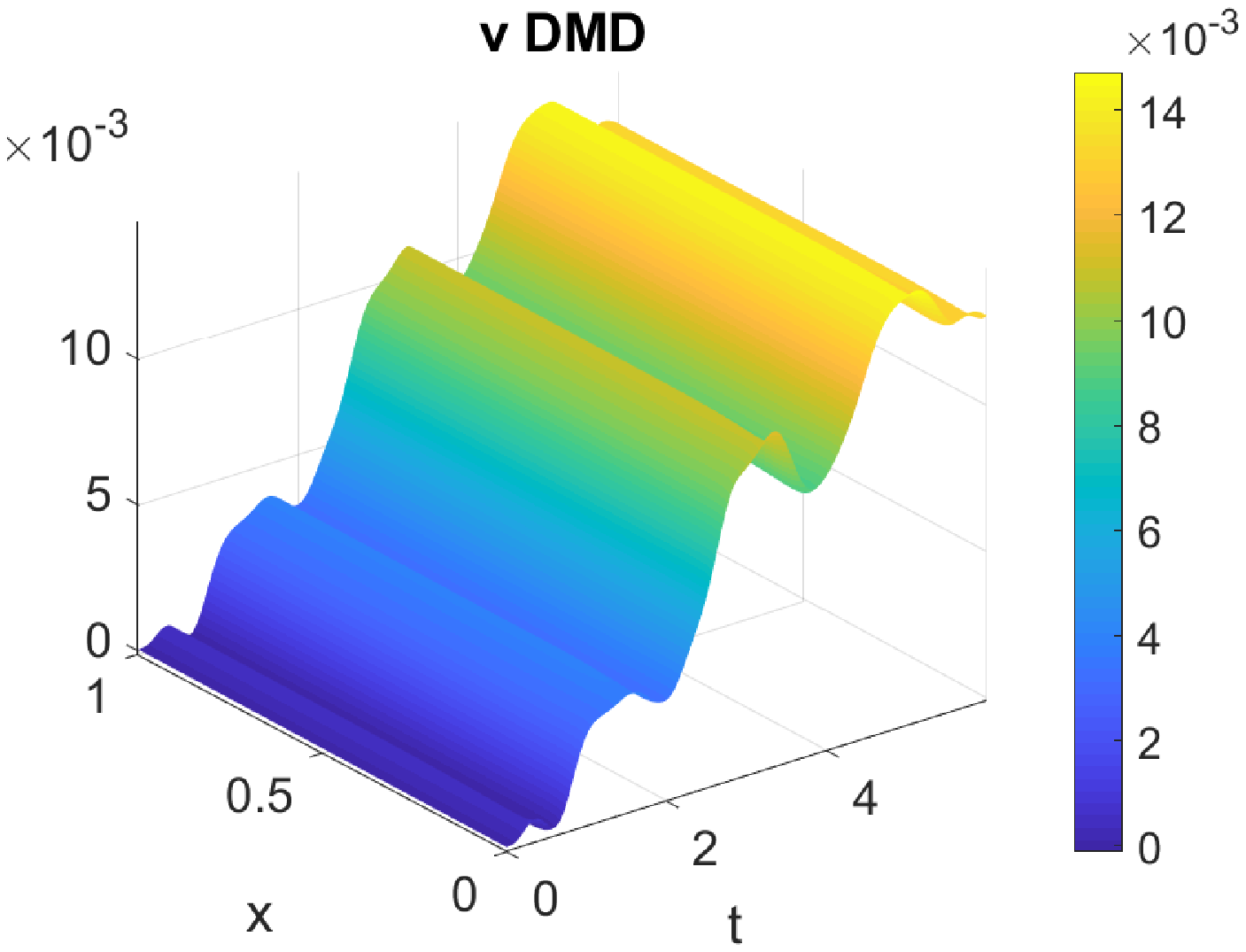}
\captionsetup{justification=justified}
\caption{FHN model, relaxation oscillations. Full model solutions $u$ and $v$ (left panels) and corresponding DMD reconstructions (right panels) for $r = 28$.}
\label{fhn_dmd_classic}
\end{figure}

\begin{figure}[htbp]
\centering
\includegraphics[scale=0.45]{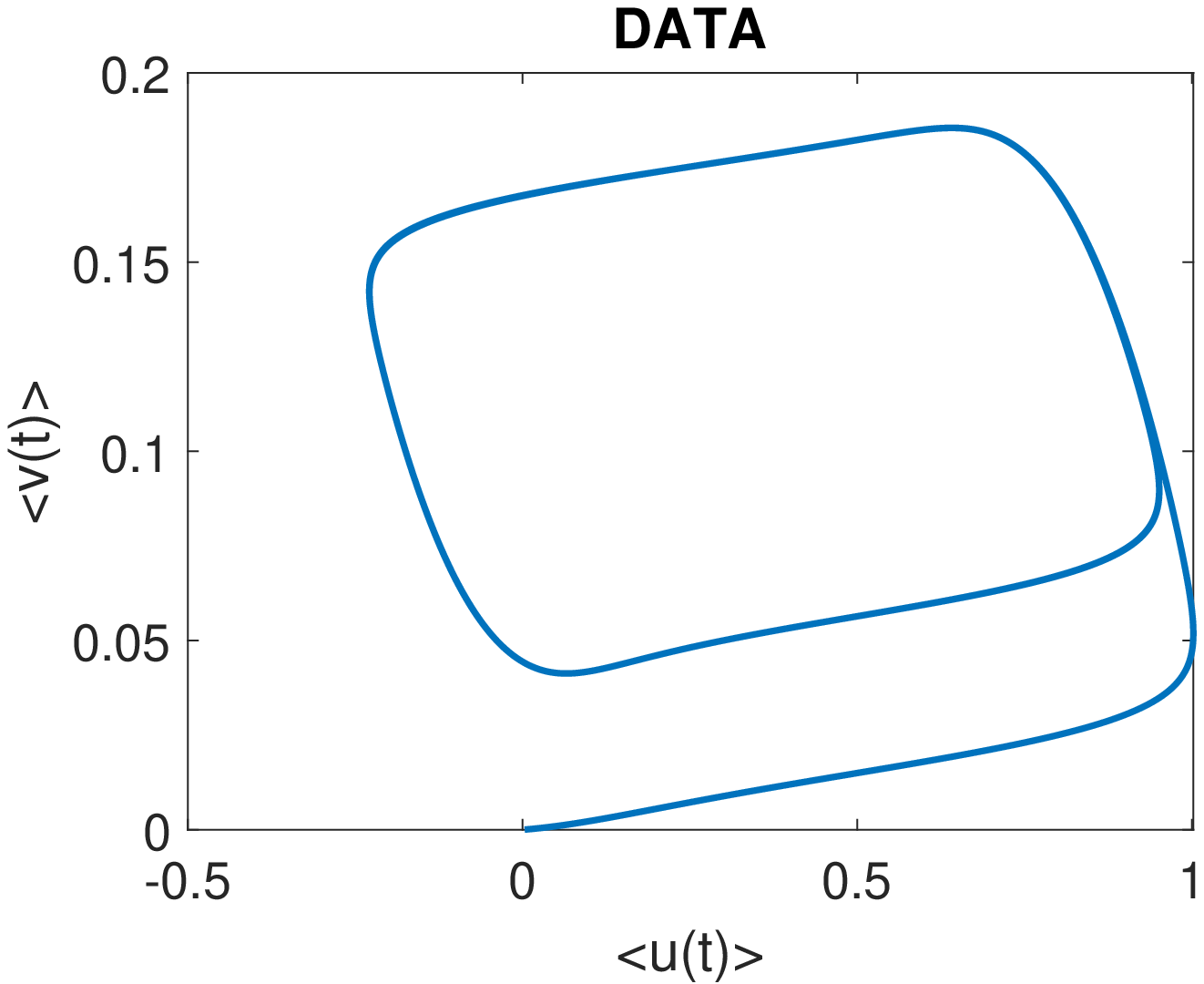}
\includegraphics[scale=0.45]{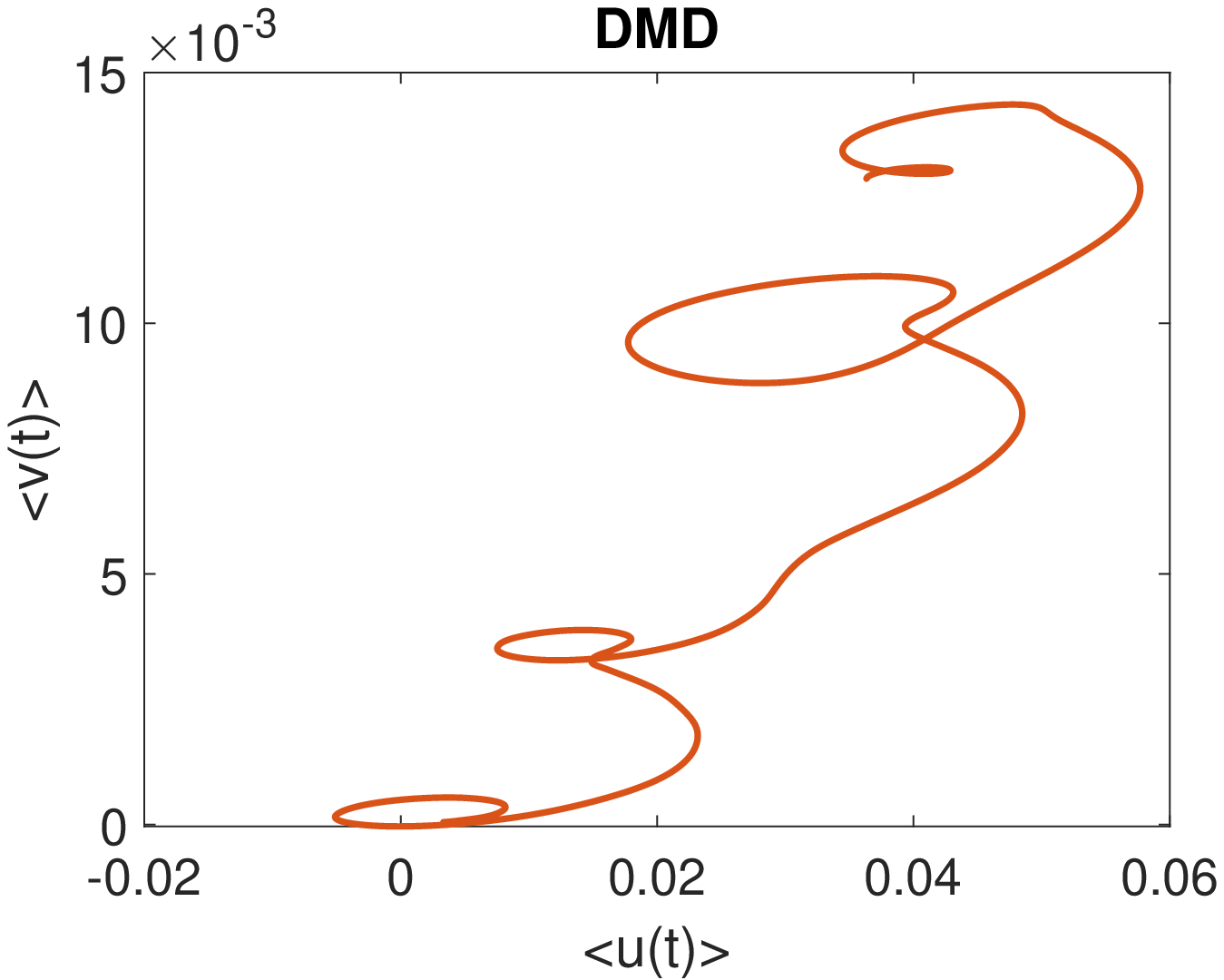}
\captionsetup{justification=justified}
\caption{FHN model, relaxation oscillations. Limit cycle corresponding to the dataset $S$ in the phase plane $(\langle u \rangle , \langle v \rangle )$(left) and its DMD approximation for $r = 28$ (right). }
\label{fhn_phase_classic}
\end{figure}

In Figure \ref{fhn_dmd_classic}, right panel, we show the DMD reconstructions of $u$ and $v$ in space and time on the domain $\Omega \times [0,T]$. By comparing the dataset $S$ reported in the left panels becomes still more evident how the DMD method fails. The effect of this DMD failure on the limit cycle reconstruction is shown in Figure \ref{fhn_phase_classic}. 

\subsubsection{$\lambda$-$\omega$ system: spiral waves}
\label{ex_lo}
Here we consider the $\lambda$-$\omega$ system from \cite{MV06, BG2006}, with nonlinear kinetics \eqref{RDPDE} given by
\begin{equation}
\label{lamda-omega}
\begin{aligned}
f(u,v) &= \rho (\lambda(u,v) u - \omega(u,v) v), \\
g(u,v) &= \rho (\omega(u,v) u + \lambda(u,v) v),
\end{aligned}
\end{equation}
where $\lambda = 1-(u^2 + v^2)$ and $\omega = -\beta(u^2+v^2), \beta >0$. On the 2D spatial domain $\Omega = [0,L]\times [0,L]$, with $L = 130$, we choose the parameter values and initial conditions from \cite{BG2006} for which spiral waves solutions arise:
$$ d_u = 4, \ d_v = 4, \ \rho = 10, \ \beta = 1,$$
$$u_0(x,y) = \frac{1}{10}(x-\frac{L}{2}), \ v_0(x,y) =\frac{1}{10} (-\frac{y}{2}+\frac{L}{4})$$
and homogeneous Neumann boundary conditions $b_u(t) \equiv 0 \equiv b_v(t)$ in \eqref{RDPDE}.
The domain $\Omega$ is discretized by $n_x = n_y = 99$ interior points, such that the total number of mesh points is $n = n_x n_y = 9801$. We integrate in time with timestep $h_t = 10^{-3}$ until $T = 50$, by using the IMEX Euler scheme in matrix-oriented form (see \cite{DSS20}). We save the snapshots every four time steps ($\kappa=4$), such that the considered dataset $S$ has dimension $2n \times (m+1)$ and $m+1 = 12500$.

\begin{figure}[htbp]
\centering
\includegraphics[scale=0.25]{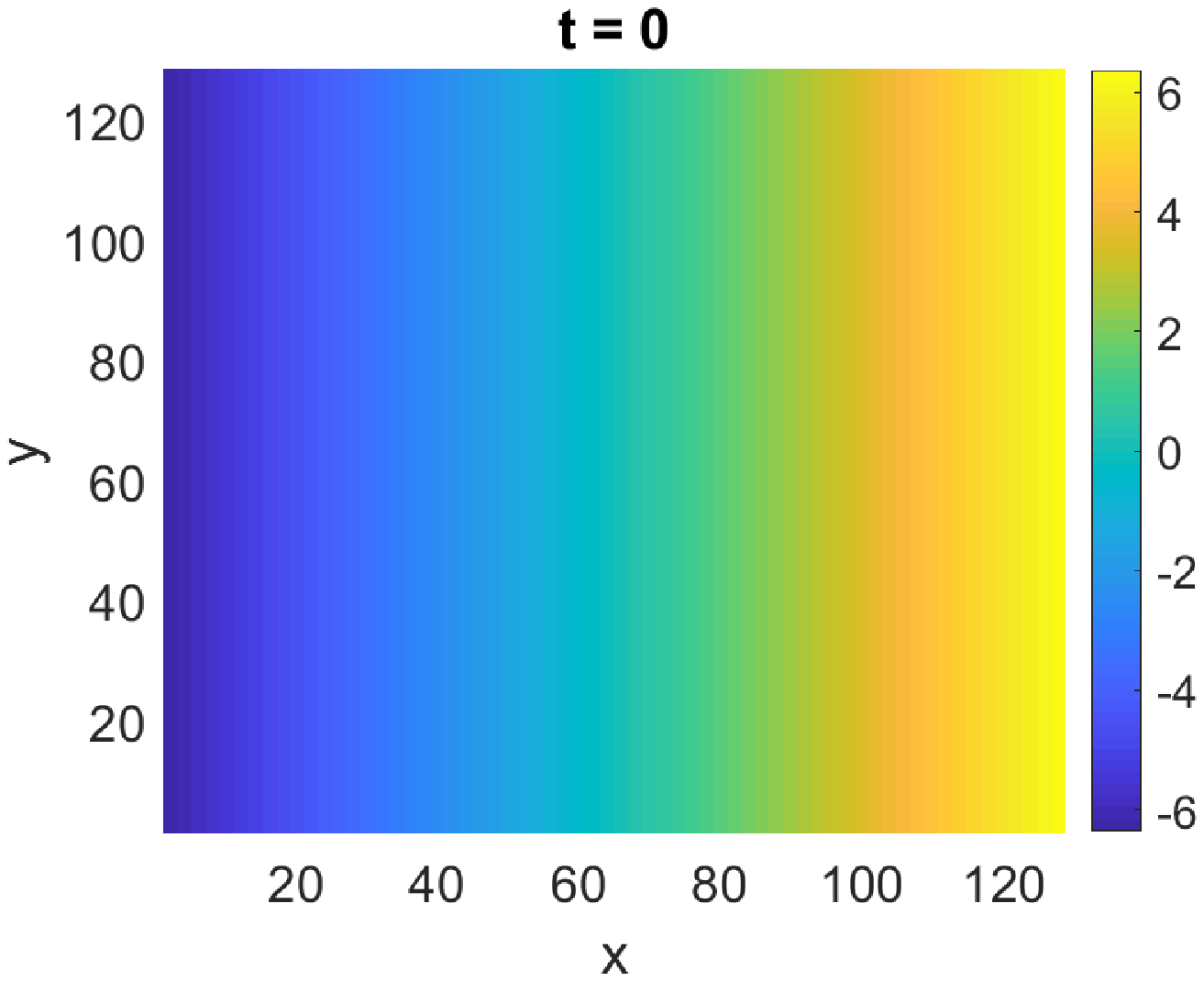}
\includegraphics[scale=0.25]{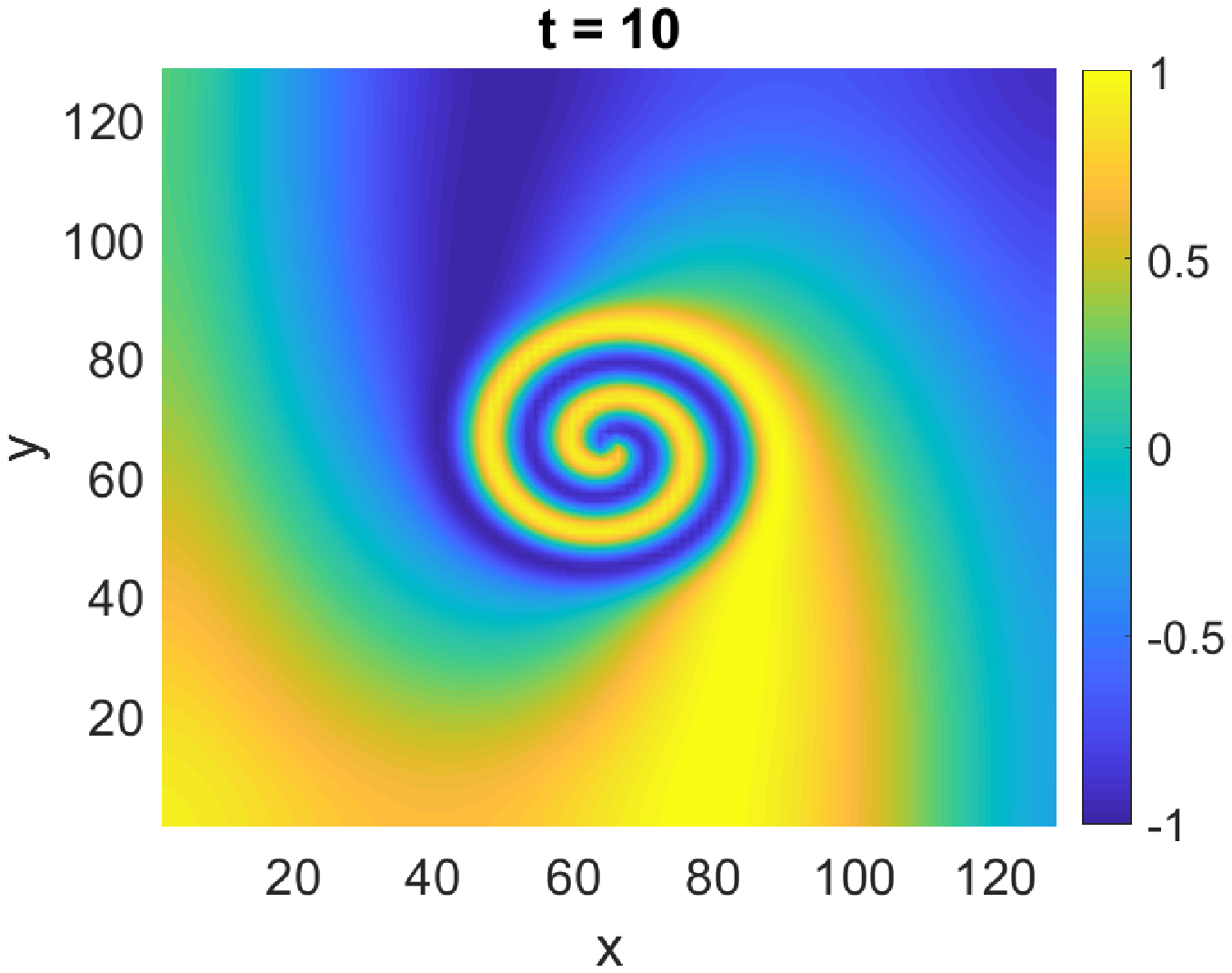}
\includegraphics[scale=0.25]{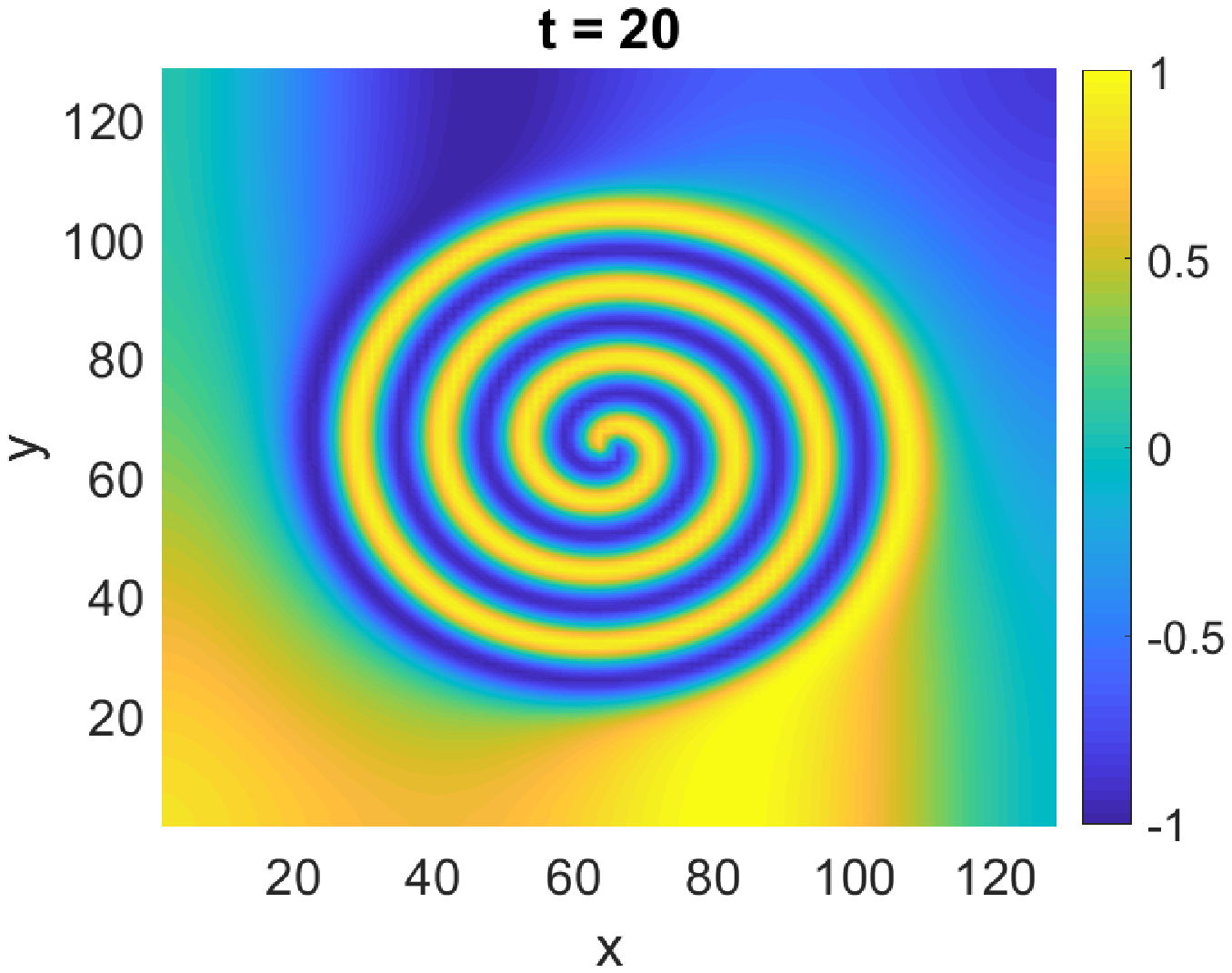}\\
\includegraphics[scale=0.25]{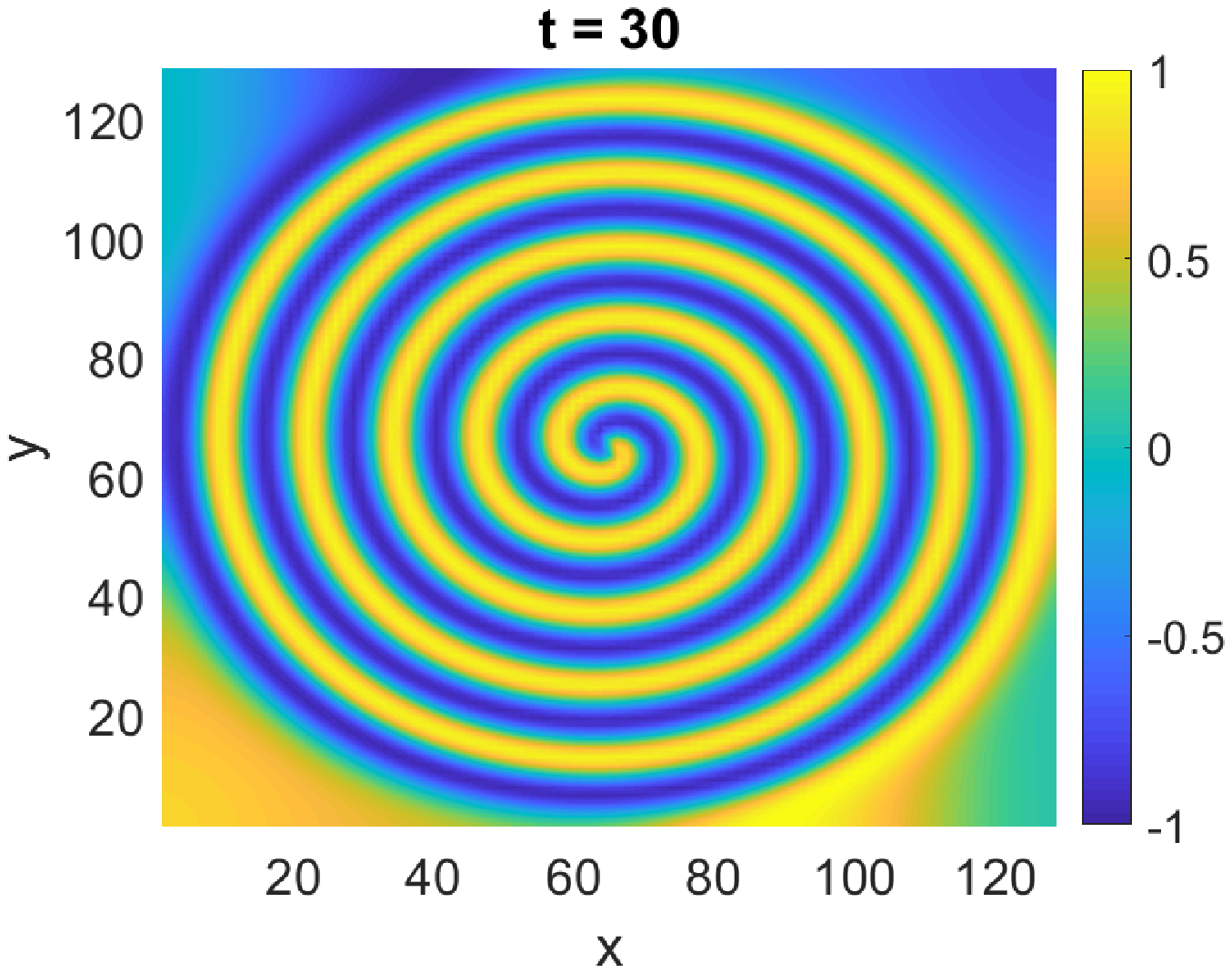}
\includegraphics[scale=0.25]{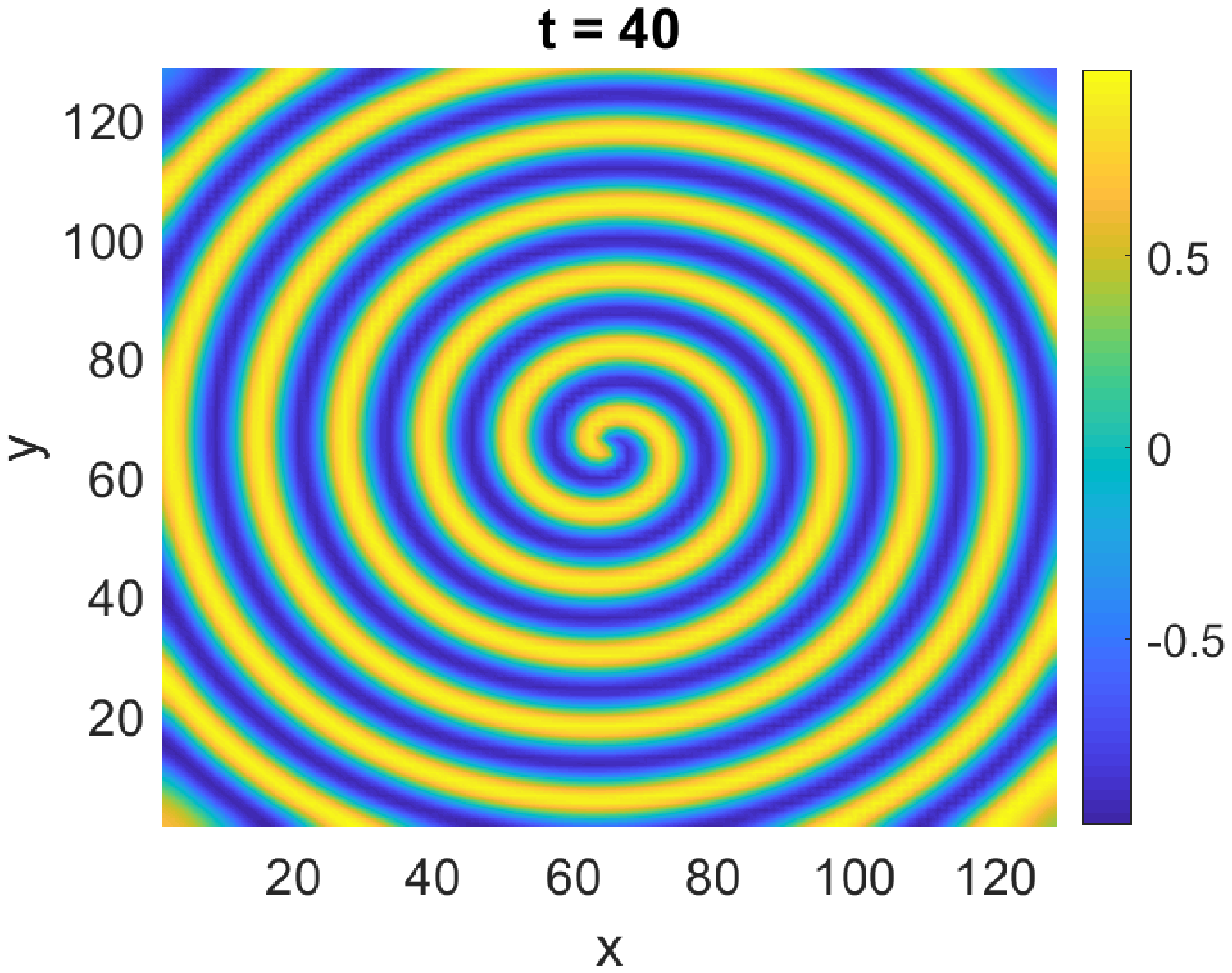}
\includegraphics[scale=0.25]{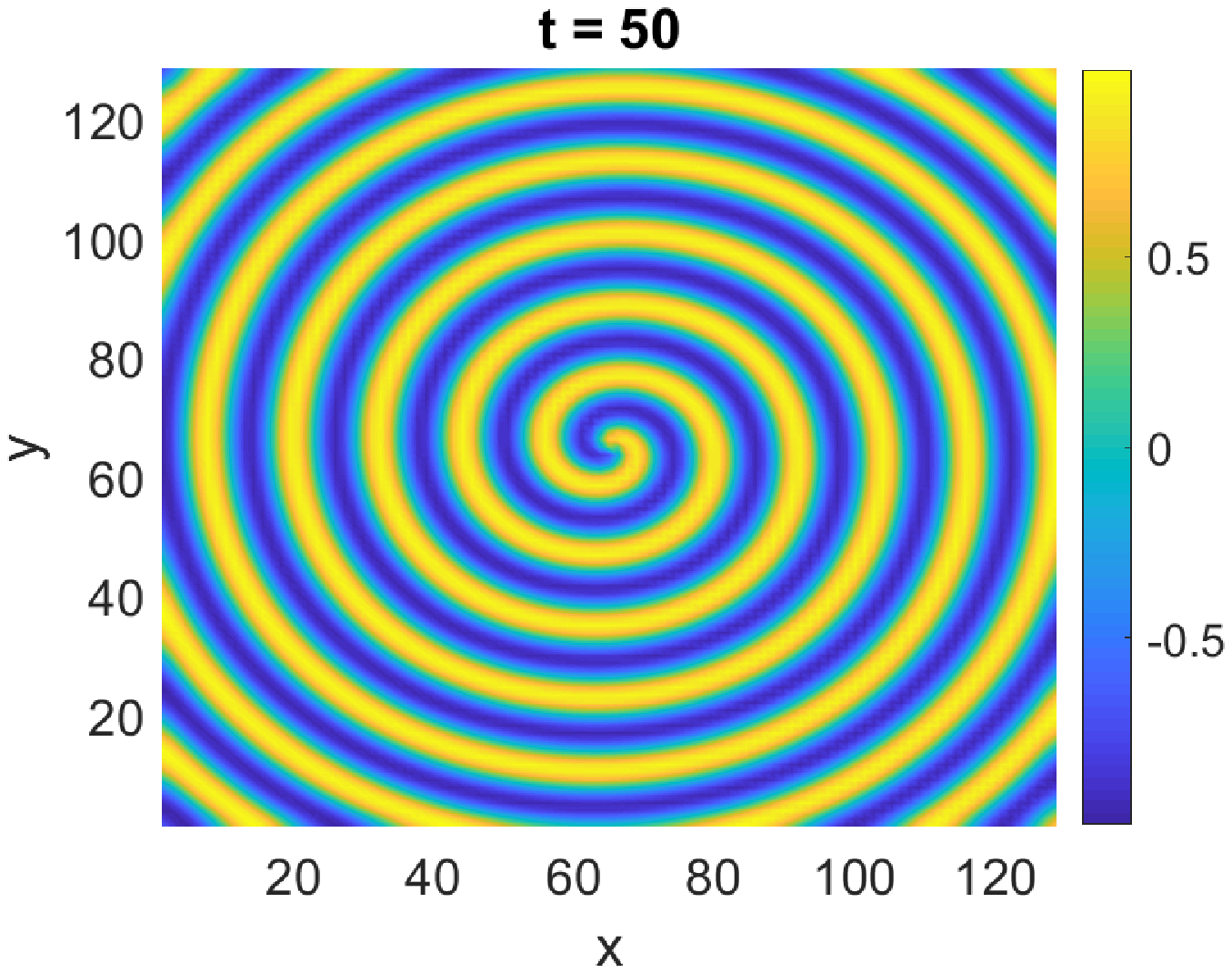}
\captionsetup{justification=justified}
\caption{$\lambda-\omega$ system. Spiral wave generated for the variable $u$ at different times in $[0, T]$ see more details in the main text.}
\label{lo_story}
\end{figure}

\begin{figure}[htbp]
\centering
 \includegraphics[scale=0.4]{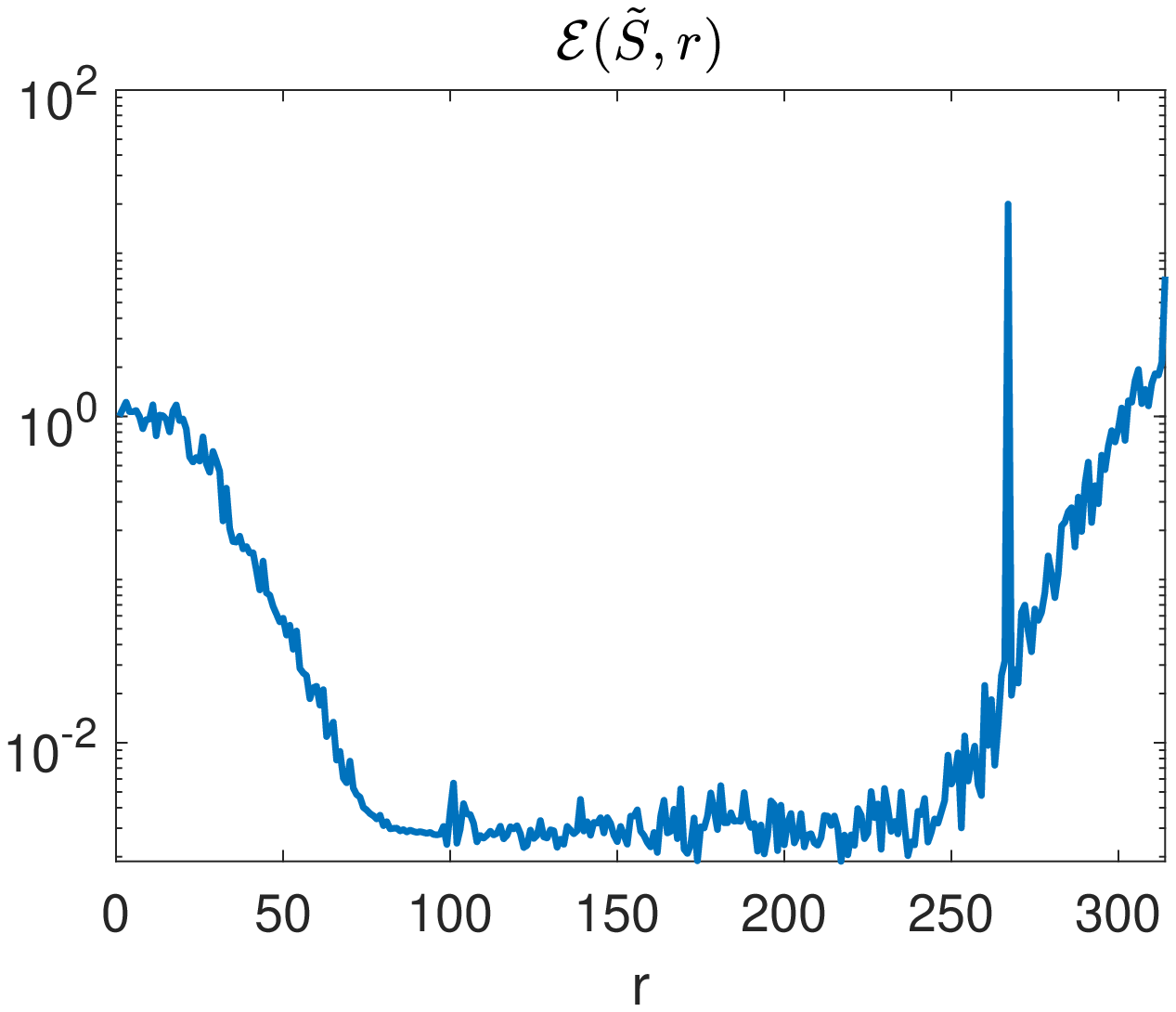}
\includegraphics[scale=0.4]{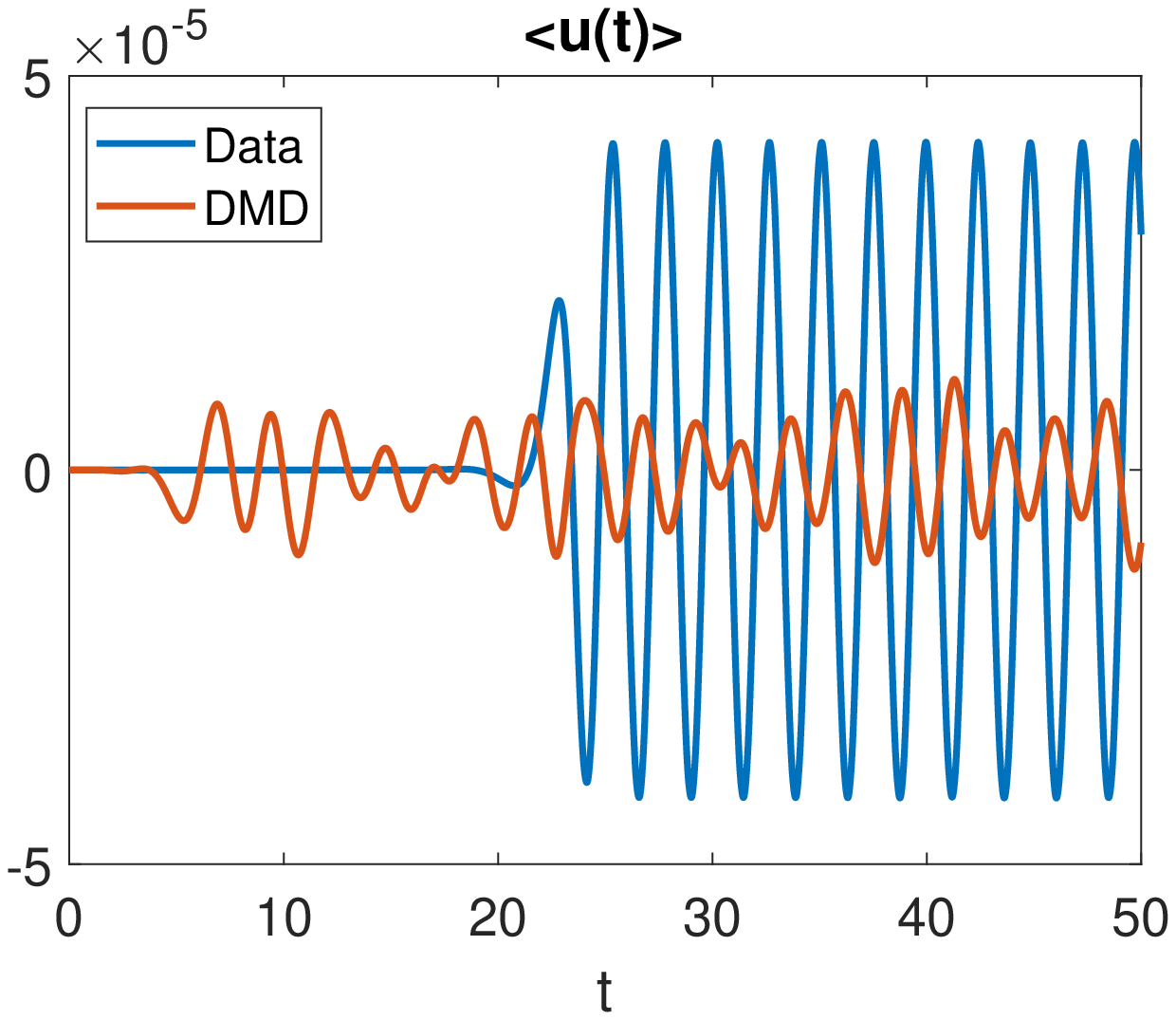}
 \includegraphics[scale=0.4]{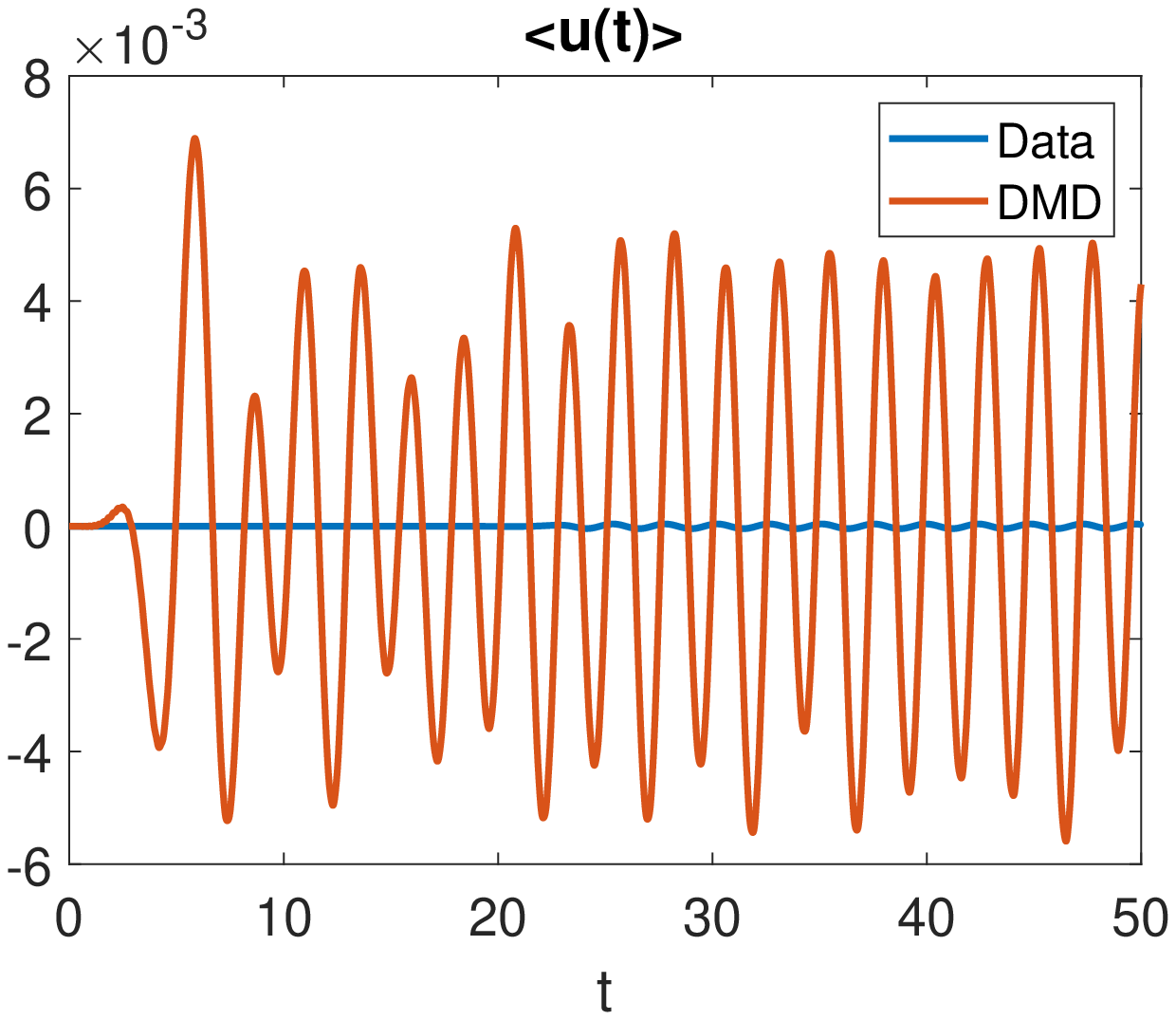}
\captionsetup{justification=justified}
\caption{$\lambda$-$\omega$ system, spiral wave. DMD relative error in \eqref{ef_data} for the dataset $S$. Comparison of the spatial means for $u$: DMD is applied with $r = 217$ (center plot) and $r = 300$ (right plot).}
\label{lo_mean_classic}
\end{figure}

\begin{figure}[htbp]
\centering
\includegraphics[scale=0.4]{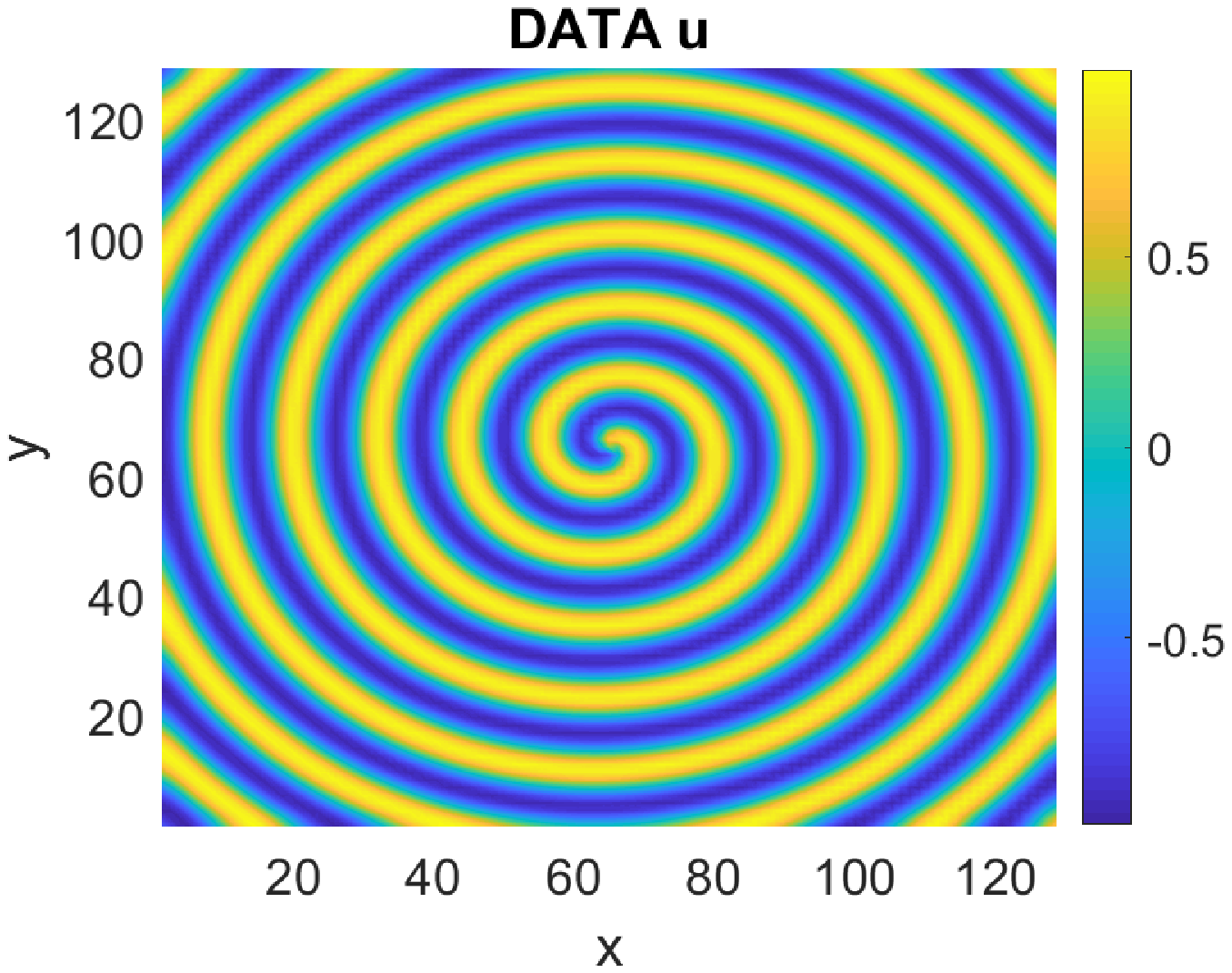}
\includegraphics[scale=0.4]{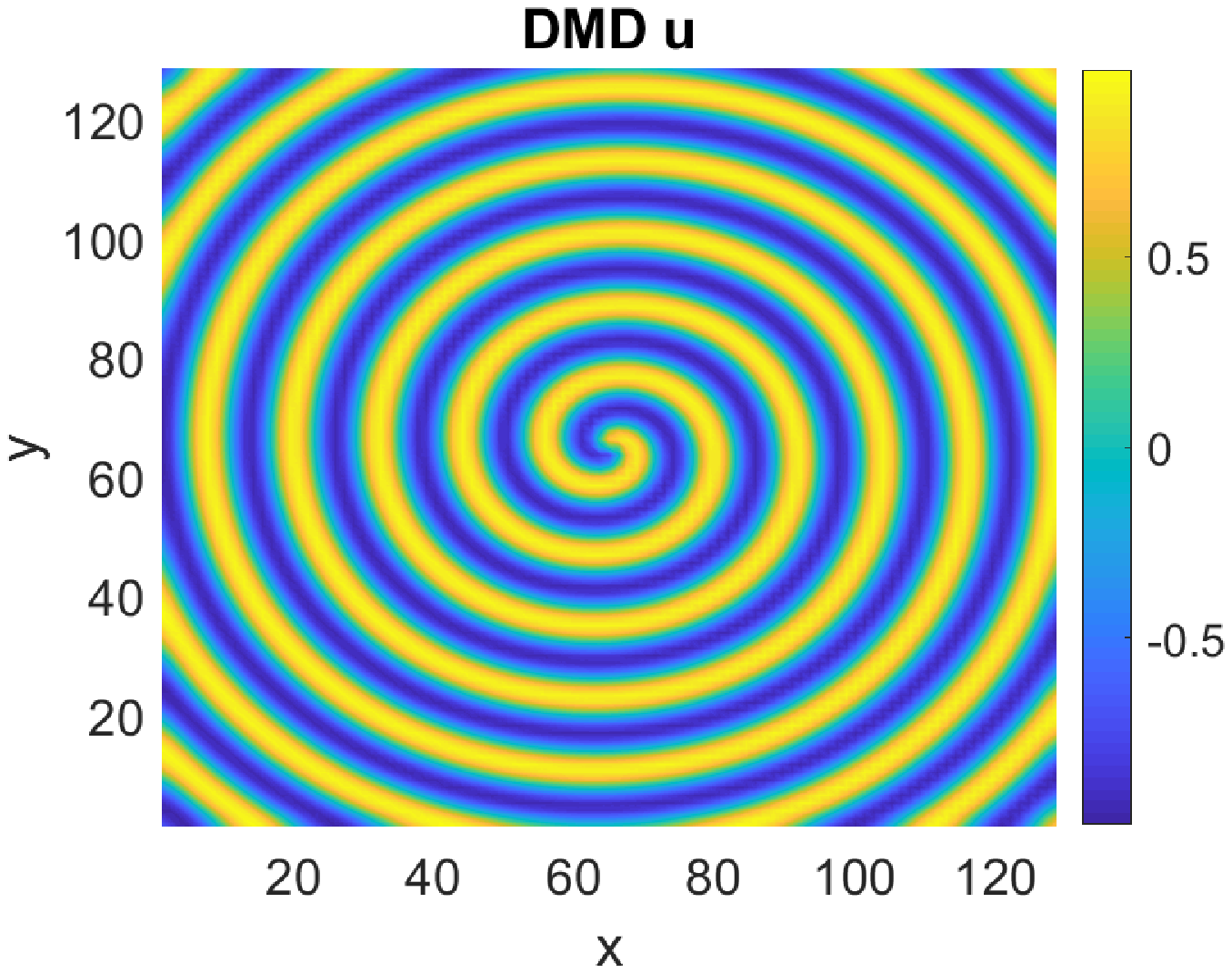}
\includegraphics[scale=0.4]{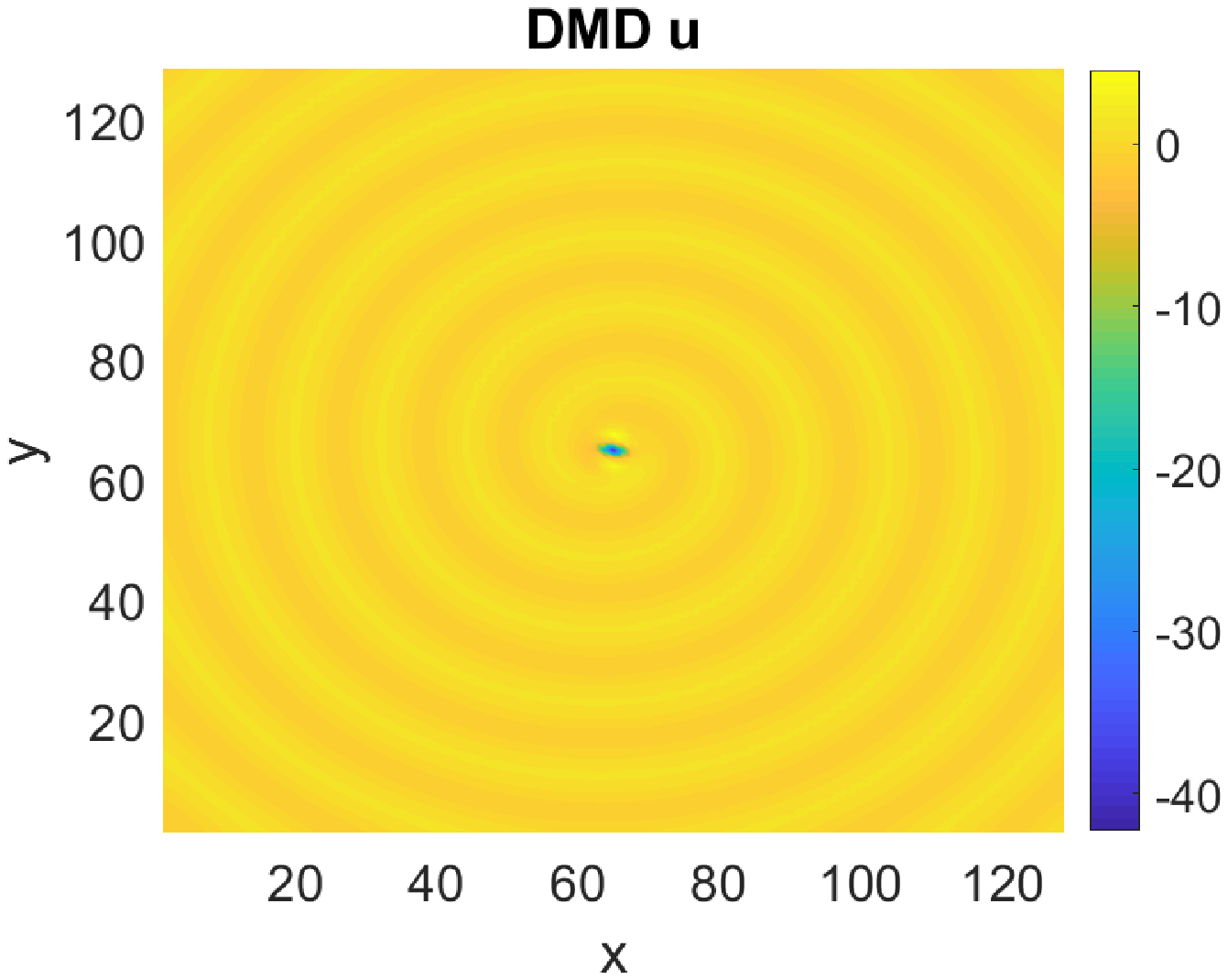}
\captionsetup{justification=justified}
\caption{$\lambda$-$\omega$ system, spiral wave. Data for the variable $u$ at the final time $T=50$ (left) and its DMD reconstructions with $r = 217$ (center plot) and $r = 300$ (right plot).}
\label{lo_dmd_classic}
\end{figure}

\begin{figure}[tbp]
\centering
\includegraphics[scale=0.4]{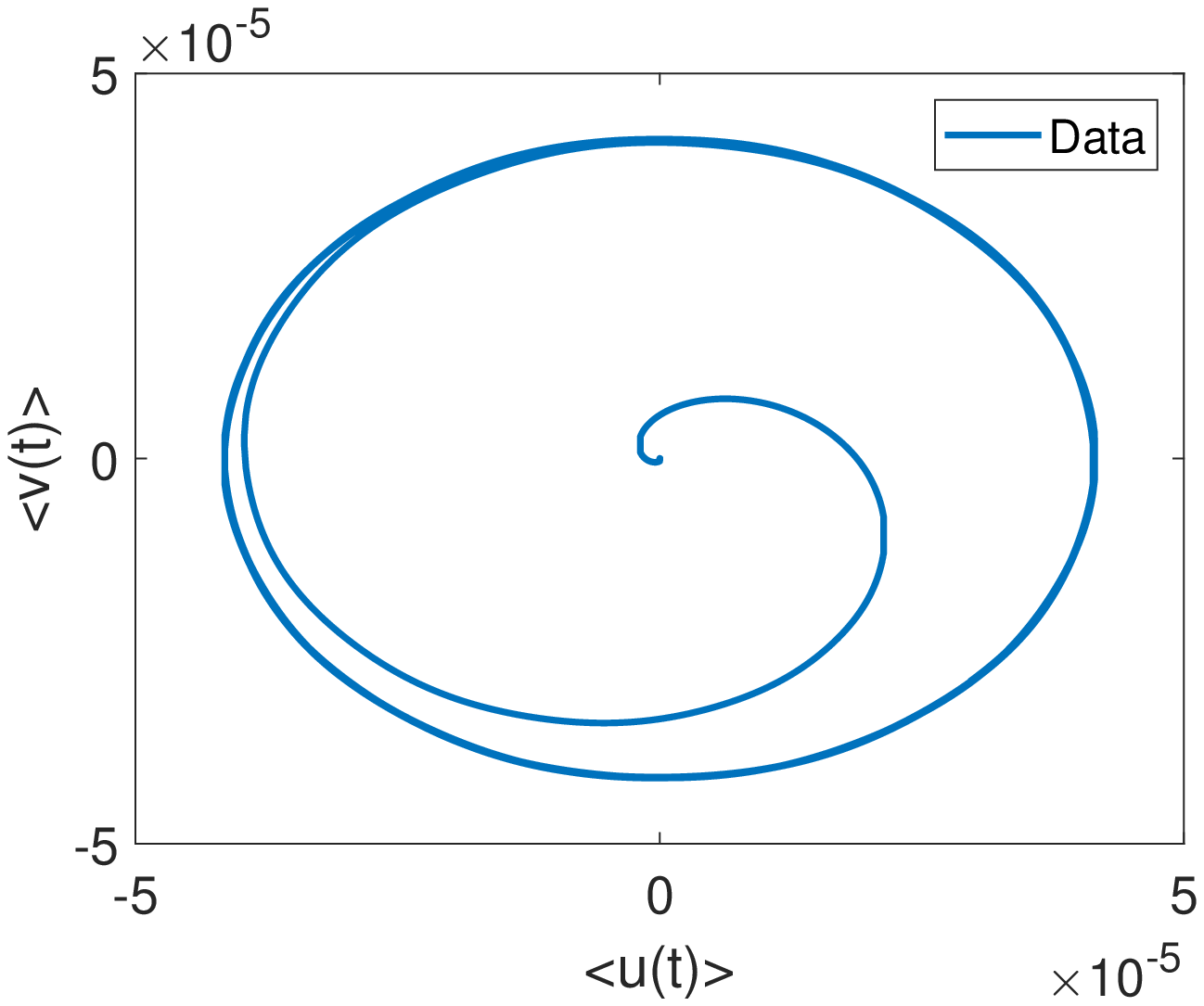}
\includegraphics[scale=0.4]{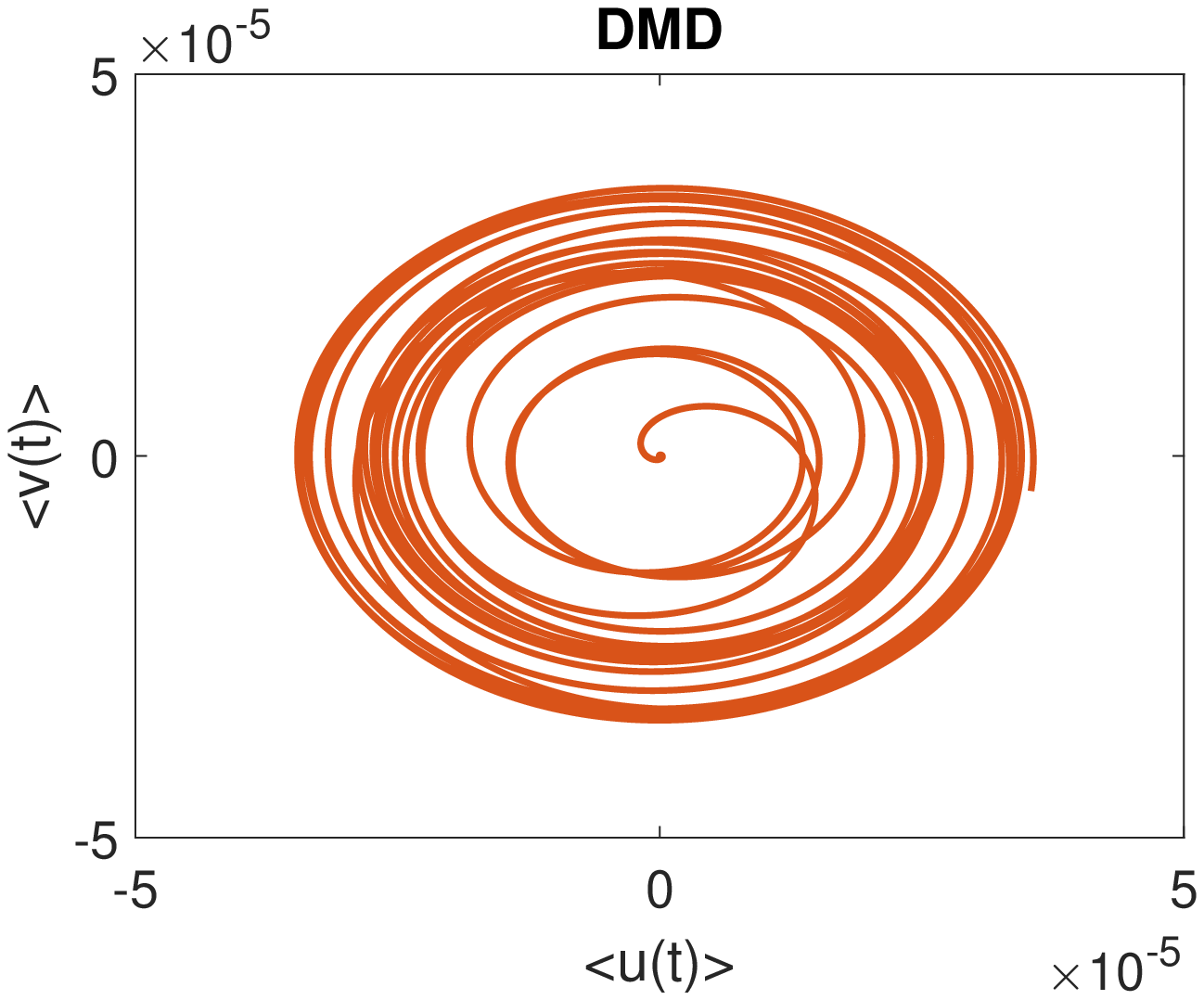}
\includegraphics[scale=0.4]{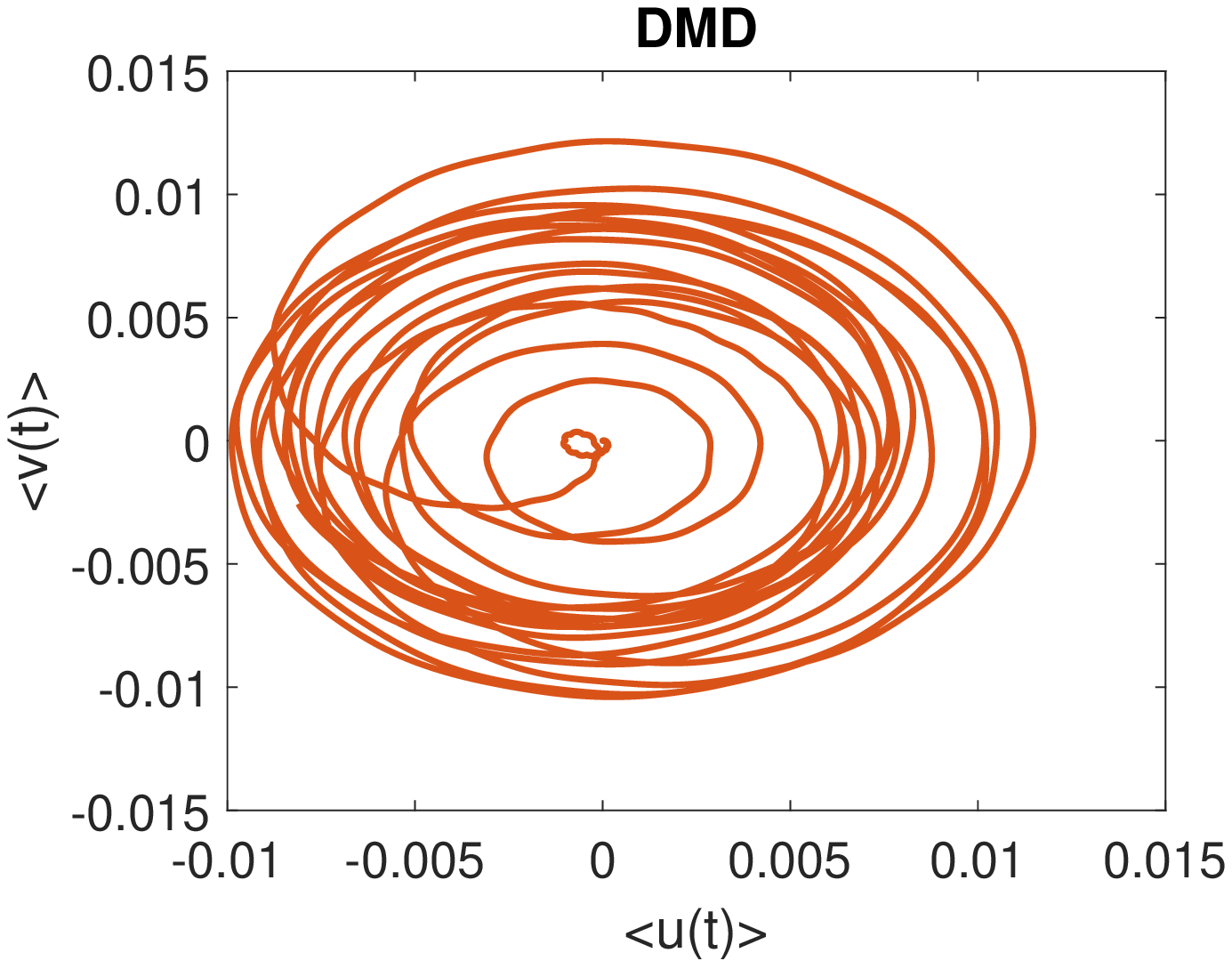}
\captionsetup{justification=justified}
\caption{$\lambda$-$\omega$ system, spiral wave. Comparison of the limit cycles in the phase plane $(\langle u \rangle, \langle v \rangle)$: data (left) and DMD with $r = 217$ (center) and $r = 300$ (right).}
\label{lo_dmd_phase}
\end{figure}

Departing from the step values of the initial data, for both $u$ and $v$, in a transient regime the numerical solution starts to form an archimedean spiral wave with ``core " (fixed point) in the center of the domain $(x_c,y_c)=(65,65)$ which arms oscillate in space and time until at a certain time, say $\bar{t}$, when the entire $\Omega$ is covered (as shown in Figure \ref{lo_story}). Thereafter, for $t \geq \bar{t}$ a new time regime arises where the spiral continues indefinitely in its oscillating dynamics such that in the phase plane we can say that a limit cycle is attained by the spatial means $(\langle u \rangle, \langle v \rangle)$. We show the dynamics of $\langle u(t) \rangle$ in Figure \ref{lo_mean_classic} (center), where the two time regimes in $[0,\bar{t}], [\bar{t}, T]$, with $\bar{t} \approx 25$ are evident, then the $u$- snapshot at the final time $T$ in Figure \ref{lo_dmd_classic} (left) and the corresponding limit cycle in Figure \ref{lo_dmd_phase}(left).

In Figure \ref{lo_mean_classic} (left) we show the relative error $\mathcal{E}(\widetilde{S},r),$ obtained by the DMD for $r = 1, \dots, R$, where the rank of the dataset $S$ is $R=314$. After an initial decreasing trend, for $r>\approx100$ an erratic plateau around  $3$e-3 can be observed where the minimum value $\mathcal{E}(\widetilde{S},r)=0.0019$ is reached for $r=217$. Moreover, this low accuracy gets worse for $r \geq 250$ when the error dramatically increases due to ill-conditioning of the matrix $\widetilde A$. We stress this DMD drawback because usually better results are expected by increasing the value of $r$ and this clearly does not happen here. 
To support this conclusion, in Figure \ref{lo_mean_classic}, we also compare the spatial mean dynamics for the variable $u$ obtained by DMD for $r= 217$ (middle plot) and $r = 300$ (right plot) when $\mathcal{E}(\widetilde{S},300)=0.8471$.

In both cases, we note that: i) the approximation of the two distinct time regimes is missed, ii) a large difference in the wave amplitude is present; iii) in the best case $r=217$ (center plot), the frequency of the oscillations is preserved, but they are in phase opposition as shown in the middle panel of Figure \ref{lo_mean_classic}.


In Figure \ref{lo_dmd_classic}, we show the full model solution (left panel) and the DMD reconstructions with rank $r = 217$ (middle plot) and $r = 300$ (right plot) at the final time $T = 50$. For $r = 300,$ DMD fails essentially in the core of the spiral, instead for $r = 217$, where the global Frobenius error is minimum, DMD seems to be in great agreement with the data, even though its time history is really different, as discussed above. 

To further support the above points (i)--(iii) concerning the defects of DMD in the time dynamics approximation, in Figure \ref{lo_dmd_phase} we report the DMD reconstructions for $r = 271$ and $r = 300$ (center and right plot, respectively) in the phase plane $(\langle u \rangle, \langle v \rangle)$. In fact, by comparing the attained limit cycles with that for the data (left plot) it is still more evident that both DMD reconstruction fails. 

To conclude, in this section we have shown two examples with oscillatory datasets where DMD fails its reconstruction.

\subsection{Example on the Turing instability}\label{sec:turing}

\begin{figure}[htbp]
\centering
\includegraphics[scale=0.45]{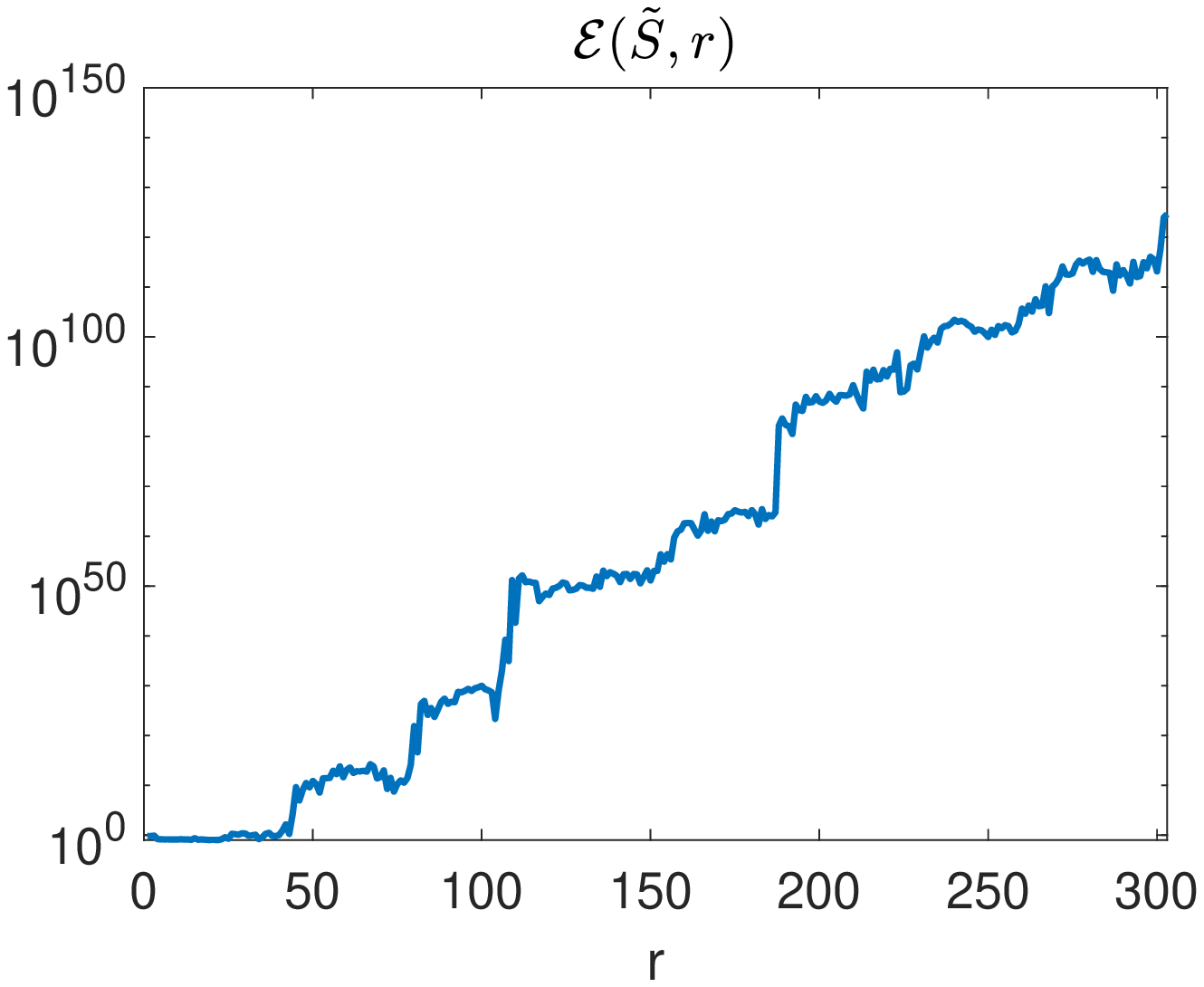}
\includegraphics[scale=0.45]{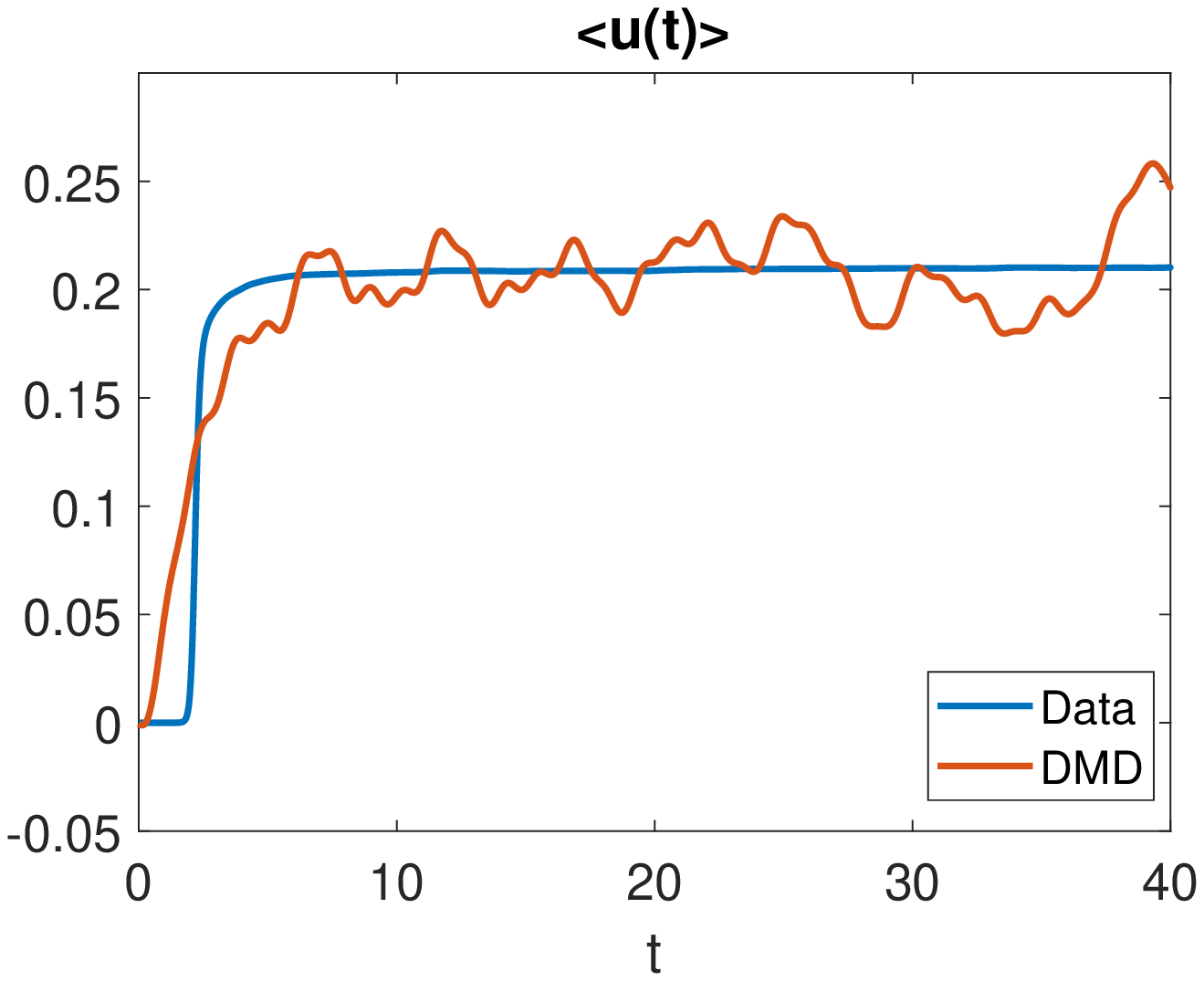}
\captionsetup{justification=justified}
\caption{DIB model, Turing instability. Left plot: relative error \eqref{ef_data}; right plot: the spatial mean for the variable $u$ in the dataset shows the two time regimes for the reactivity (until $t \approx 4$) and the stabilizing zone, whereas the DMD reconstruction for $r = 22$ destabilizes early and exhibits an oscillating behaviour for long times.} 
\label{dib_mean_classic}
\end{figure}

In this section, we deal with a RD model that exhibits the so-called Turing instability. In this case, the initial data are spatially random perturbations of the equilibrium of the model in absence of diffusion, say $(u_e,v_e)$. This peculiar dynamics presents essentially two time regimes: i) the {\it reactivity} zone where $(u_e,v_e)$ destabilizes from the initial conditions because of diffusion and ii) the {\it stabilizing} regime where the solution starts to be attracted towards a steady state spatially structured pattern, known as {\it Turing pattern} of different morphologies, see e.g. \cite{Murray03book, Neubert97, Neubert02}. The challenges for low rank techniques, like Proper Orthogonal Decomposition (POD) and Discrete Empirical Interpolation Method (DEIM), to reconstruct both regimes have been already highlighted in \cite{AMS23}. For the DMD performance an initial study is reported in \cite{BMS21} where the authors have considered an uncoupled approach, that is they have reconstructed independently the unknowns. 
This example focuses on the DIB morphochemical model that is an important realistic application for electrochemical phase formation modelling (\cite{DIB13, DIB15, SLB19}). The kinetics in \eqref{RDPDE} are given by
\begin{equation}
\label{DIB_kin}
\begin{aligned}
f(u,v) &= \rho \big{(}A_1(1-v)u - A_2u^3-B(v-\alpha) \big{)}, \\
g(u,v) &= \rho \big{(} C(1+k_2 u)(1-v)[1-\gamma(1-v)]-Dv(1+k_3u)(1+\gamma v) \big{)}.
\end{aligned}
\end{equation}
If $D = \frac{C(1-\alpha)(1-\gamma+\gamma \alpha)}{\alpha(1+\gamma \alpha)}$, there exists the homogeneous equilibrium $(u_e,v_e) = (0,\alpha)$ that can undergo Turing instability \cite{DIB13}. Here, we consider the parameter values taken from \cite{AMS23}:
$$A_1 = 10, \ A_2 = 1, \ \alpha = 0.5, \ B = 66, \ C = 3, \ \gamma = 0.2,$$ 
$$d_u = 1, \ d_v = 20, \ k_2 = 2.5, \ k_3 = 1.5, \ \rho = \frac{25}{4}.$$
The initial conditions are spatially random perturbation of the homogeneous equilibrium, given by
$u_0(x,y) = u_e + 10^{-5} {\tt rand}(x,y), \quad v_0(x,y) = v_e + 10^{-5} {\tt rand}(x,y)$.
We discretize the spatial domain $\Omega =[0,20] \times [0, 20]$ with $n_x = n_y = 100$ spatial meshpoints, such that $n = n_x n_y = 10000$ and, for stability reasons, we consider the timestep $h_t = 10^{-3}$ until the final time $T = 40$. We save the snapshots every four time steps, such that the dataset is $S \in \RR^{2n \times 10000}$.
In the left panel of Figure \ref{dib_mean_classic}, we show the DMD relative error $\mathcal{E}(\widetilde{S},r)$ for $r=1, \dots, R,$ for $R=303$ corresponding to the rank of $S$. The error dramatically increases for large values of $r$ and indeed  blows up for $r \geq 50$. Its minimum is obtained for $\mathcal{E}(\widetilde{S},22) = 0.1008$ (almost $10\%$). Nevertheless the time dynamics of the spatial mean reconstruction for $r=22$ exhibits an oscillating behaviour around the mean $\langle u(t) \rangle$ of the dataset, as shown in the right panel of Figure \ref{dib_mean_classic}. 

\begin{figure}[htbp]
\centering
\includegraphics[scale=0.4]{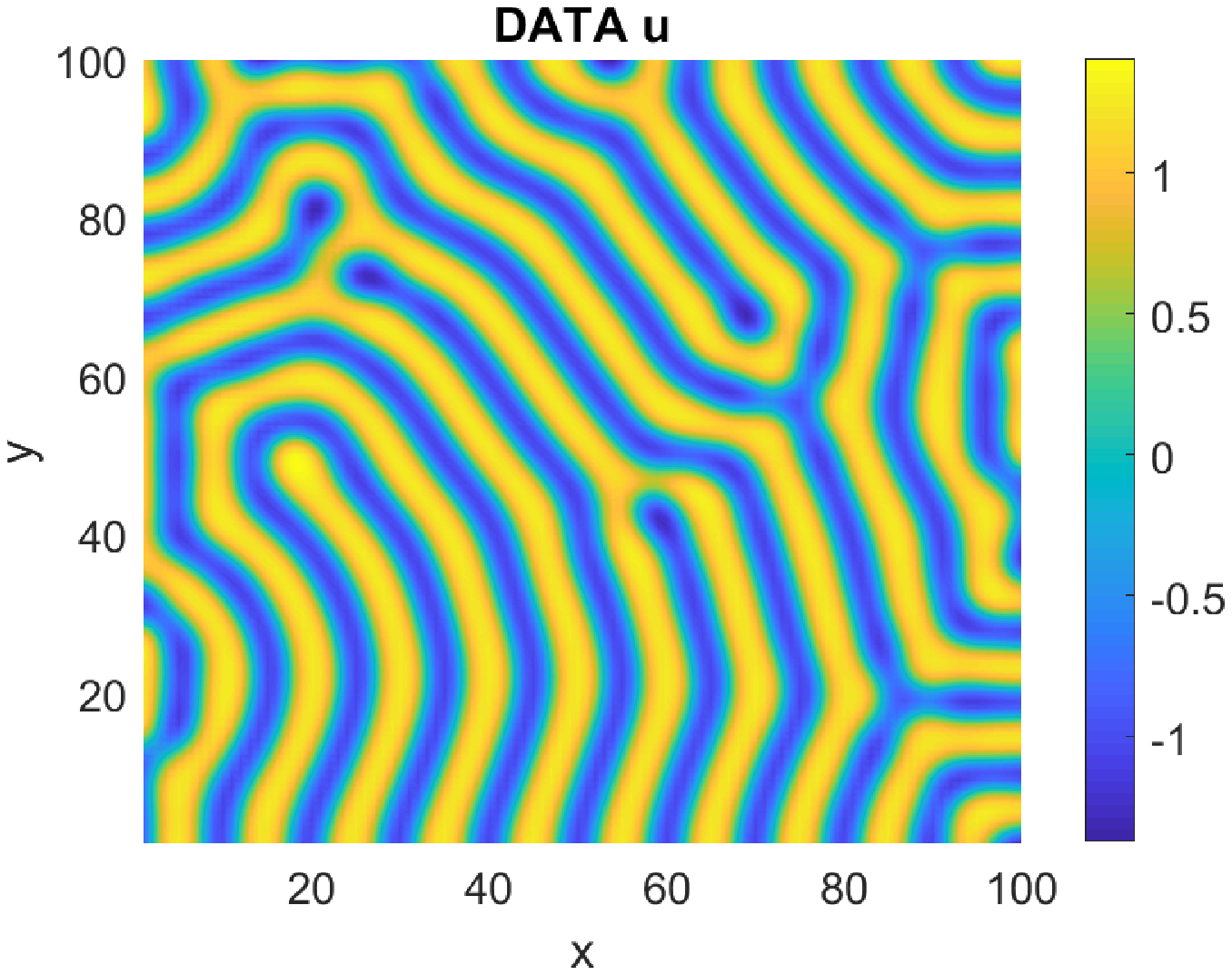}
\includegraphics[scale=0.4]{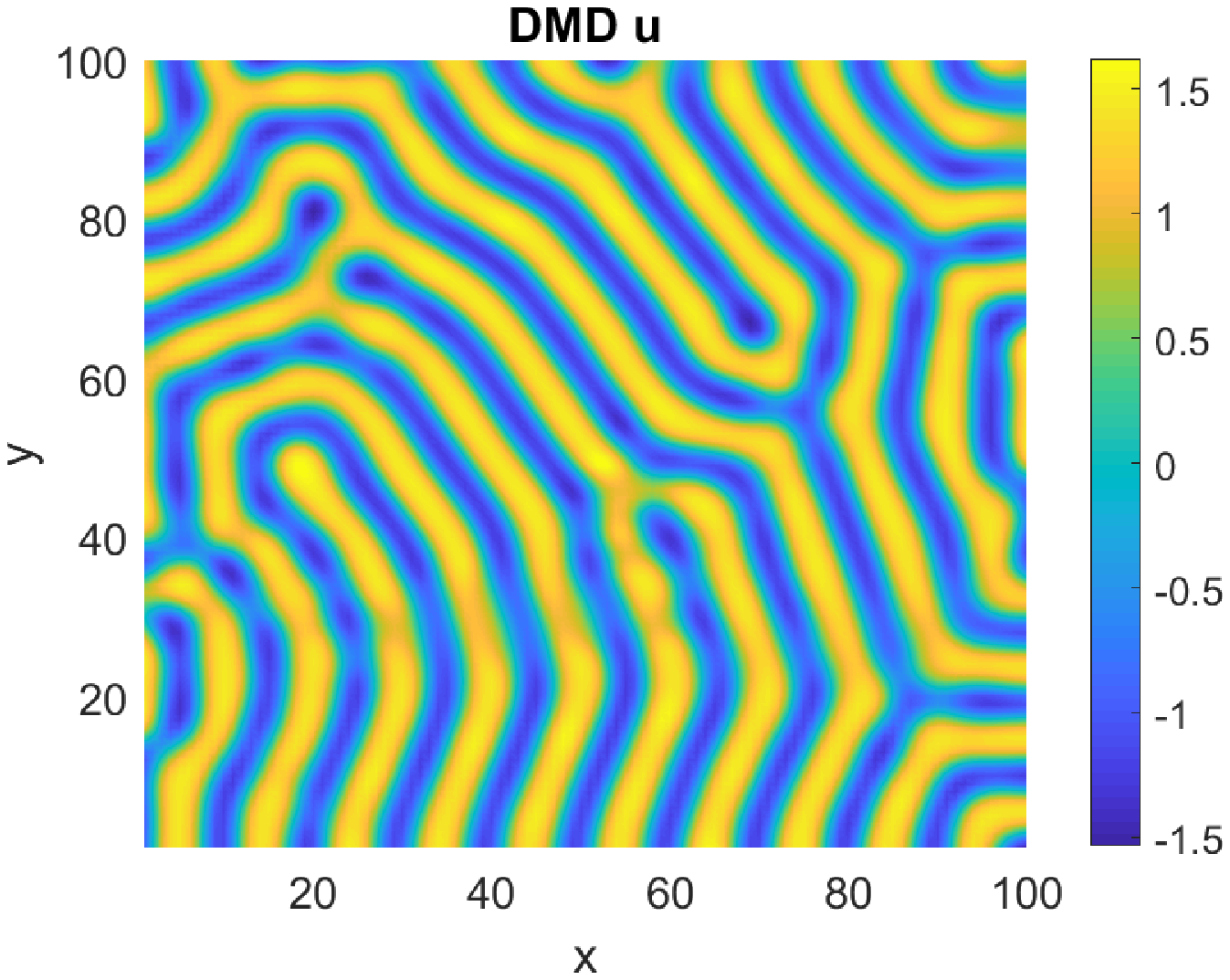}
\includegraphics[scale=0.4]{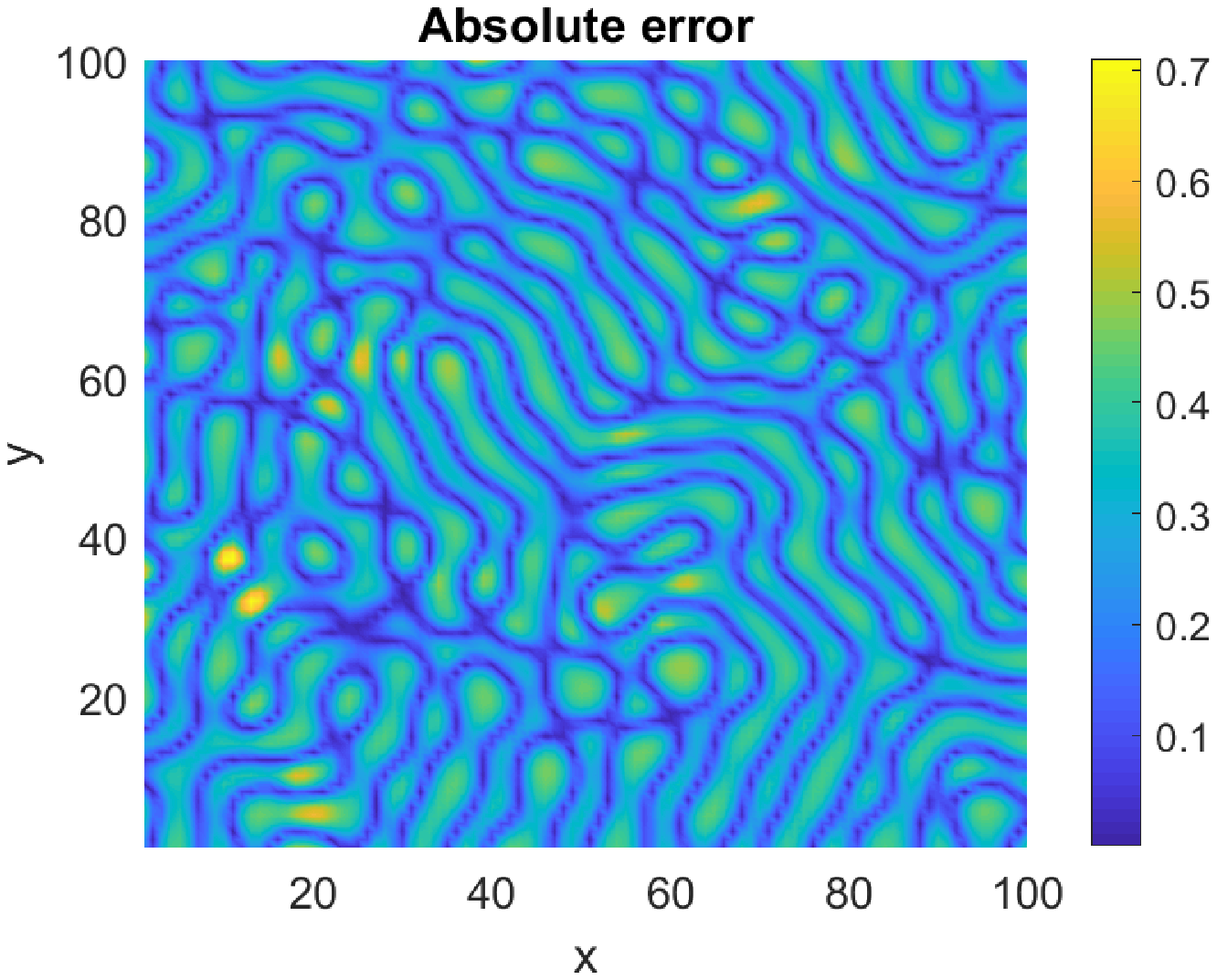}\\
\includegraphics[scale=0.4]{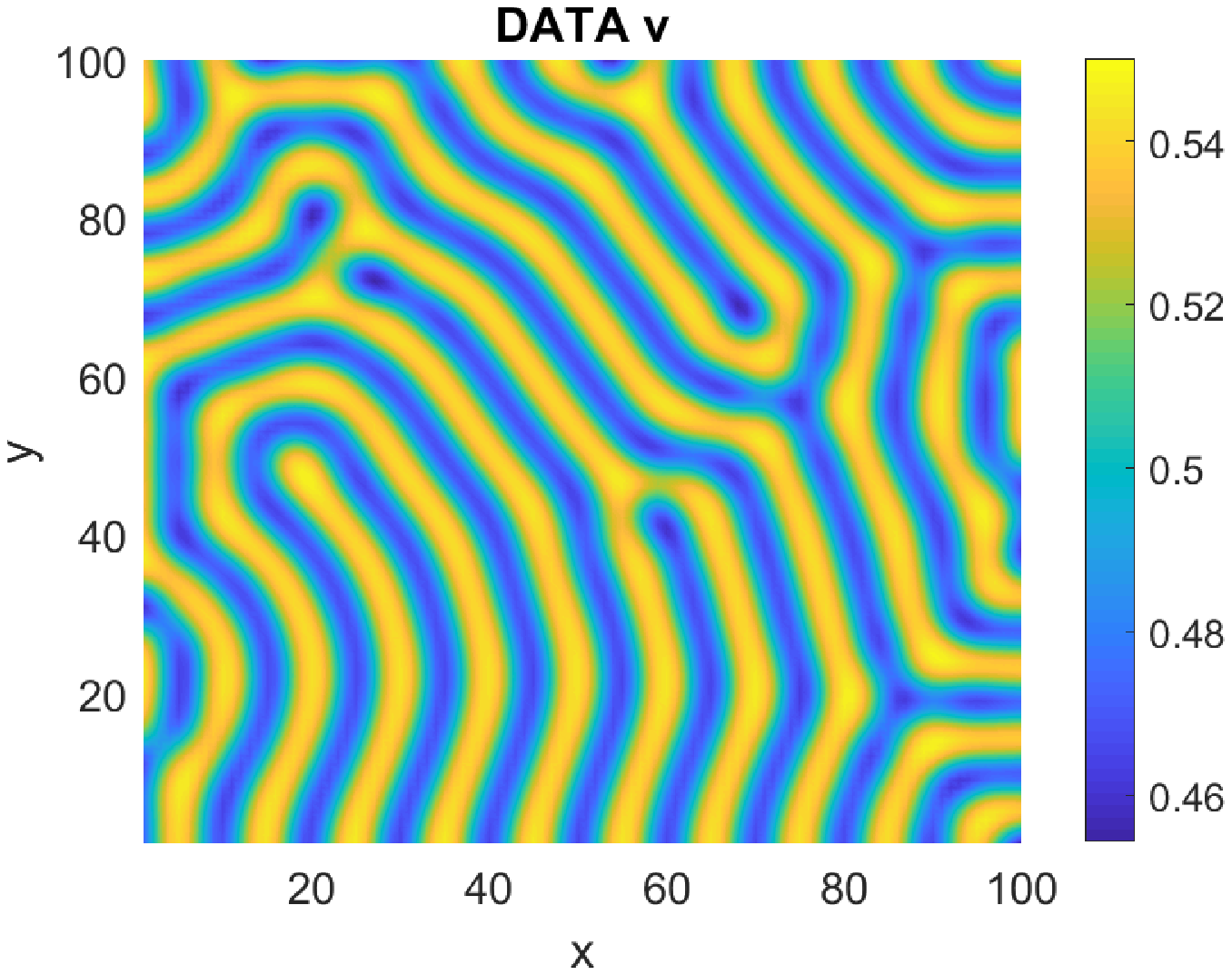}
\includegraphics[scale=0.4]{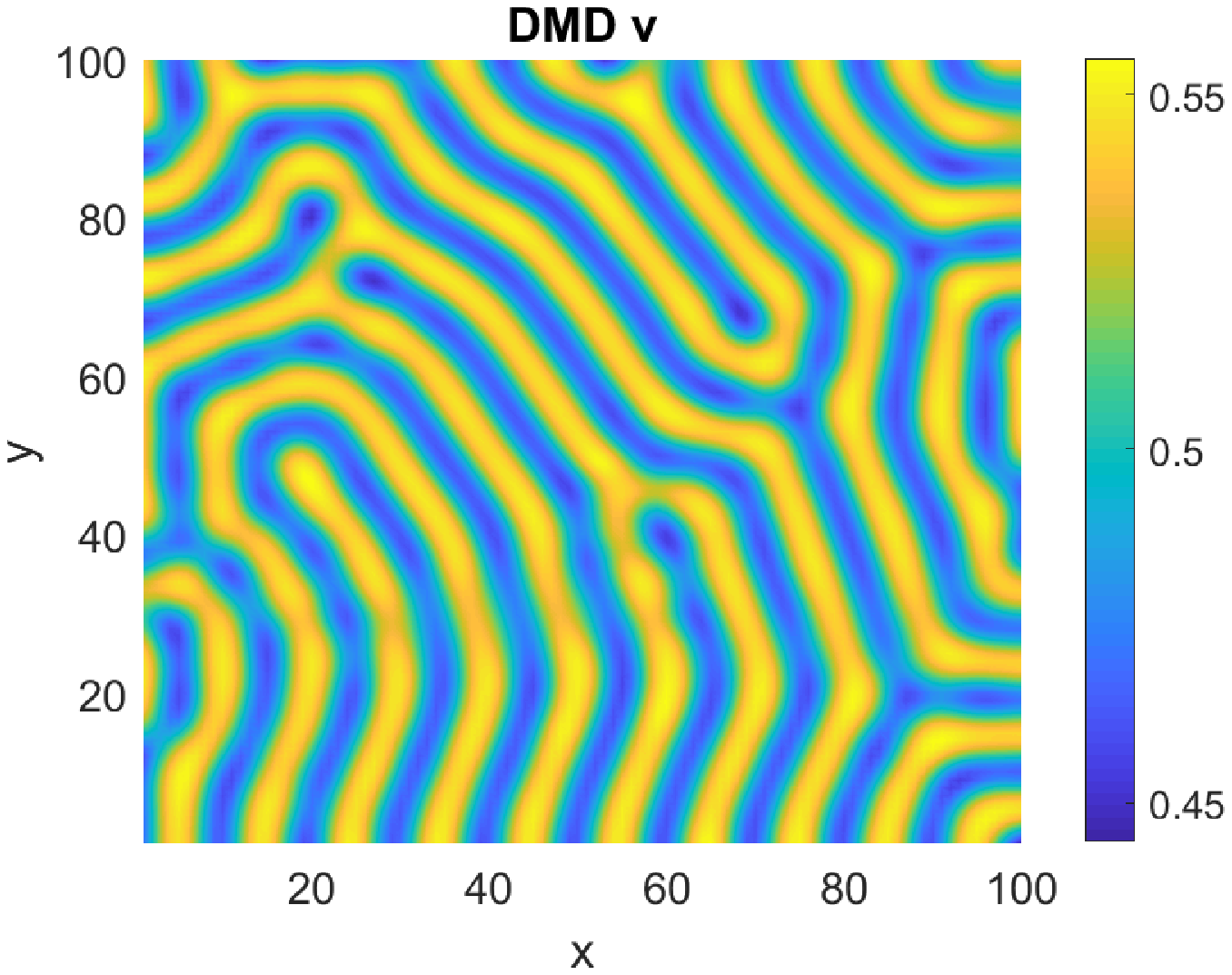}
\includegraphics[scale=0.4]{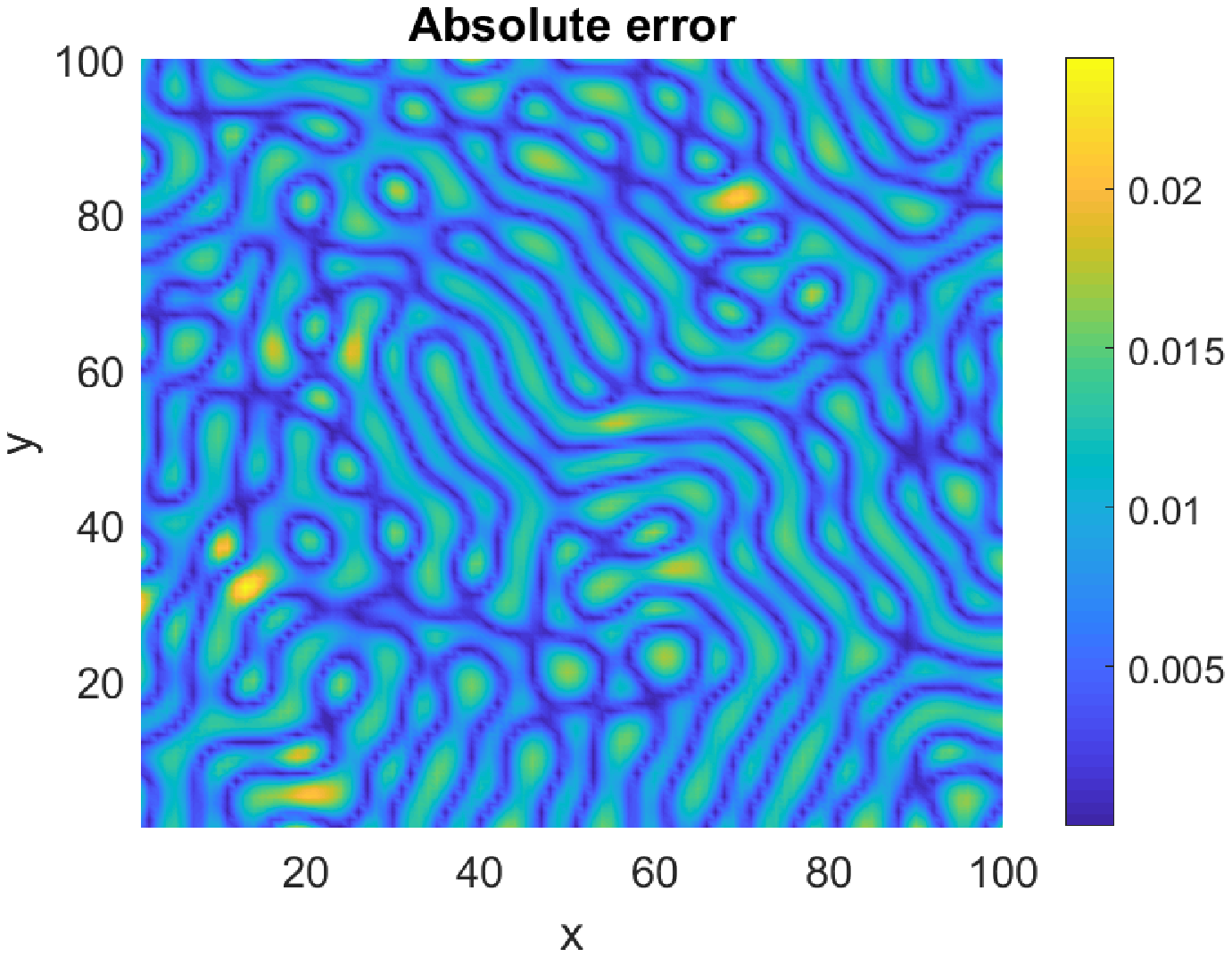}
\captionsetup{justification=justified}
\caption{DIB model: Turing instability. Full model solutions (left) and DMD reconstruction (center) with $r = 22$. Right plot: absolute error of the DMD reconstruction with respect to data.}
\label{dib_dmd_classic}
\end{figure}

In Figure \ref{dib_dmd_classic}, we show the full model solution, i.e. the labyrinth Turing patterns both for the variables $u$ and $v$ attained at the final time $T = 40$ (left plots) and their DMD reconstructions with $r= 22$ (center plots). The DMD approximates quite well the shape of the final patterns, although the amplitude is not correct, as highlighted by computing the spatial absolute errors between them that are reported in the right panels of Figure \ref{dib_dmd_classic}.
%

In conclusion, also in this example we have shown how DMD does not approximate properly the dataset. In addition, the relative error gets worse when the rank increases and it never gets lower than $10\%$. It is worth noting that a similar bad DMD behaviour for the DIB model was already discussed in \cite{BMS21}. Even if here, we consider a different DMD implementation, that is a coupled approach based on a randomized version of DMD to improve the computational efficiency. This improvement is not able to remove the DMD drawbacks for Turing pattern approximation.

%
\subsection{Example with spatio-temporal oscillatory dynamics in the Turing-Hopf instability}
\label{ex_dib_hopf}
The last class of problems in exam exhibits a Turing-Hopf instability, that is an interplay between Turing and Hopf instabilities (\cite{DIB15}). In particular, the solutions of \eqref{RDPDE} are oscillatory patterns both in space and time. 
We consider the DIB morpho-chemical model, whose kinetics in \eqref{RDPDE} are the same as in \eqref{DIB_kin}, but the new dynamics arises for different model parameter values given by:
$A_2 = 30, \ B = 109, \ C = 2.794, \ \rho = 50.$
The initial conditions are again spatially random perturbation of the homogeneous equilibrium, as in the previous section.

We discretize the rectangular spatial domain $\Omega =[0,100] \times [0, 70]$ with $n_x = n_y = 100$ spatial meshpoints, thus $n = n_x n_y = 10000$. We integrate in time (IMEX-Euler in matrix oriented form) with time step $h_t = 10^{-4}$ and final time $T = 4.5$. We emphasize that, to the best of author's knowledge, the IMEX Euler scheme in the matrix form has never been applied to this kind of problems. This approach reduces significantly the computational execution time with respect to a standard vector form, therefore it allows to speed-up the offline stage, that is the construction of the dataset $S$.
We save the snapshots every four time steps ($\kappa=4$, see discussion in Section \ref{sec:fom}), such that the dataset is $S \in \RR^{2n \times (m+1)}$  with $m +1= 11250$. 

\begin{figure}[htbp]
\centering
\includegraphics[scale=0.45]{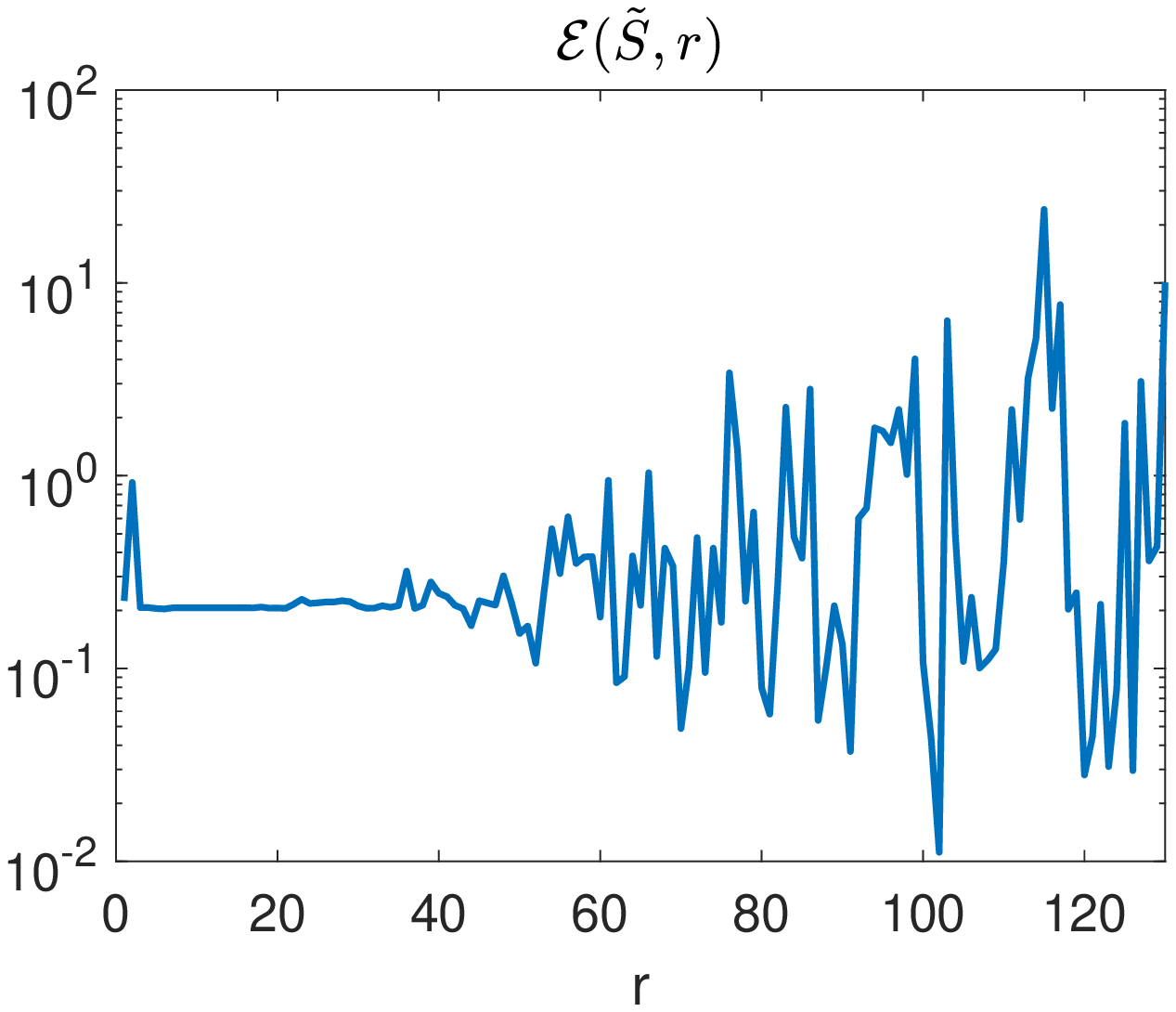}
\includegraphics[scale=0.45]{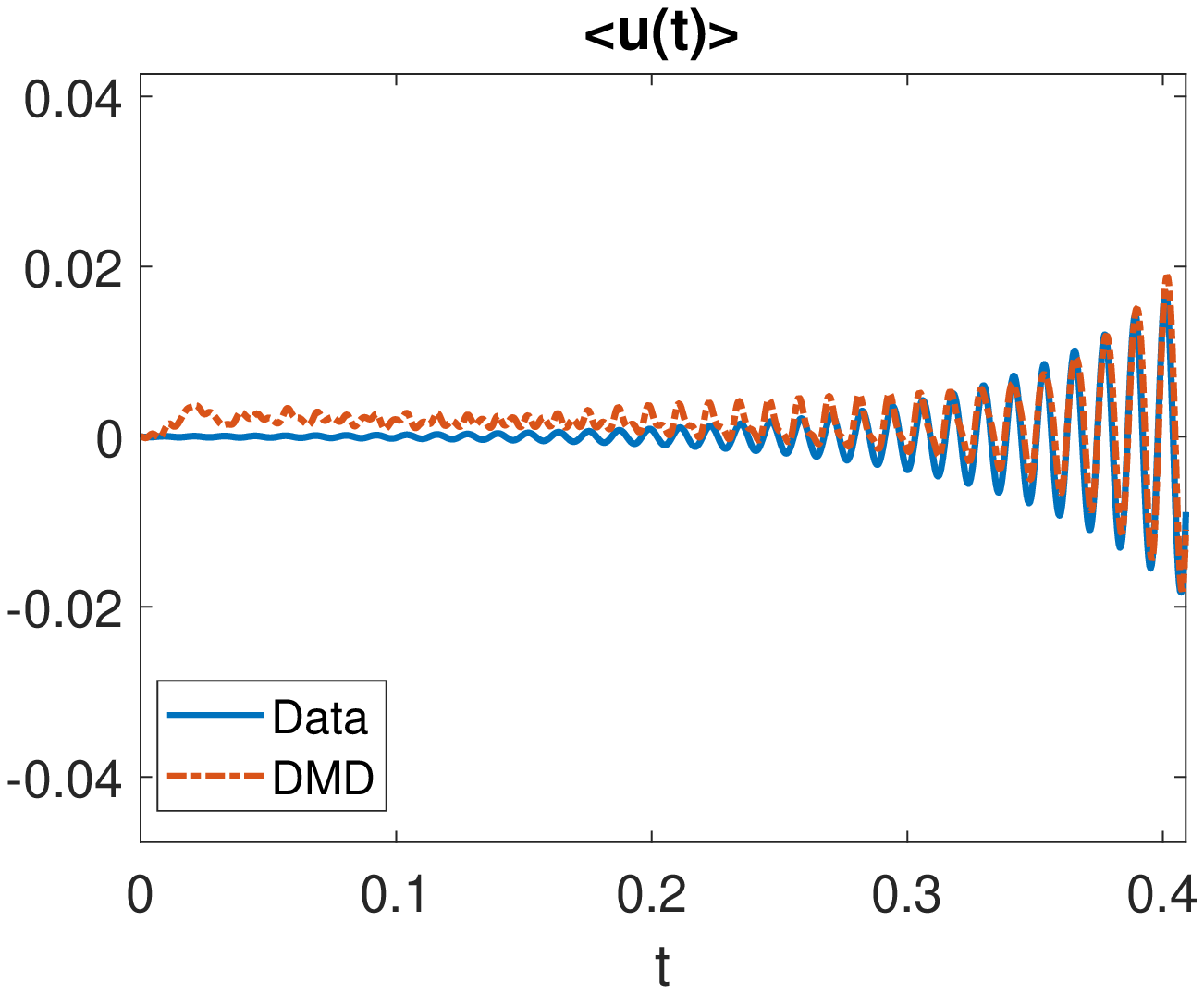}

\includegraphics[scale=0.45]{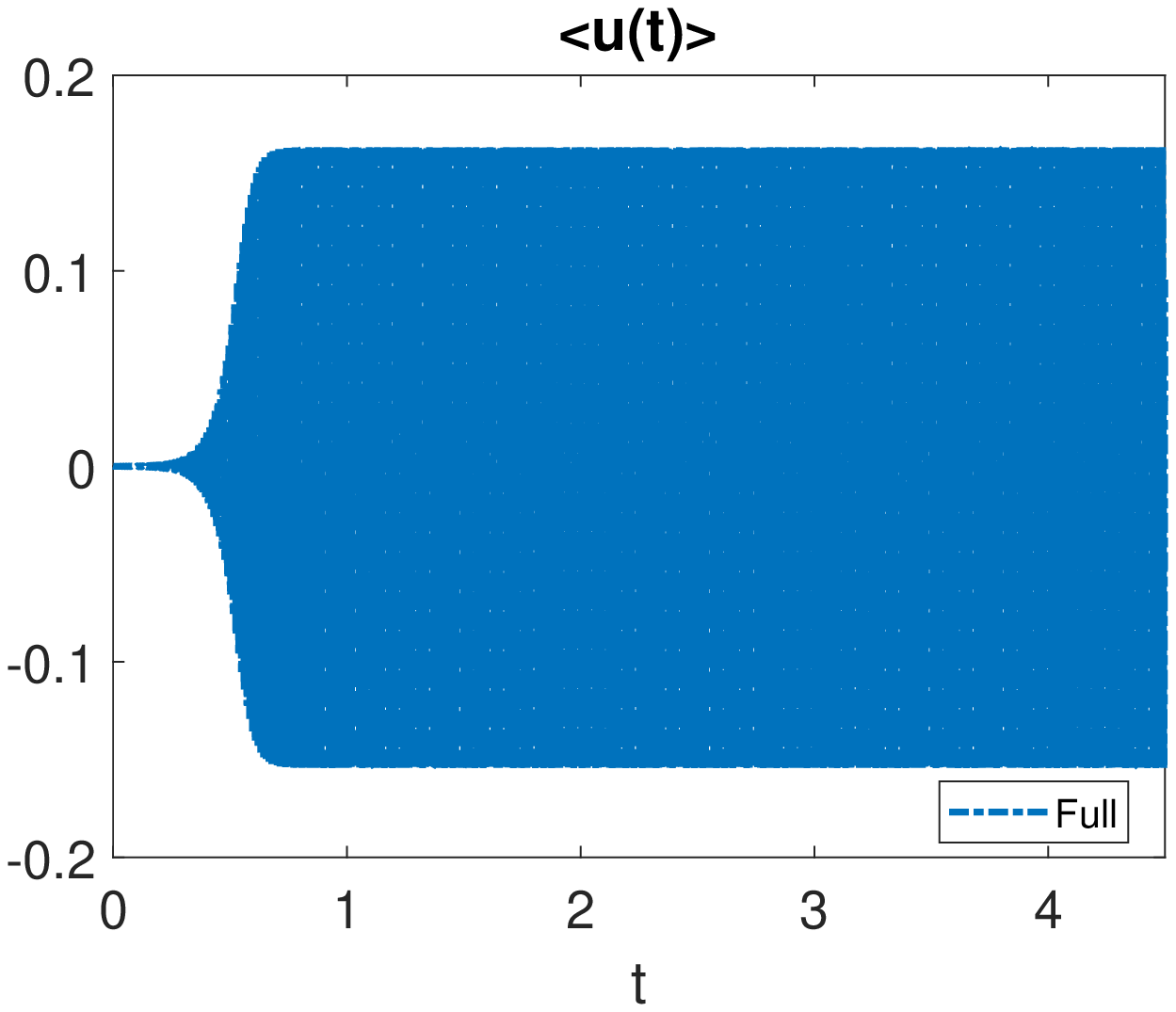}
\includegraphics[scale=0.45]{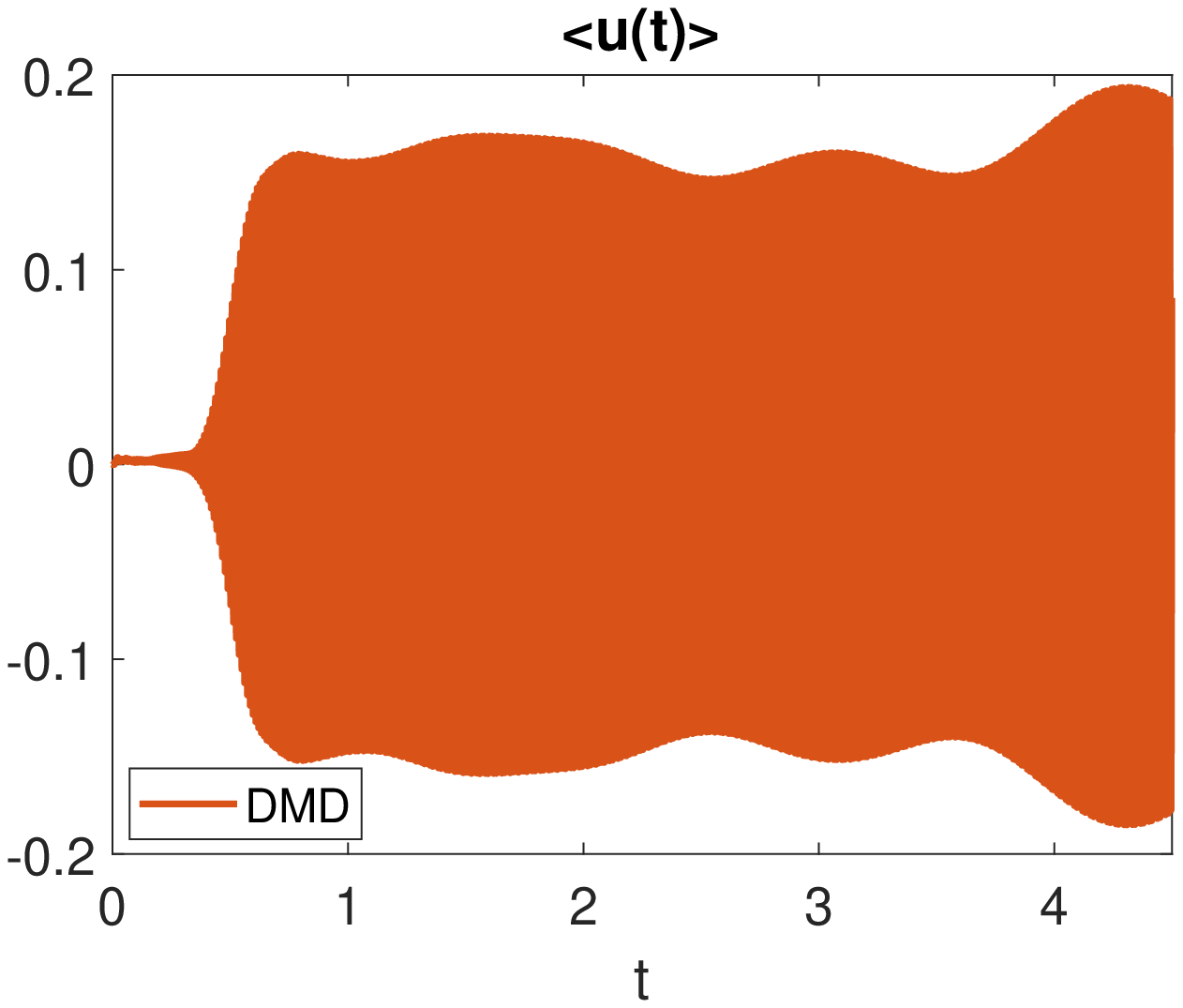}

\captionsetup{justification=justified}
\caption{DIB model: Turing-Hopf instability. Top left plot: DMD relative error in \eqref{ef_data} for the dataset $S$. Top right plot: zoom on the spatial means. Bottom left plot: spatial mean for the full model solution $u$; bottom right plot: spatial mean for the DMD reconstruction with $r = 102$. }
\label{dib_hopf_mean_classic}
\end{figure}
\begin{figure}[htbp]
\centering
\includegraphics[scale=0.4]{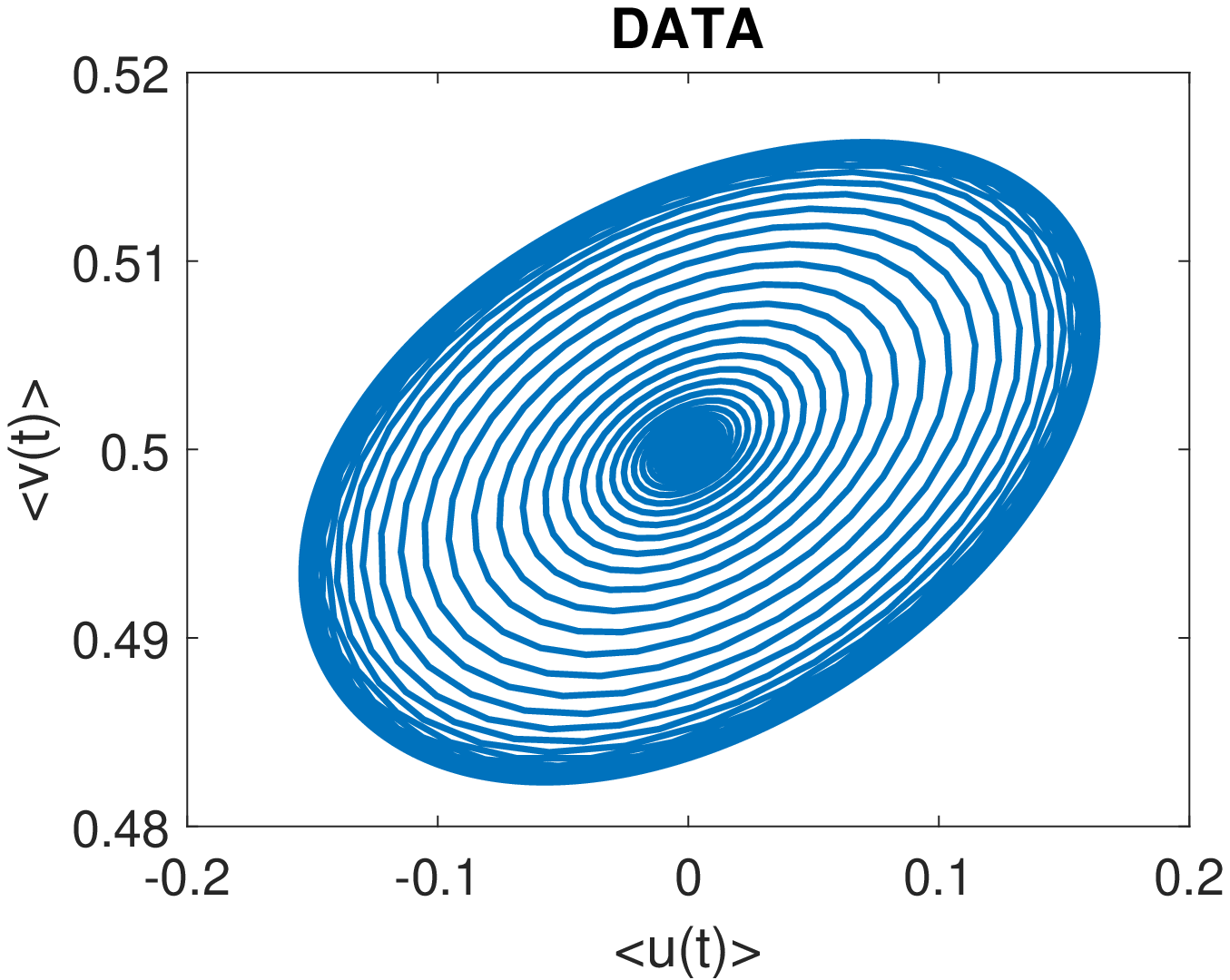}
\includegraphics[scale=0.4]{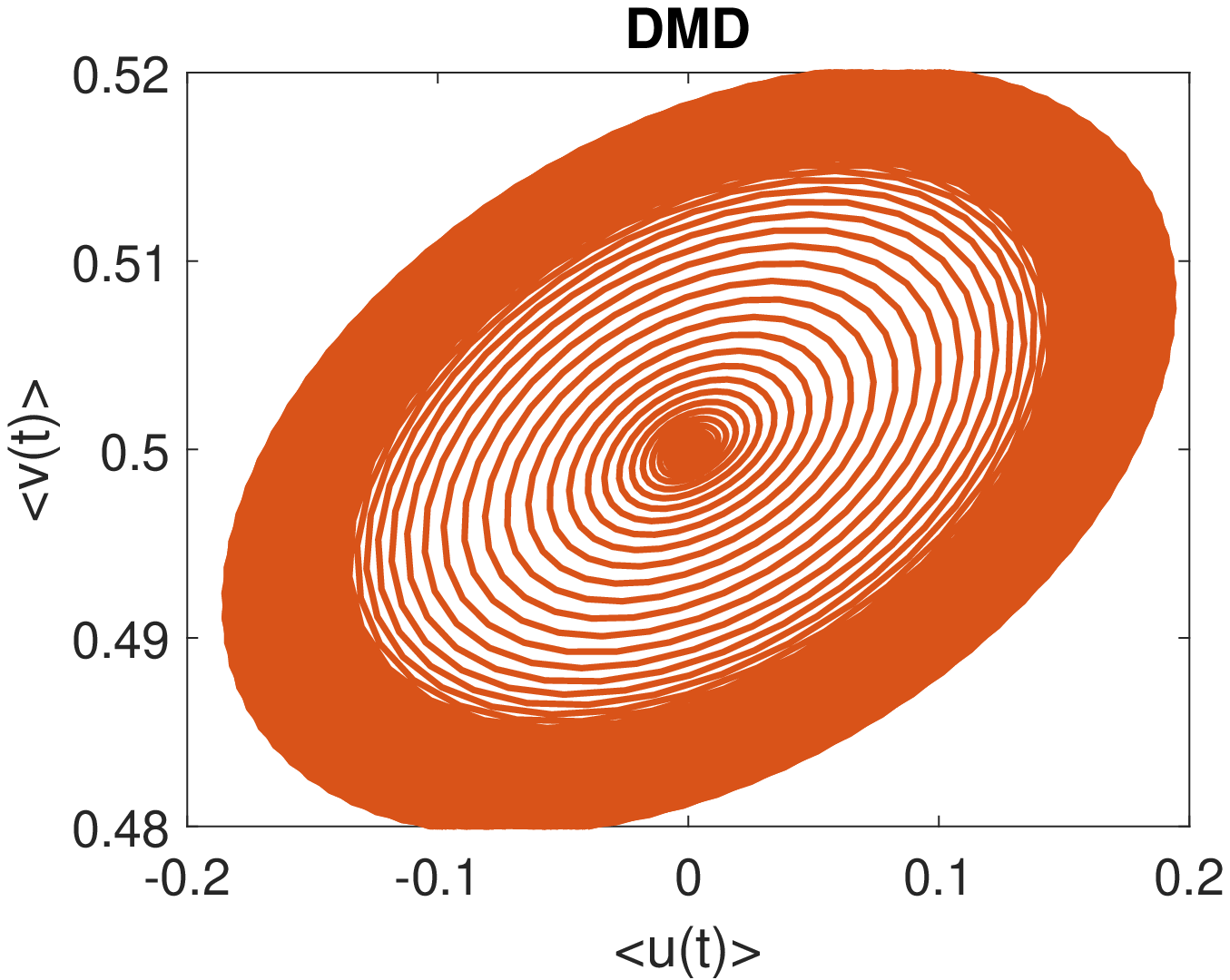}
\captionsetup{justification=justified}
\caption{DIB model: Turing-Hopf instability. Comparison in the phase plane of data (left) and DMD with $r = 102$ (right).} 
\label{dib_hopf_phase}
\end{figure}

%
%
In the top left panel of Figure \ref{dib_hopf_mean_classic}, we show the relative error $\mathcal{E}(\widetilde S, r)$ for $r = 1, \dots, R$, where $R=130$ is the rank of the snapshot matrix $S$. We can observe that the behaviour is very erratic and the minimum is $\mathcal{E}(\widetilde{S},102) = 0.0111$, reached for $r=102$. Different time dynamics of the spatial means are obtained, as shown in the bottom plots of Figure \ref{dib_hopf_mean_classic}, for the data $u$ (left) and for DMD reconstruction with $r=102$ (right). In the top right picture we report a zoom of both spatial means over $0<t<0.5$ to show that DMD does not match the mean of the dataset also in the transient regime. Instead for $t>\simeq1$, DMD catches the frequency but not the amplitude of the spatial mean.
To further emphasize the above significant difference, in Figure \ref{dib_hopf_phase} we also report the corresponding limit cycles obtained in the phase plane $(\langle u \rangle, \langle v \rangle)$ for the data (left plot) and DMD reconstruction with $r = 102$ (right plot). 

\begin{figure}[htbp]
\centering
 \includegraphics[scale=0.4]{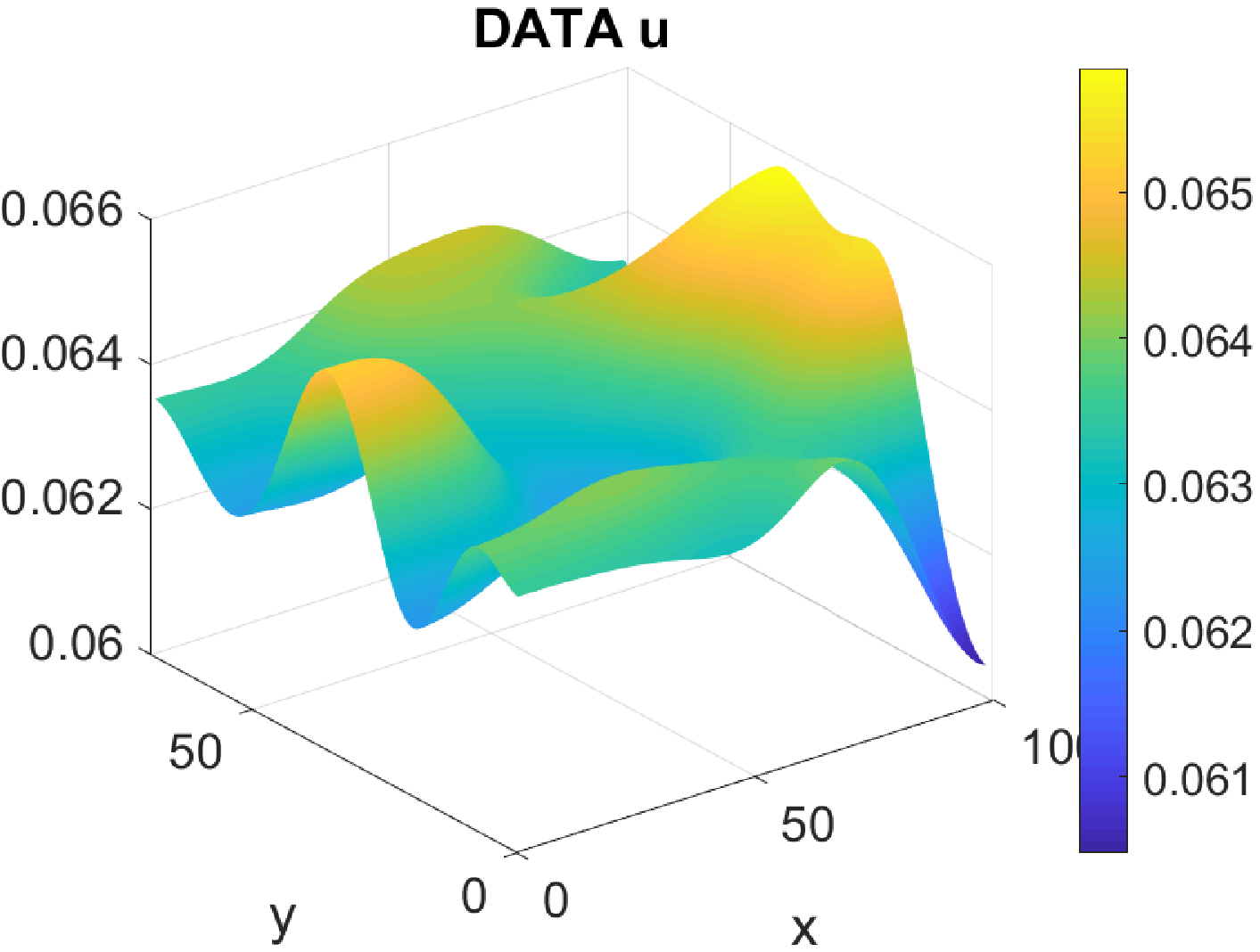}
 \includegraphics[scale=0.4]{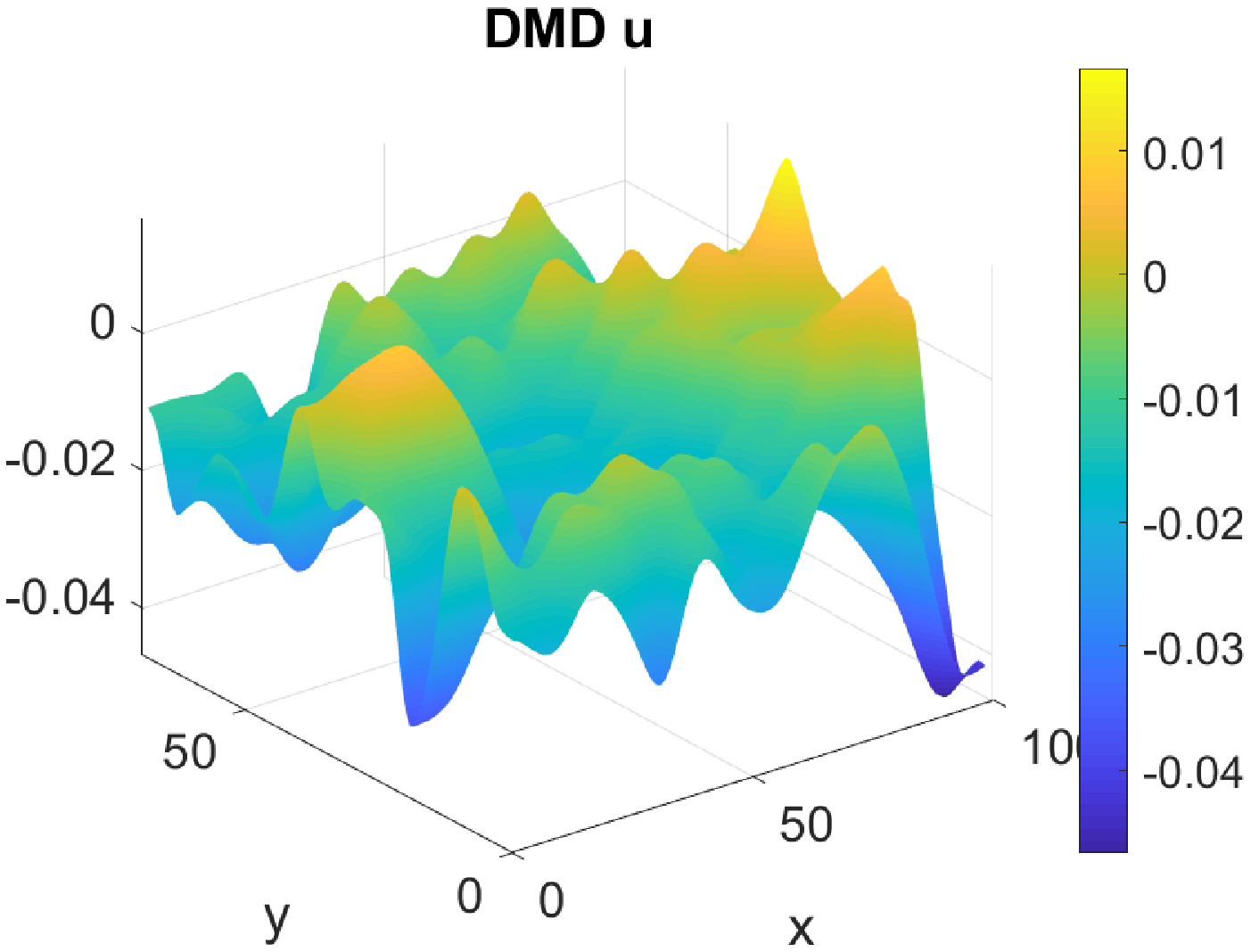}
 \includegraphics[scale=0.4]{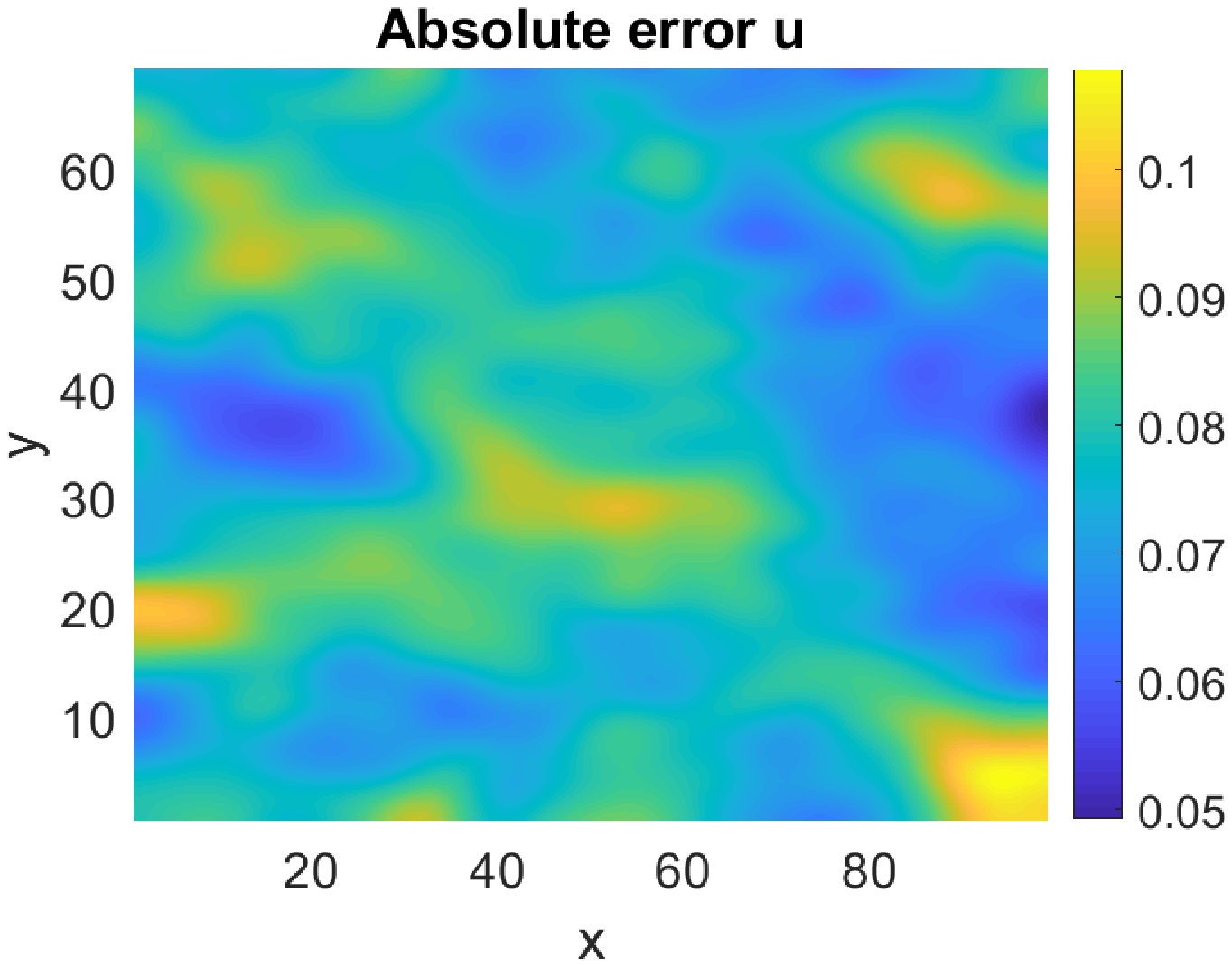}\\
 \includegraphics[scale=0.4]{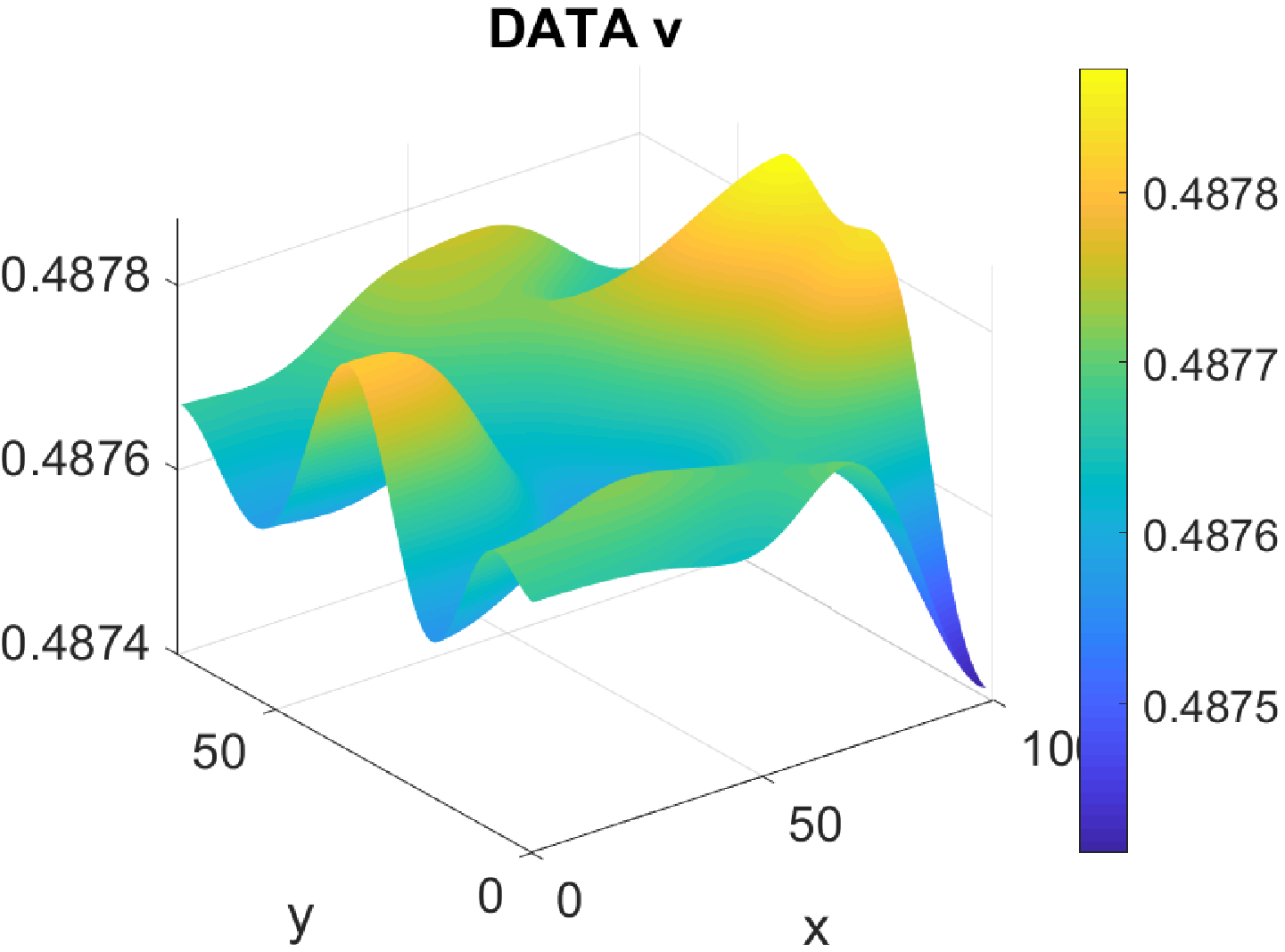}
 \includegraphics[scale=0.4]{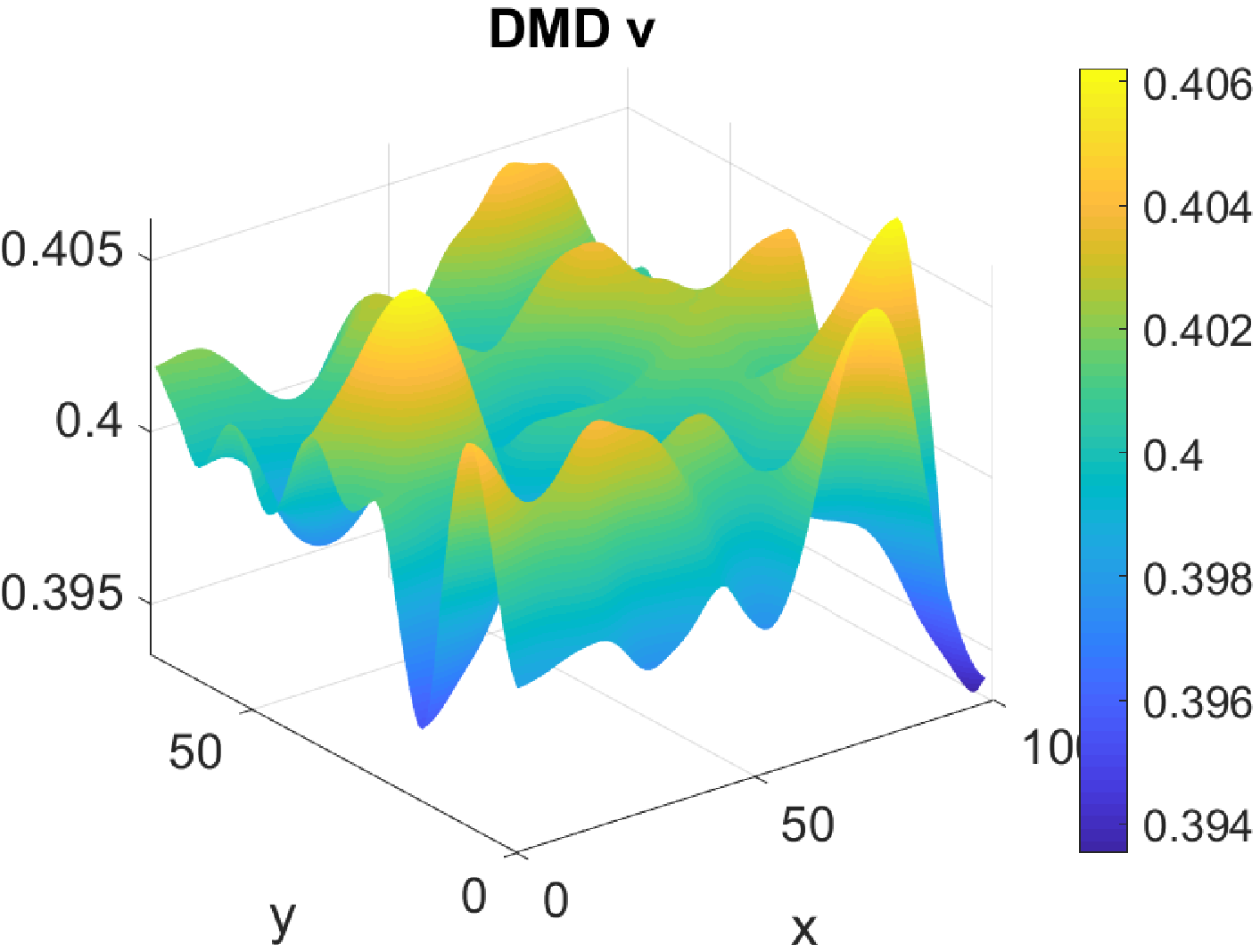}
 \includegraphics[scale=0.4]{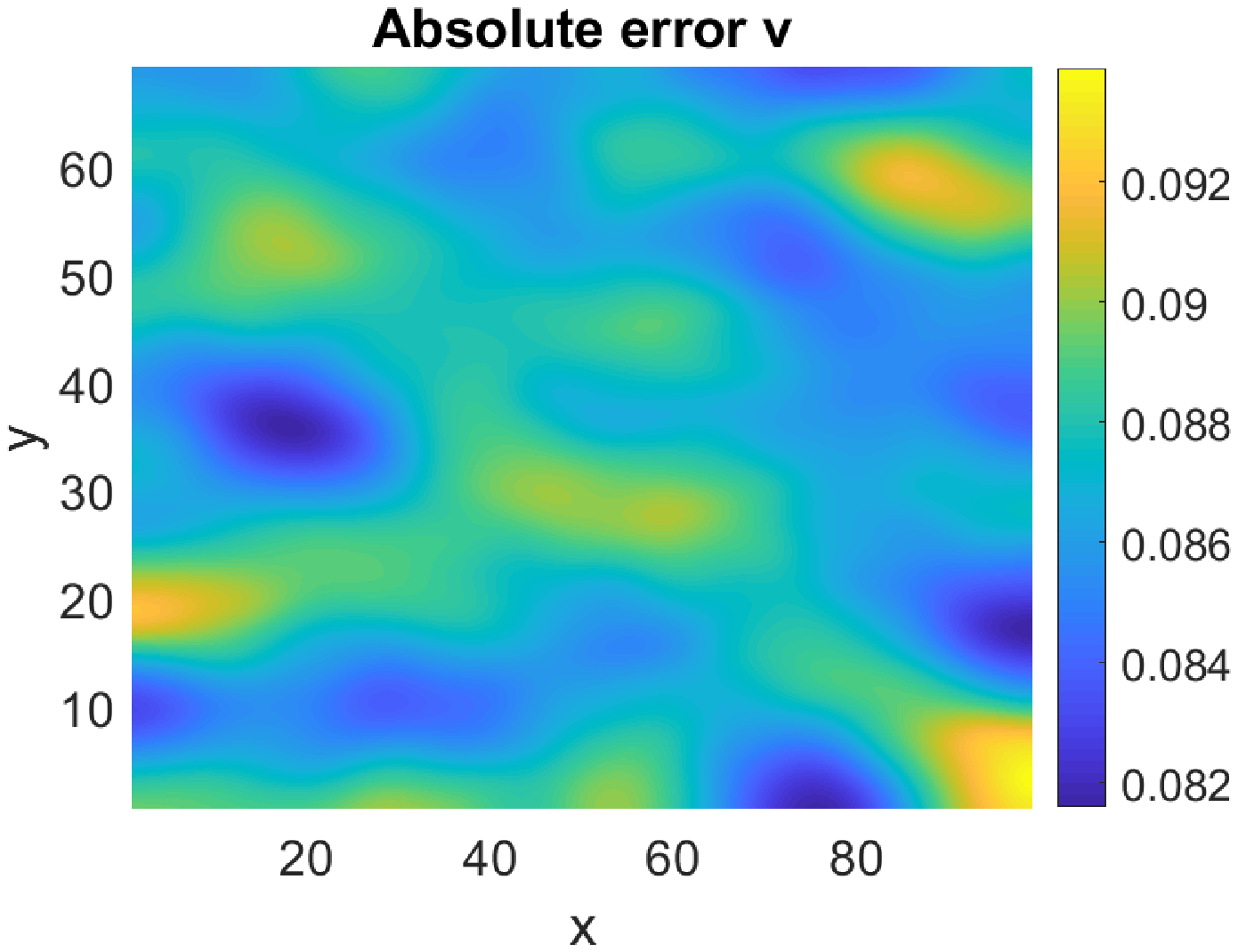}
\captionsetup{justification=justified}
\caption{DIB model: Turing-Hopf instability. Full model solutions (left panels) and DMD reconstruction (center panels) with $r = 102$. Right plots: absolute error of the DMD reconstruction with respect to data.} 
\label{dib_hopf_dmd_classic}
\end{figure}

To complete the discussion on this example, where the more complicated spatio-temporal oscillating dynamics is presented, In Figure \ref{dib_hopf_dmd_classic} we report the pattern solutions at the final time $T = 4.5$: the left panels concern the full model solutions $u$ and $v$, the middle ones are for the corresponding DMD reconstructions. As expected from the previous results on the temporal dynamics, DMD does not approximate accurately also the final patterns both for $u$ and $v$, as confirmed quantitatively from the absolute errors with respect to the data shown in the right panels.
%

%

\section{The piecewise DMD method}
\label{sec:pdmd}
In the previous section, we have shown a selection of examples where DMD fails its reconstruction during the time dynamics and  in addition the relative error with respect to the dataset increases when the rank does. It is worth remarking that we have tested the Higher Order DMD (HODMD, \cite{LV17, LV17b}) on the datasets discussed in the previous section without any improvement on the quality of the DMD approximation.
  In this section, we propose a new approach to tackle safely the spatio-temporal features of the peculiar solution dynamics discussed so far, that is datasets with oscillating behaviours and pattern formation by Turing instability.


The main idea here is to propose a {\it piecewise} version of the DMD algorithm, that we will denote by pDMD.
In fact, we suppose that for oscillatory and Turing spatio-temporal dynamics, the main assumption underlying the original DMD, that is a ``global'' linear fitting over the full temporal horizon, is not sufficient to recognize different ``phenomena'' arising along the time pathways. Therefore, instead of performing a DMD approximation on the whole time interval $[0, T],$ we propose to decompose it and the corresponding dataset into $N \geq 1$ parts, as follows.

Consider $\nu =\lceil \frac{m+1}{N}\rceil \geq \nu^*$ and the dataset decomposition
$S = \cup_{i=1}^{N} S_i$, where $S_i$ is the submatrix of $\nu$ columns of $S$ defined by $S_i = [S_{:,(i-1)\nu+1}, \dots, S_{:,i\nu}] \in \RR^{2n \times \nu}$ for $i = 1, \dots, N$. In practice, $S_i$ corresponds to consider those snapshots of $S$ belonging to the time interval $[t_{(i-1)\nu}, t_{i\nu -1}].$ 

We suppose that $\nu^* \geq 10$, such that a minimum number of snapshots in each subset $S_i$ is guaranteed and a maximum value for $N$ can be chosen.
Then, we apply the DMD technique of rank $r_i$, using the randomized version based on the QB decomposition recalled in Algorithm \ref{alg:qb}, on each subset $S_i, i = 1,\ldots, N,$ with the convention that for $N = 1$, we recover the original dataset $S_1\equiv S$. In Algorithm \ref{Alg:PiecDMD} below, we present in details the piecewise DMD, defined as pDMD, returning in output, not only the reconstructed snapshots, but also the vector $\r = [r_1, r_2, \dots, r_N] \in \RR^N$ accounting for the ranks considered on each dataset $S_i$. We can suppose to fix {\it a priori} the values in $\r$ or to estimate them. 


Our idea stems from the ``divide and conquer'' approach well known in the numerical analysis framework to reduce the ``global error'' in the approximation under exam, as, for example, the piecewise interpolation and composite quadrature rules. Then, due to the basic meaning of the DMD, recalled in Section \ref{sec:dmd}, in each subinterval $[t_{(i-1)\nu}, t_{i \nu-1}]$, for $i=1,\ldots, N$ we will get the best linear fit of the form \eqref{hid:dyn} on the dataset portion therein. Hence, this local linearization can help to look at the different solution regimes along time by applying multiple separated linear fittings. In fact, we will show that local/piecewise linear fit can follow better the switches between these regimes instead of the global ($N=1$) approach that could miss them. 
%

We summarize the proposed method pDMD in Algorithm \ref{Alg:PiecDMD} and below we comment it step by step.
\begin{algorithm}[htbp]
\caption{Piecewise DMD (pDMD)}
\label{Alg:PiecDMD}
\begin{algorithmic}[1]
\STATE {\bf INPUT} Dataset $S \in R^{2n \times (m+1)}$ in $[0,T]$, threshold $\overline{tol} >0$
\STATE {\bf OUTPUT} $\widetilde{S}^{N}$ piecewise reconstruction in $[0,T],$ ${\bf r}$ the ranks used in each partition
\STATE Choose an initial number $N$ of partitions
\STATE Split the datasets $S = \cup_{i=1}^{N} S_i$
\FOR{$i = 1, \ldots, N$}
\STATE set the target rank $r_i=\mbox{rank}(S_i)$  
\STATE compute the (randomized) DMD solution $\widetilde{S}_i$ of rank $r_i$ using Algorithm \ref{alg:qb}
\STATE compute the error $err(i)$ defined in \eqref{e_max} 
\IF{ $err(i) > \overline{tol}$}
\STATE    $N = N+1$ 
\STATE   go to step 4
\ENDIF
\ENDFOR
\STATE $\widetilde S^{N} = \cup_{i = 1}^{N} \widetilde{S}_i, {\bf r}=[r_1,\ldots, r_N].$
\end{algorithmic}
\end{algorithm}

\noindent
\paragraph{\bf Inputs} The inputs of the algorithm are the dataset $S$ and a desired threshold $\overline{tol}$ for the pDMD.
\paragraph{\bf Initialization} We choose an initial value $N$ to construct the partition of the dataset. One can easily start with $N = 1$, which corresponds to the standard DMD discussed in Section \ref{sec_fail_dmd}. We then split the dataset $S\in \RR^{2n \times (m+1)}$ in $N$ parts and build all the matrices $S_i \in \RR^{2n \times \nu}, i=1,\ldots,N, \nu=\lceil\frac{m+1}{N}\rceil$.
\paragraph{\bf pDMD} For each subset $S_i$ we fix the rank $r_i$ using the QB decomposition in Algorithm \ref{alg:qb}. We then compute the following relative error
\begin{equation}
\label{e_max}
err(i) = \max_{(i-1)\nu +1 \leq k \leq i \nu} \frac{\|\x_k - \tilde{\x}_k \|_{\infty}}{\| \x_k \|_{\infty}}
\end{equation}
that computes the worst approximation in each interval $[t_{(i-1)\nu +1}, t_{\nu +1}]$. This error is faster to compute than \eqref{ef_data}.
It is worth noting that the computation of \eqref{e_max} is always possible since we are dealing with DMD reconstruction on the training dataset which is the focus of this work.

If $err(\bar{i})>\overline{tol}$ for some $\bar{i}\in\{1,\ldots, N\}$, we do not compute DMD for {$i>\bar{i}$}, but we increase the value of $N$ and restart the method with a new finer partition. That error indicator is, indeed, meant to save computational time whenever possible. In step 10 of Algorithm \ref{Alg:PiecDMD}, we propose to increase $N$ by $1$, but clearly other choices can be used such as e.g. $N=N+\widetilde{N}$ with $\widetilde{N}\in\mathbb{N}$. Note that $\widetilde{N}=1$ corresponds to the choice in Algorithm \ref{Alg:PiecDMD}. We iterate until the desired convergence is reached.
\paragraph{\bf Output}The Algorithm returns the pDMD reconstruction $\widetilde{S}^{N}$ that is the union of the DMD subsets $\widetilde{S}_i$, for $i = 1, \dots, N$ and the ranks ${\bf r}=[r_1,\ldots, r_N]$ used in each partition.

\begin{rmk}[Choice of the rank]
It is important to note that the rank values in this algorithm plays a crucial role. One can always assume to work with a full rank approximation in each $S_i$, as set in Algorithm \ref{Alg:PiecDMD} but this might be computationally expensive, especially when $m$ is large and $N$ is still small. In some simulations, to avoid the computation of the rank for large matrices, we choose the target rank in step 6 of Algorithm \ref{Alg:PiecDMD} as $r_i = \min\{\nu, 200\}$, where $\nu$ is the number of snapshots in $S_i$. One could also use randomized rank revealing methods as proposed in e.g. \cite{HIQM22,MQH19}, but this is out of the scope of this paper.
\end{rmk}%
%
%
%
\section{Numerical experiments based on the Piecewise DMD method}
\label{sec:testpdmd}
In this section, we apply the proposed pDMD Algorithm \ref{Alg:PiecDMD} to the datasets generated from the RD--PDE models presented in Section \ref{sec_fail_dmd} where several drawbacks of the classical DMD have been discussed. 
For each test, we will show different error indicators, as follows.
First of all, we consider the relative error $\mathcal{E}_p$ in Frobenius norm between the dataset $S$ and its piecewise DMD reconstruction $\widetilde{S}^{N}$ defined by: 
\begin{equation}
\label{rel_err}
\mathcal{E}_p(\widetilde{S}^{N}, {\bf r}) = \frac{\| S - \widetilde{S}^{N} \|_F}{\| S \|_F},
\end{equation}
depending from the number $N$ of partitions used to split the whole dataset $S$ and {${\bf r} = [r_1, r_2, \dots, r_N] \in \RR^N$ the vector of all ranks considered, such that $r_i$ is used for the subset $S_i$. We observe that for $N = 1$ we recover the error corresponding to  the "global" DMD approach, i.e. $\mathcal{E}_p(\widetilde{S}^{1},r_1) = \mathcal{E}(\widetilde{S},r_1)$ defined in \eqref{ef_data}. \\
In our numerical tests we will provide the behaviour of the error  $\mathcal{E}_p(\widetilde{S}^{N}, {\bf r})$ of the pDMD method for different choices of $N$, that are obtained by the inner computations of the following Algorithm \ref{Alg:PiecDMD_N}.
\begin{algorithm}[htbp]
\caption{Convergence of pDMD }
\label{Alg:PiecDMD_N}
\begin{algorithmic}[1]
\STATE {\bf INPUT} Dataset $S \in R^{2n \times (m+1)}$ in $[0,T]$, threshold $tol>0$,
\STATE {\bf OUTPUT} $\widetilde{S}^{N}$ piecewise reconstruction in $[0,T]$,  ${\bf r}$ the ranks used in each partition, the error $\mathcal{E}_p(\widetilde S^{N},{\bf r})$
\STATE Choose an initial number $N$ of partitions
\STATE Split the datasets $S = \cup_{i=1}^{N} S_i$
\WHILE{$\mathcal{E}_p(\widetilde S^{N},{\bf r})>tol$}
\STATE compute $\widetilde S^{N}$ and ${\bf r}$ from Algorithm \ref{Alg:PiecDMD} 
\STATE N = N+1
\ENDWHILE
\STATE $\widetilde S^{N} = \cup_{i = 1}^{N} \widetilde{S}_i, \; {\bf r}=[r_1,\ldots, r_N], \; \mathcal{E}_p(\widetilde S^{N},{\bf r})$
\end{algorithmic}
\end{algorithm}

As second indicator, for some values of $N$ identified by  Algorithm \ref{Alg:PiecDMD_N}, we will also check the relative error in Frobenius norm over time calculated by:
\begin{equation}
\label{error_time}
\epsilon_k(\widetilde{S}^N,{\bf r}) = \frac{\| \x_k - \tilde{\x}_k \|_F}{\| \x_k \|_F}, \quad k = 0, \dots, m
\end{equation}
where $\{ \tilde{\x}_k \}_{k=0}^m$ are the snapshots reconstructed by the pDMD with $N$ partitions.

We recall that the target rank value $r_i$ for the DMD reconstruction is fixed on each subset $S_i \in \RR^{2n \times \nu}$ and it depends on $N$, because each $r_i \leq \nu=\lceil\frac{m+1}{N}\rceil$. 
Of course, the computational load of the pDMD can depend on how large are the $r_i$ values used in the algorithm. For this reason, in the next simulations we will visualize: i) for a fixed $N$, the target rank vector $\r \in \RR^N$ as a measure of complexity along the subsets $S_i, i=1, \dots, N$ (that is on the time subintervals of the piecewise technique) and ii) the maximum rank 
\begin{equation}
\label{rank_N}
\tilde{r}(N) = \max_{i = 1, \dots, N} r_i = \| \r \|_\infty
\end{equation}
needed by pDMD by varying the partition size $N$ of the original dataset $S$ until the optimal value identified by the Algorithm \ref{Alg:PiecDMD_N}.


%
\subsection{FitzHugh-Nagumo model}
%
In this section, we apply the pDMD to the subset $S$ generated by the FitzHugh-Nagumo model introduced in Section \ref{ex_FHN}. To start the Algorithm \ref{Alg:PiecDMD_N}, we consider the thresholds $\overline{tol}= 10^{-1},\, tol=10^{-6}$ and $N=1$.

In the left panel of Figure \ref{fhn_err_piecewise}, we show the behaviour of the relative error $\mathcal{E}_p(\widetilde{S}^N,{\bf r})$ with respect to $N$. The plot starts from $N = 17$, because this is the first $N$ value that satisfies the condition $err(i) \leq  \overline{tol}$ for all $i=1, \dots, N$, see step 9 of Algorithm \ref{Alg:PiecDMD}. We note that the value of $\overline{tol}$ in this example is a rather mild request. We opted for this choice to show a more complete history of the error $\mathcal{E}_p(\widetilde{S}^N,{\bf r})$. We note that, by choosing for instance $\overline{tol}=10^{-3},$ $N=87$ would have been the first acceptable value.  
Incrementing the partition size $N$ by one, Algorithm \ref{Alg:PiecDMD_N} stops for $N=147$ with stopping criteria $\mathcal{E}_p(\widetilde{S}^N,{\bf r})<tol$, 
We can observe that the relative error $\mathcal{E}_p(\widetilde{S}^N,{\bf r})$ is almost decreasing with respect to $N$, although there are few little jumps, still remaining in the same order of magnitude. We remark that for $N=147$, we have datasets of dimension $\nu = 41,$ whereas for $N=17$, $\nu = 353$. 

In the right panel of Figure \ref{fhn_err_piecewise}, we show 
how the relative error \eqref{error_time} changes in time for $N = 17$ and $N = 147$, corresponding to the maximum and minimum of $\mathcal{E}_p$ that is $\mathcal{E}_p(\widetilde{S}^{17},{\bf r}) = 0.0124$ and $\mathcal{E}_p(\widetilde{S}^{147},{\bf r}) = 8.1308 \times 10^{-7}$, respectively. In the right plot, we can appreciate that for larger $N$ the error uniformly decreases, especially in the peaks. 


%
\begin{figure}[hbpt]
\centering
 \includegraphics[scale=0.45]{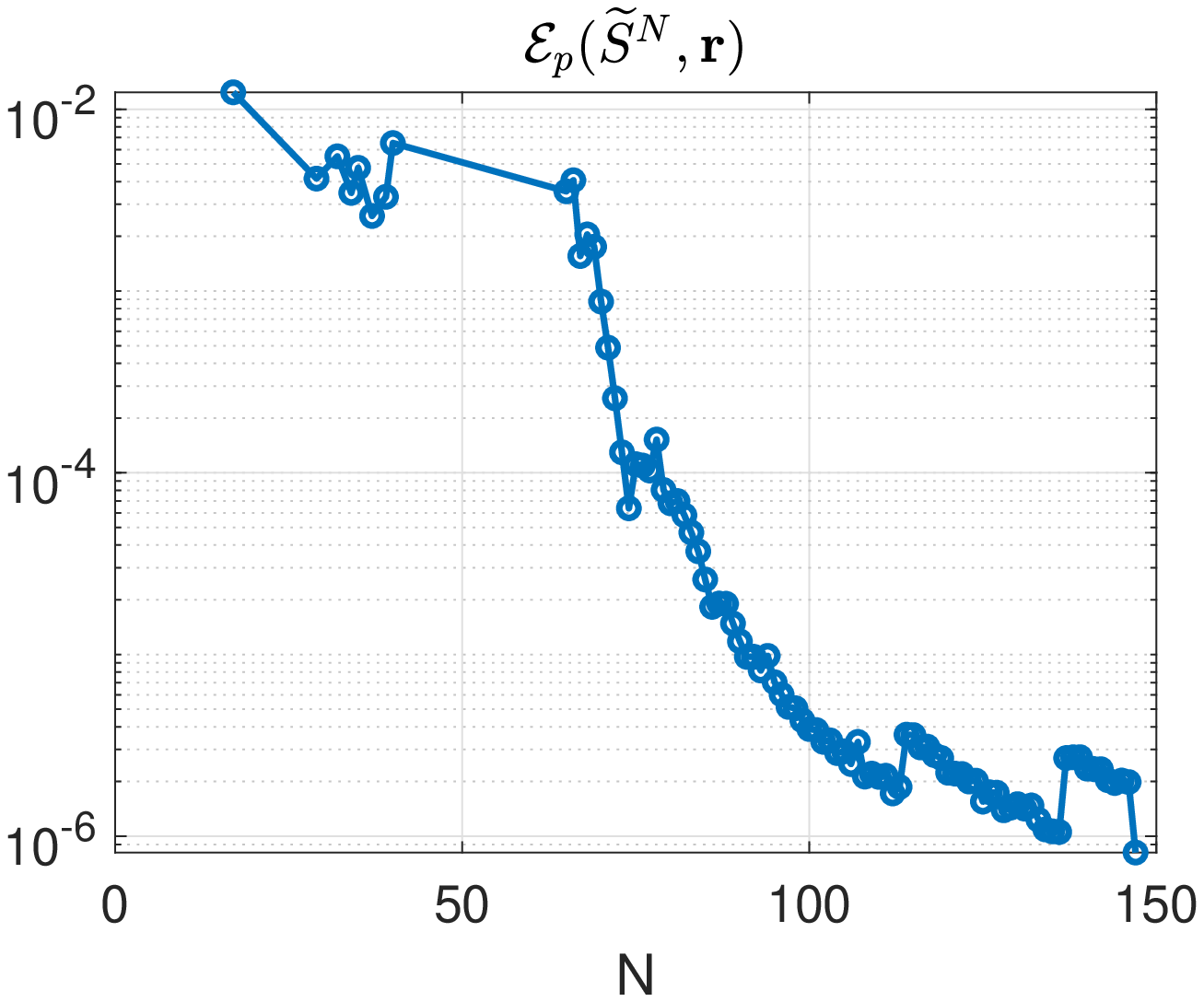}
 \includegraphics[scale=0.45]{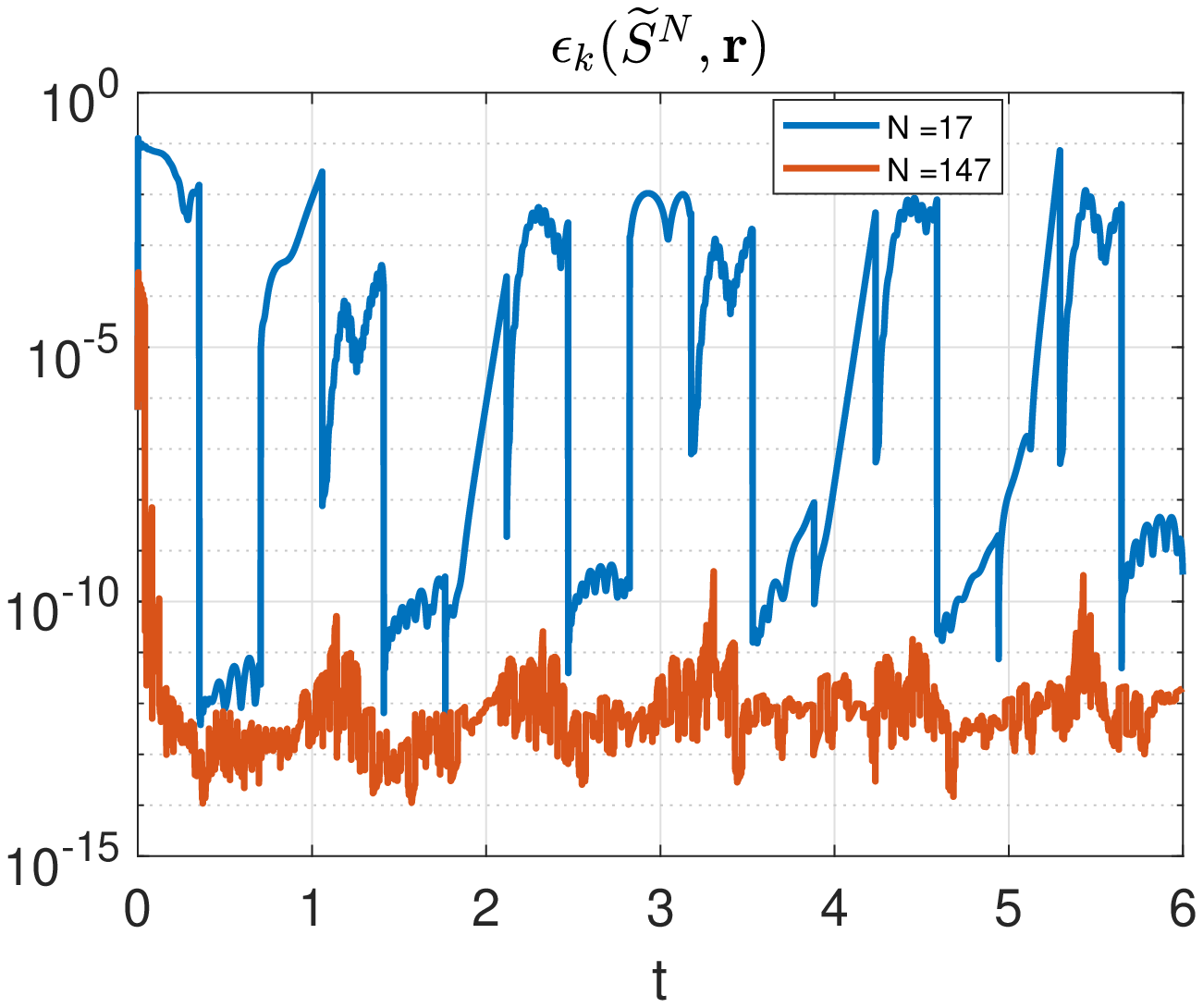}
\captionsetup{justification=justified}
\caption{FHN model. Left plot: relative error $\mathcal{E}_p(\widetilde{S}^N,{\bf r})$ in \eqref{rel_err} for increasing values of the partition size $N$. Right plot: relative error $\epsilon_k(\widetilde{S}^N,{\bf r})$ over time $t$ in \eqref{error_time} for $N$ corresponding to the minimum ($N=147$) and maximum ($N=17$) errors in the left plot.}
\label{fhn_err_piecewise}
\end{figure}
\noindent
Then, we compare the spatial mean \eqref{mean} for the variables $u$ and $v$ obtained by the pDMD reconstruction with $N = 147$, with respect to the data. We choose the value $N = 147$, since it is the value that satisfies the condition $\mathcal{E}_p(\widetilde{S}^N,{\bf r}) < tol$. The results are shown in Figure \ref{fhn_mean_piecewise}, left and middle panels, respectively. It is evident that pDMD carefully matches the data, as also confirmed looking at the phase plane in the right panel. 

%
\begin{figure}[htbp]
\centering
\includegraphics[scale=0.4]{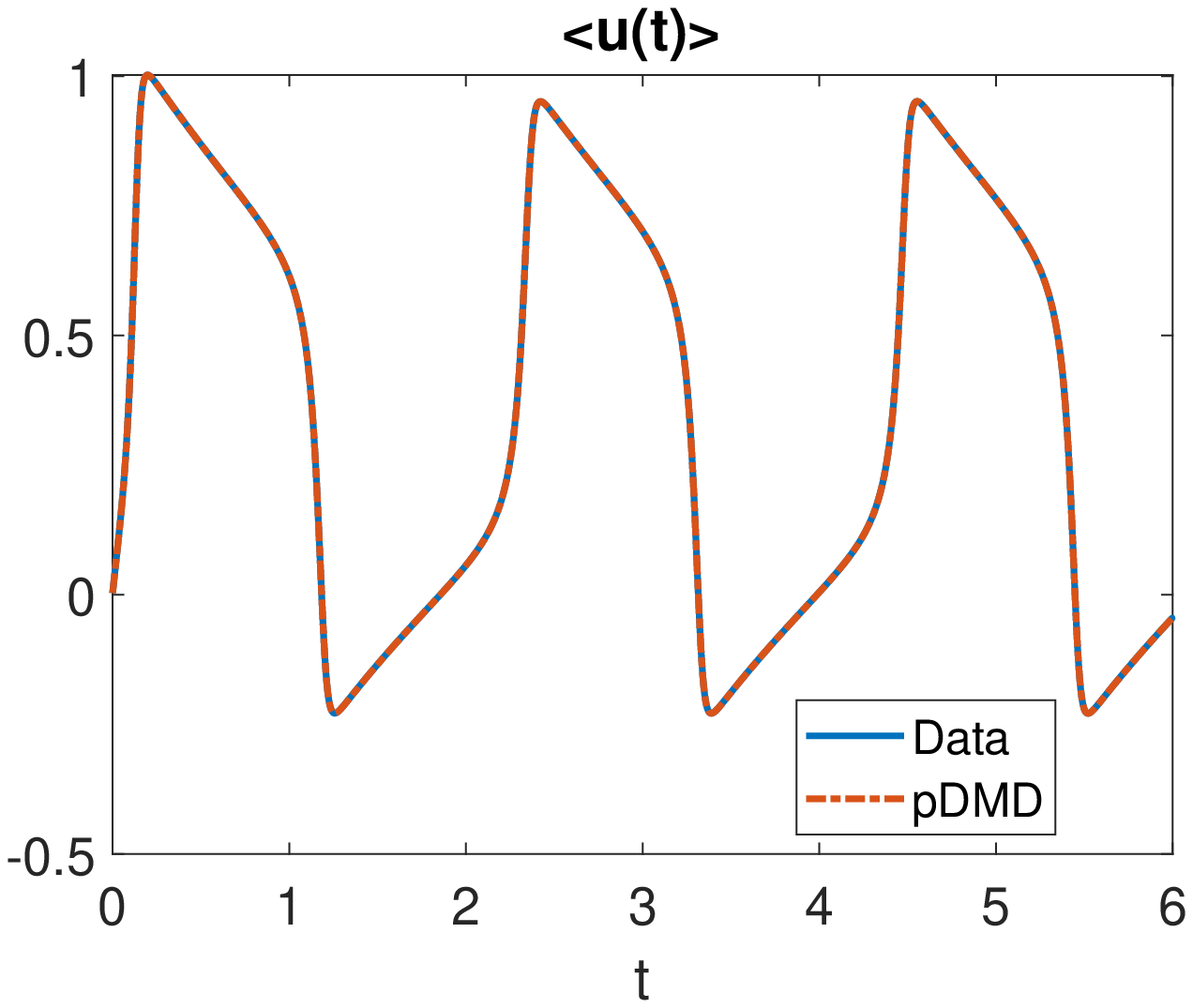}
 \includegraphics[scale=0.4]{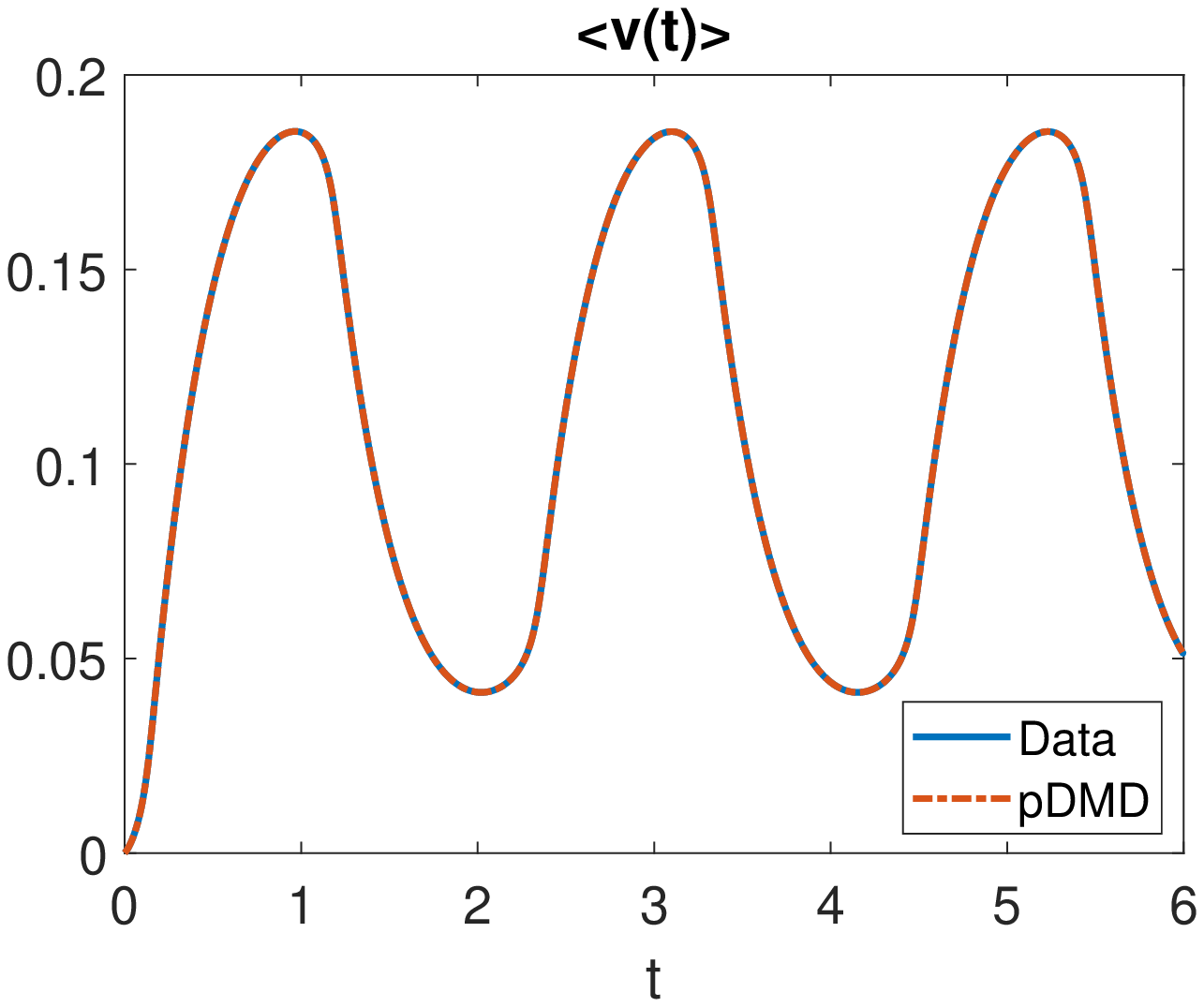}
 \includegraphics[scale=0.4]{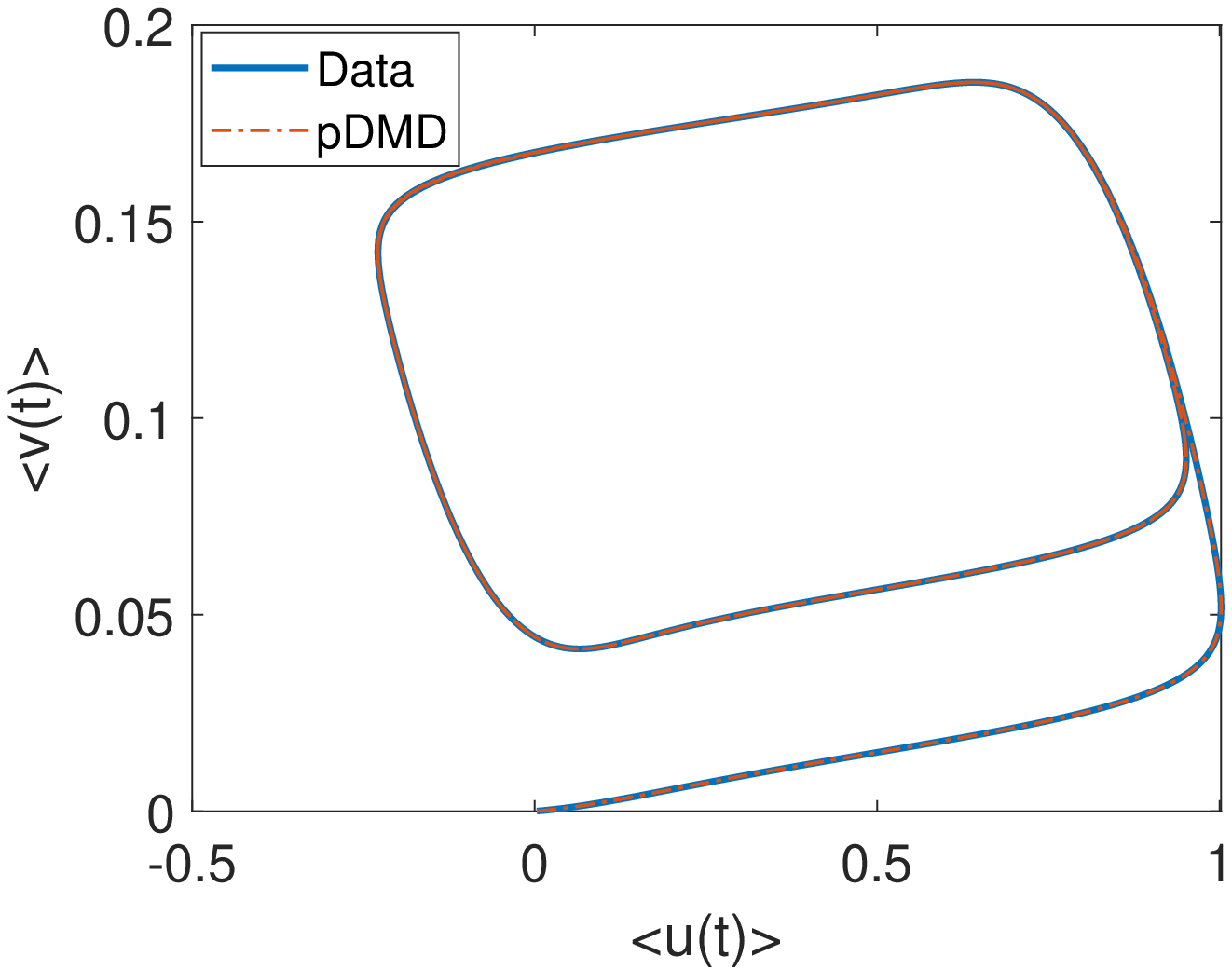}
\captionsetup{justification=justified}
\caption{FHN model. Comparison of the spatial mean for the variable $u$ (left plot) and $v$ (center plot). Right plot: phase plane for the DMD reconstruction with $N = 147$, for which the error has its minimum that is  $\mathcal{E}_p(\widetilde{S}^N,{\bf r}) = 8.1308 \times 10^{-7}$.}
\label{fhn_mean_piecewise}
\end{figure}
\noindent
Finally, in the left plot of Figure \ref{fhn_rank_pdmd}, we show the target rank vector $\r$ for $N = 17, 147$, to show how the pDMD ranks $r_i$ change in each subset $S_i$. It is evident that, the peaks of the rank correspond to the maxima and minima of the spatial mean. Moreover, we show  in the right plot of Figure \ref{fhn_rank_pdmd} the maximum rank $\tilde{r}(N)$, defined in \eqref{rank_N}. We can see that, for all $N$, the maximum rank is always less or equal to $\tilde{r}(17)=29 \ll 51 = {\tt rank}(S)$ the original dataset and for the last $N$ value $\tilde{r}(147) = 12$ holds. This indicates that pDMD is also convenient from the computational point of view since we deal with problems of significant small size.
%
\begin{figure}[htbp]
\centering
 \includegraphics[scale=0.45]{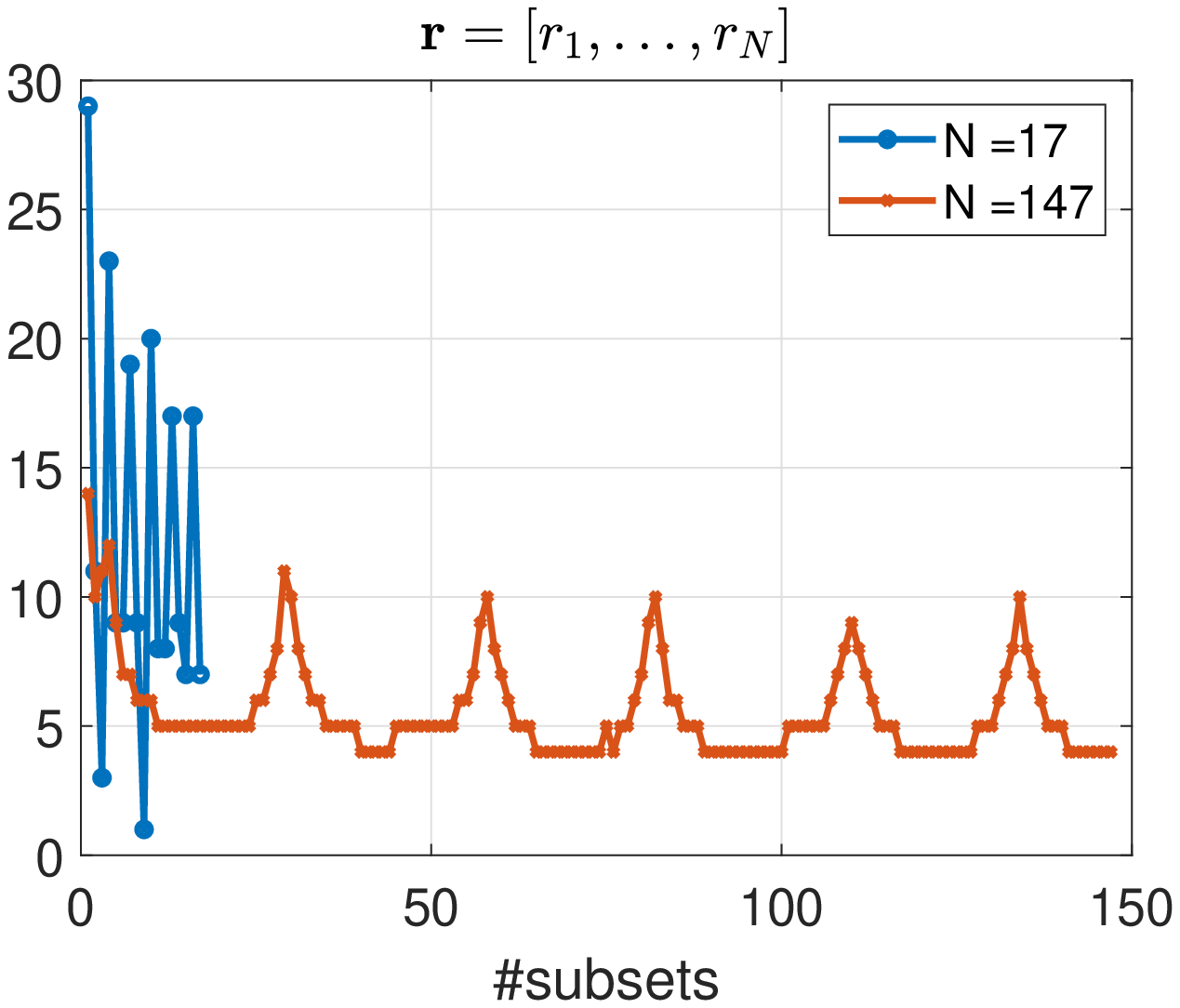}
 \includegraphics[scale=0.45]{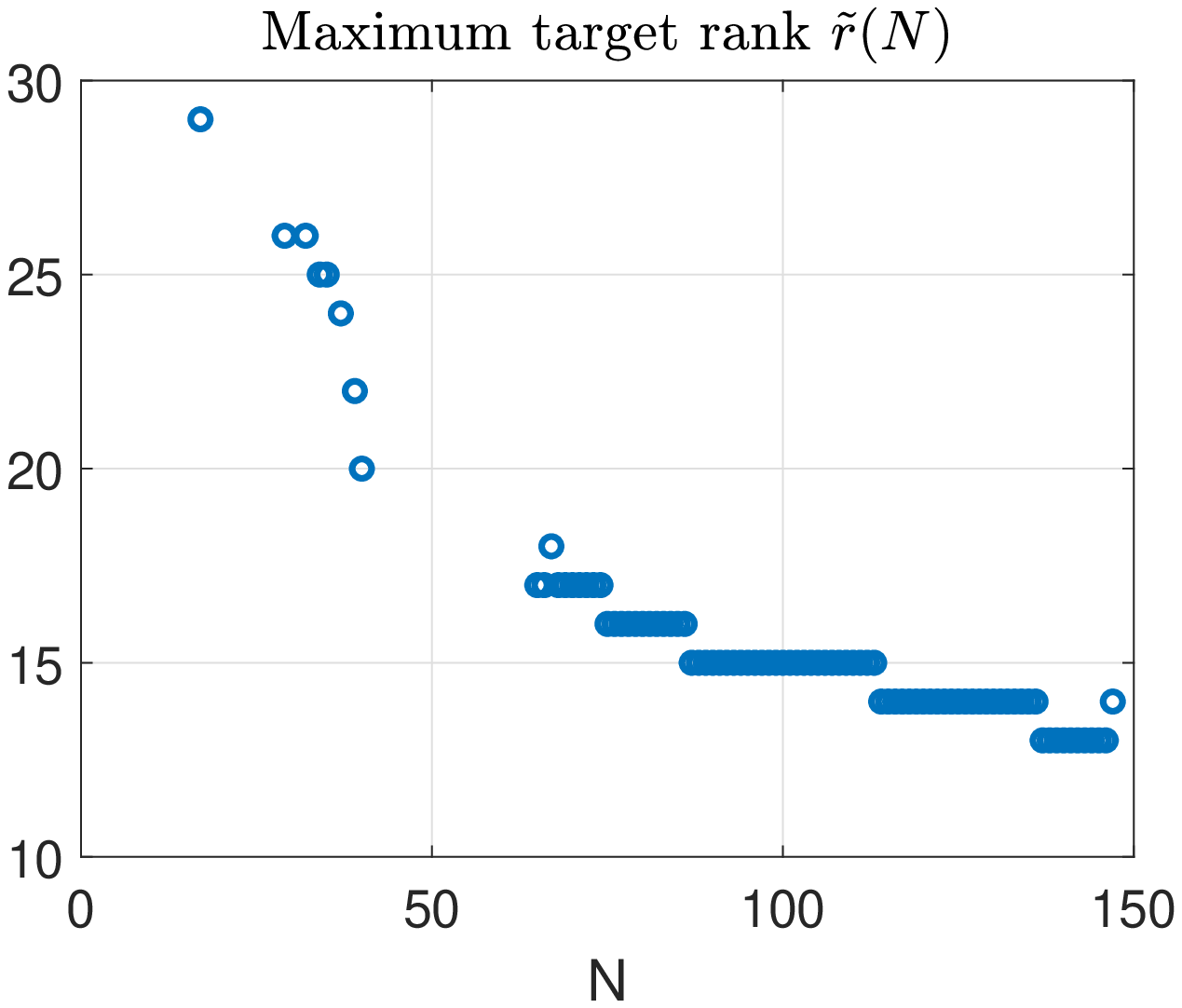}
\captionsetup{justification=justified}
\caption{FHN model. Left plot: vector $\r=[r_1, \dots, r_N]$ of target ranks used on the subsets $S_i$, for $N=17,147$. Right plot: maximum target rank \eqref{rank_N} used by pDMD with respect to $N$.}
\label{fhn_rank_pdmd}
\end{figure}
%
%
%
%
%
\subsection{$\lambda$-$\omega$ RD-PDE system}
We apply pDMD to the same dataset $S$ generated for the $\lambda$-$\omega$ system discussed in Section \ref{ex_lo} and we consider $\overline{tol}=10^{-2}, N=1$ to start the Algorithm \ref{Alg:PiecDMD}. In Figure \ref{lo_err_piecewise} (left panel) is shown the error $\mathcal{E}_p(\widetilde{S}^N,{\bf r})$ defined in \eqref{rel_err} when the partition size $N$, i.e. number of submatrices $S_i$ of $S$, is increased. For $N = 1$, corresponding to the "global" DMD case, as expected from our results in Section \ref{ex_lo}, the condition in step $9$ of the algorithm is not satisfied and the first acceptable value is $N = 16$ which corresponds to $\nu = 1042$.

We stop Algorithm \ref{Alg:PiecDMD_N} when $tol=10^{-6}$, for $N=48$, that is $\nu=261$. 
We can see that the error exhibits a sligthly oscillating behaviour for $16 \leq N \leq 47$, such that $10^{-5} \leq \mathcal{E}_p(\widetilde{S}^N,{\bf r}) \leq 10^{-6}$, that is the error remains within the same order of magnitude.
As for the previous test, we consider the values of $N$ for which the error has its maximum and minimum, that are $\mathcal{E}_p(\widetilde{S}^{16},{\bf r}) = 9.4858 \times 10^{-6}$ and $\mathcal{E}_p(\widetilde{S}^{48},{\bf r}) = 1.8478 \times 10^{-7}$. In any case, these error approximation levels are much lower than the best obtained by the ``global'' DMD in Section \ref{ex_lo}. In Figure \ref{lo_err_piecewise}, right panel, we show the errors in time \eqref{error_time} for $N = 16$ and $N = 48$. We can observe in the error behaviour two time regimes, corresponding to those of the spiral wave dynamics, discussed before in Section \ref{ex_lo}. In fact, in both cases, the error rapidly decays immediately after the initial phase until $\bar{t} \approx 25$, then an almost constant oscillating trend is present in $[\bar{t},T]$. We note also that for larger $N$ the error uniformly decreases along all the interval $[0,T]$.
 

%
\begin{figure}[htbp]
\centering
\includegraphics[scale=0.45]{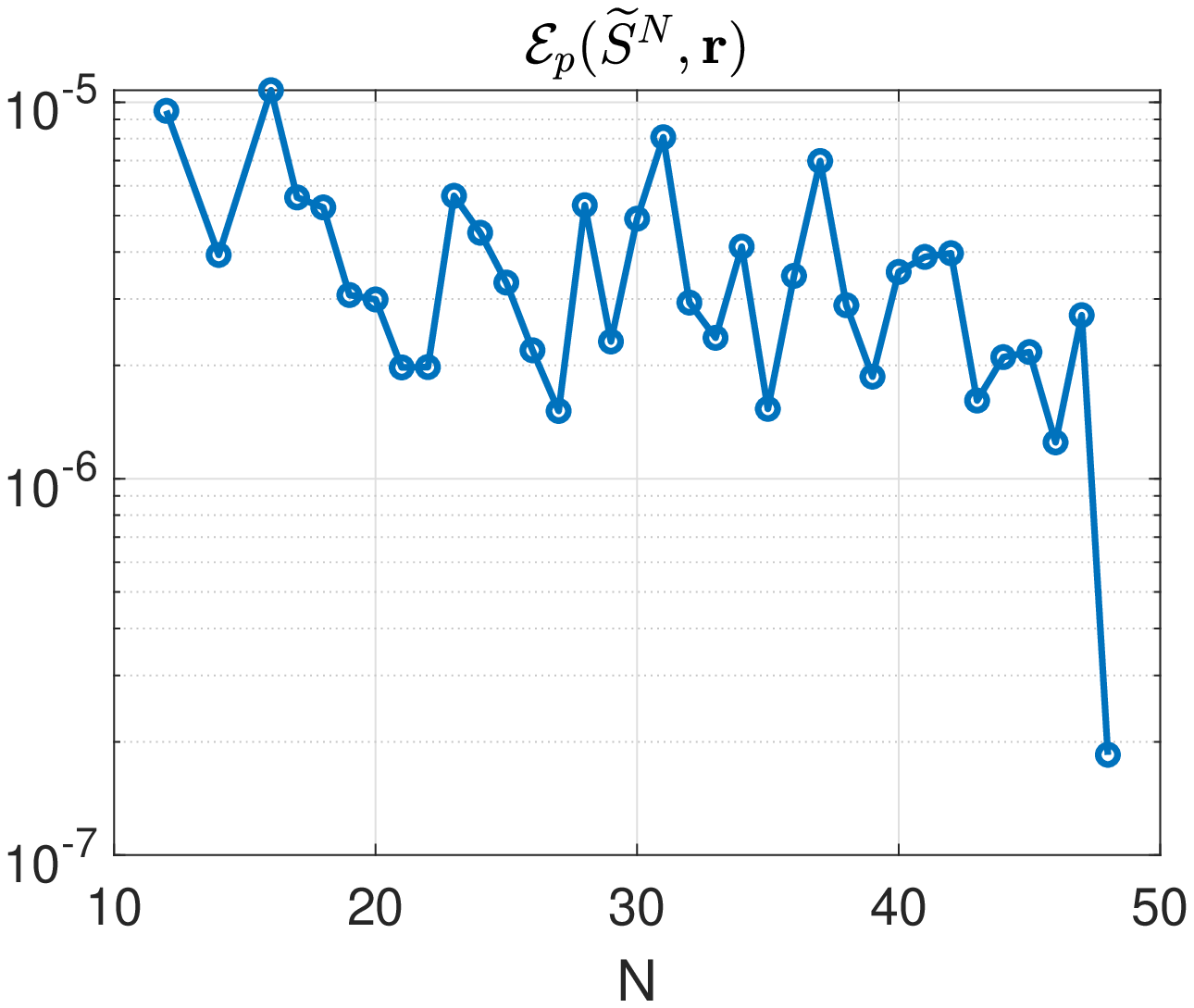}
\includegraphics[scale=0.45]{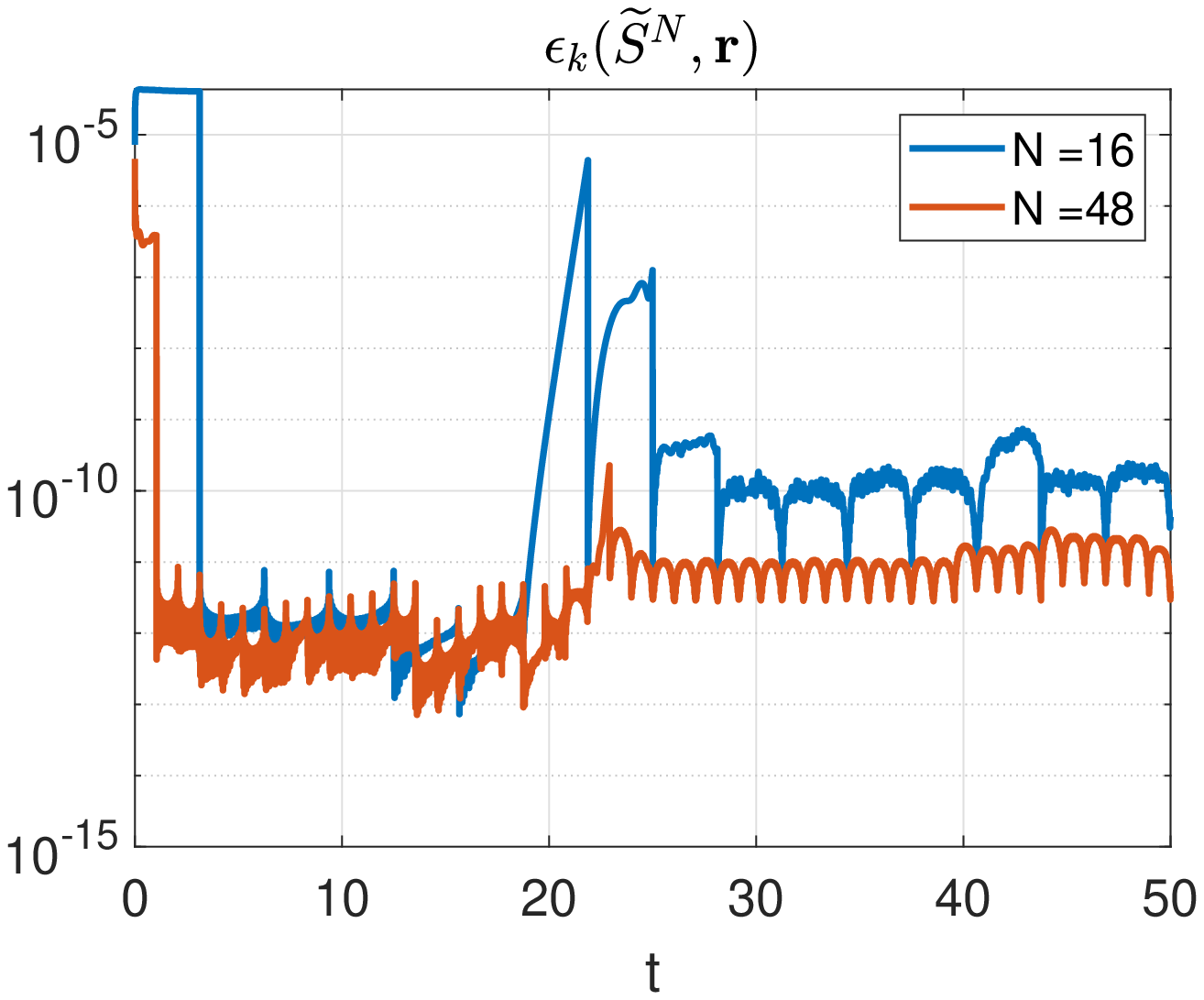}
\captionsetup{justification=justified}
\caption{$\lambda$-$\omega$ system, spiral wave. Left plot: relative error $\mathcal{E}_p(\widetilde{S}^N,{\bf r})$ for $
\overline{tol}=10^{-2}$ in \eqref{rel_err}. Relative error $\epsilon_k(\widetilde{S}^N,{\bf r})$ over time in \eqref{error_time} (right plot) for two significant values of $N$ as discussed in the main text.}
\label{lo_err_piecewise}
\end{figure}
\noindent
To further show how the pDMD overcomes the drawbacks of the original DMD highlighted in Section \ref{ex_lo}, we compare the time dynamics of the spatial mean for $N=48$ 
with respect to the data. Figure \ref{lo_mean_piecewise} shows these comparisons for $u$ (left plot) and $v$ (middle plot). As for the FHN model, the oscillating time dynamics obtained by pDMD matches perfectly the dataset, as also confirmed looking at the reconstructed limit cycle in the right panel. 

\begin{figure}[t]
\centering
 \includegraphics[scale=0.4]{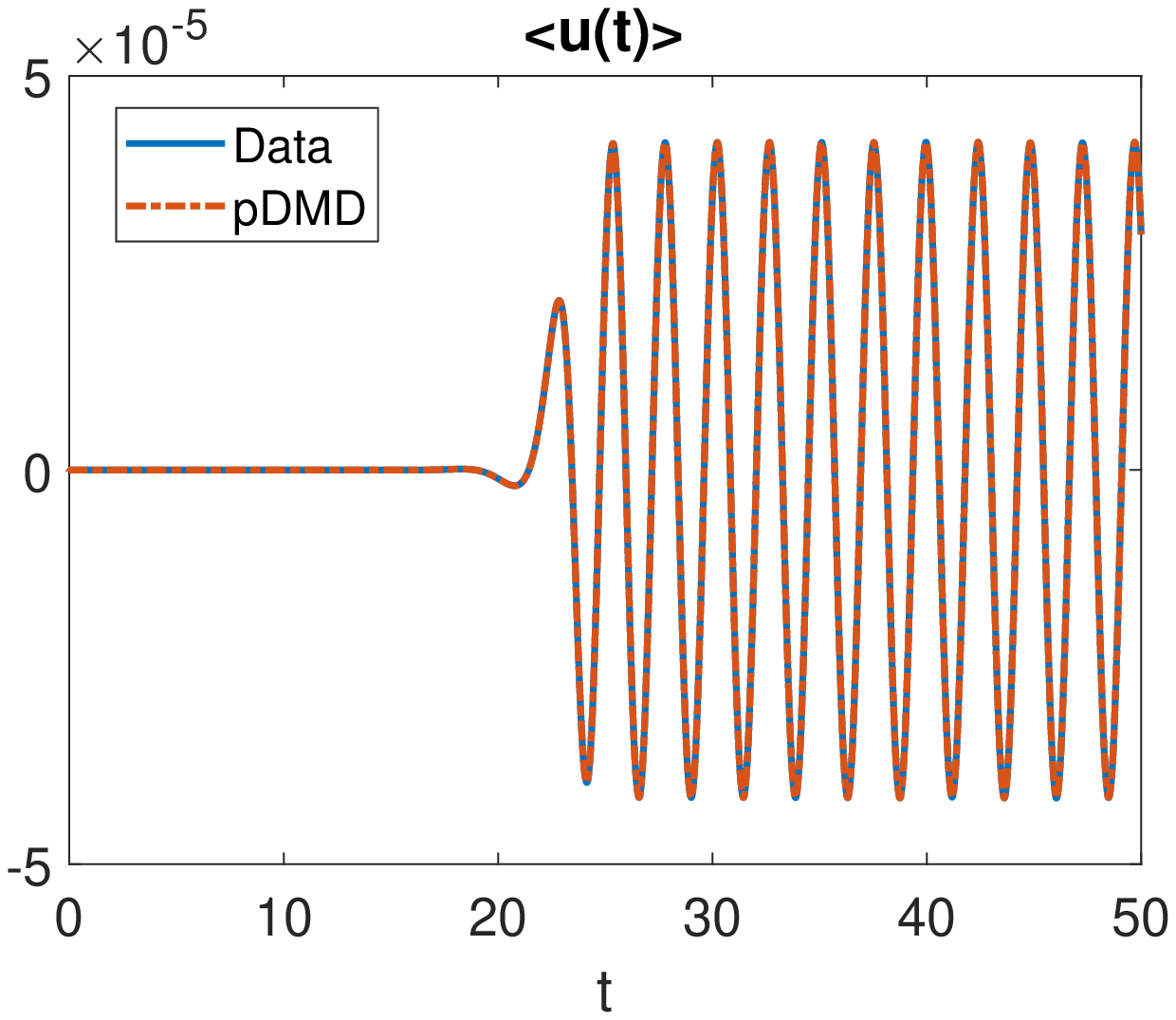}
 \includegraphics[scale=0.4]{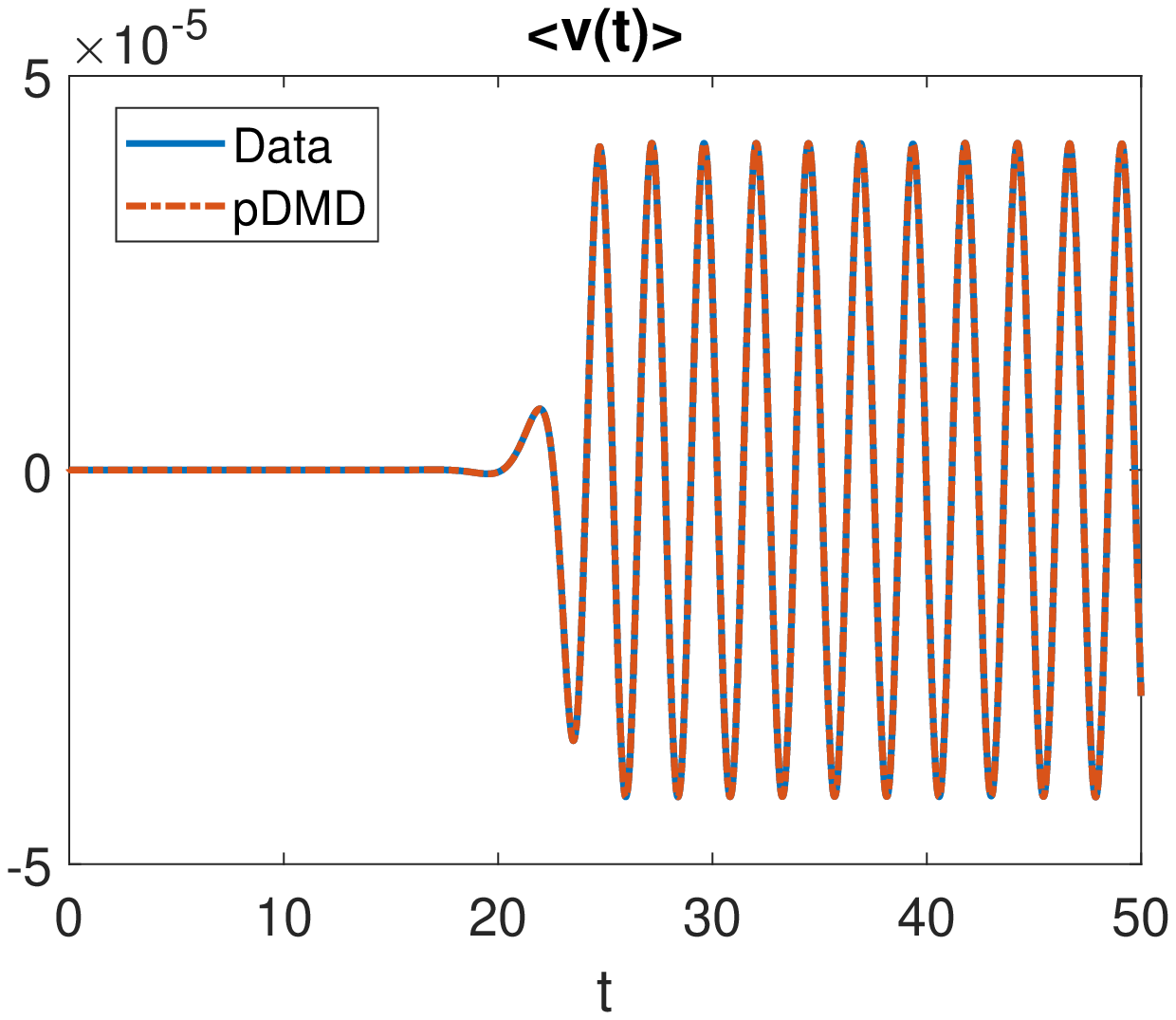}
\includegraphics[scale=0.4]{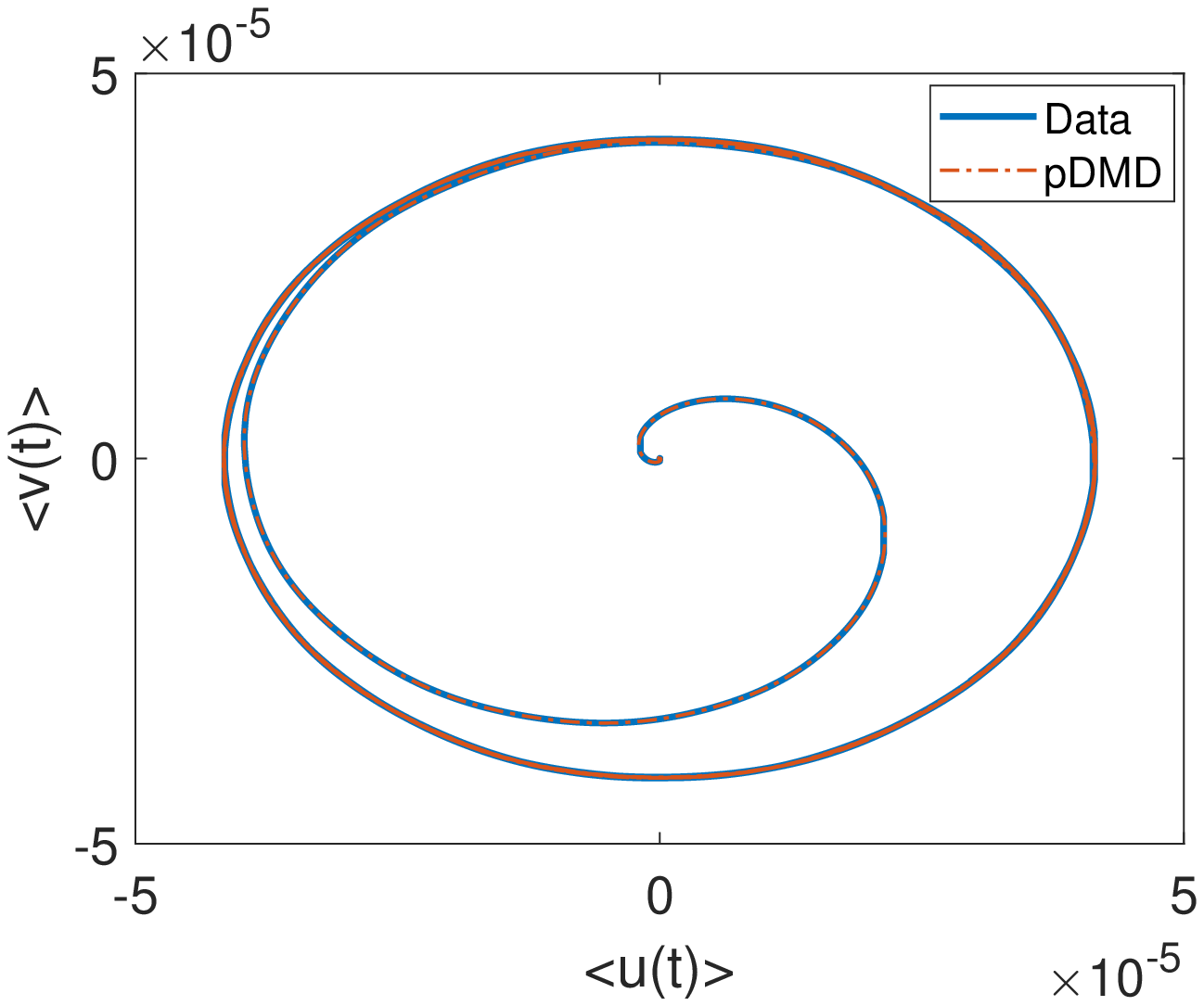}
\captionsetup{justification=justified}
\caption{$\lambda$-$\omega$ system, spiral waves. pDMD reconstruction for the spatial mean for the variable $u$ (left plot) and $v$ (right plot) with $N = 48$, for which the error has its minimum that is $\mathcal{E}_p(\widetilde{S}^N,{\bf r}) = 1.8478 \times 10^{-7}$. Corresponding limit cycle in the phase plane is shown in the right plot, a very good agreement is obtained with respect to the best approximation by the original DMD in Section \ref{ex_lo}.}
\label{lo_mean_piecewise}
\end{figure}
\noindent
Finally, in the left panel of Figure \ref{lo_rank_pdmd}, we show the target rank vector $\r$ for $N = 16, 48$, to show how the pDMD ranks $r_i$ changes along the subsets $S_i$.
We note that, for each $N$, the largest value of the rank is always required in the first subset, but in general lower $r_i$ values are required for larger $N$. This result is also confirmed by the behaviour of the maximum $\tilde{r}(N)$,  shown in Figure \ref{lo_rank_pdmd}, right plot, that monotonically decays with respect to $N$.
Moreover, it is worth noting that the maximum rank needed by pDMD is $\tilde{r}(16)= 56$ which is much smaller than the rank of the dataset $S$, that is $314$.

\begin{figure}[t]
\centering
\includegraphics[scale=0.45]{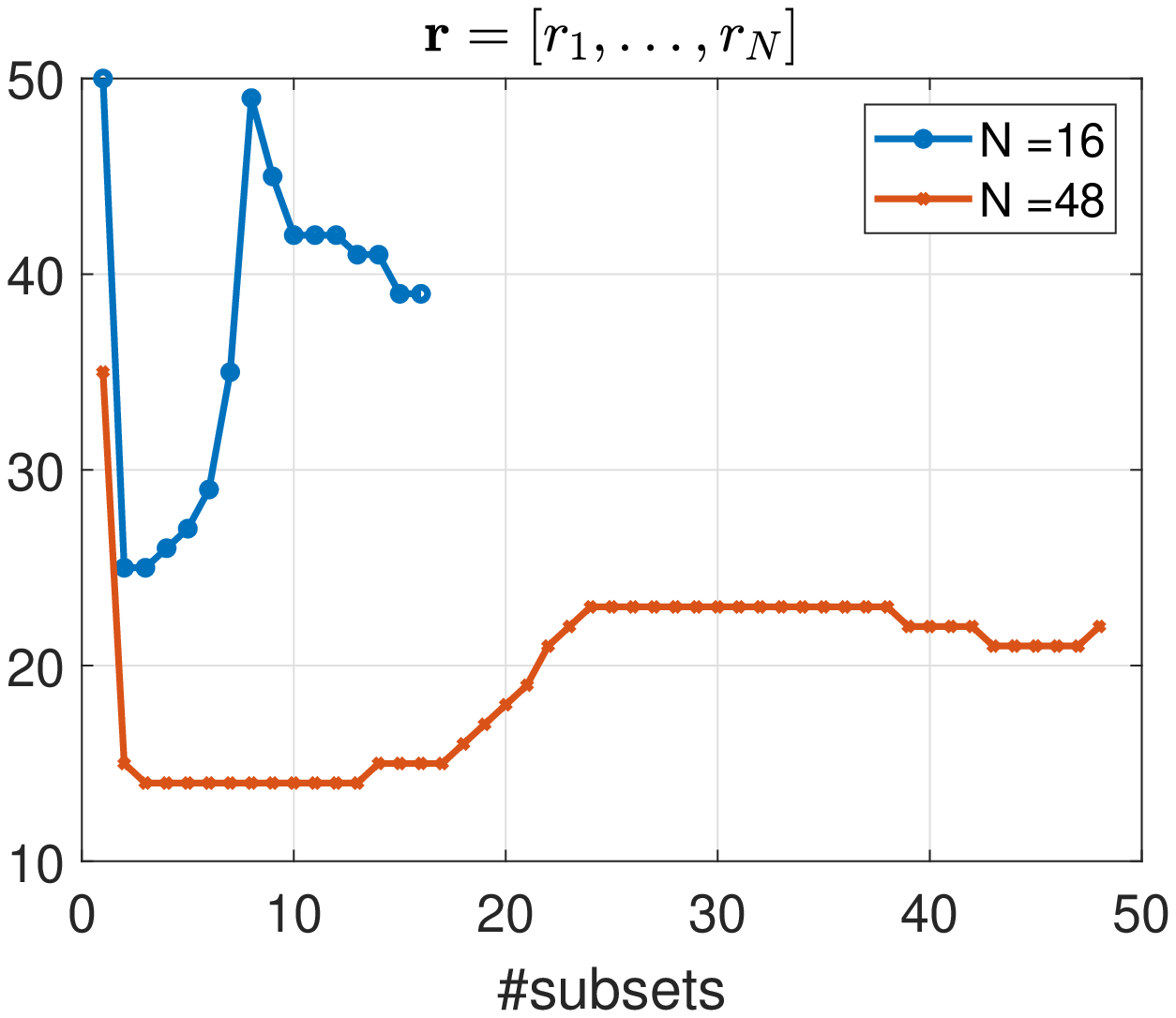}
\includegraphics[scale=0.45]{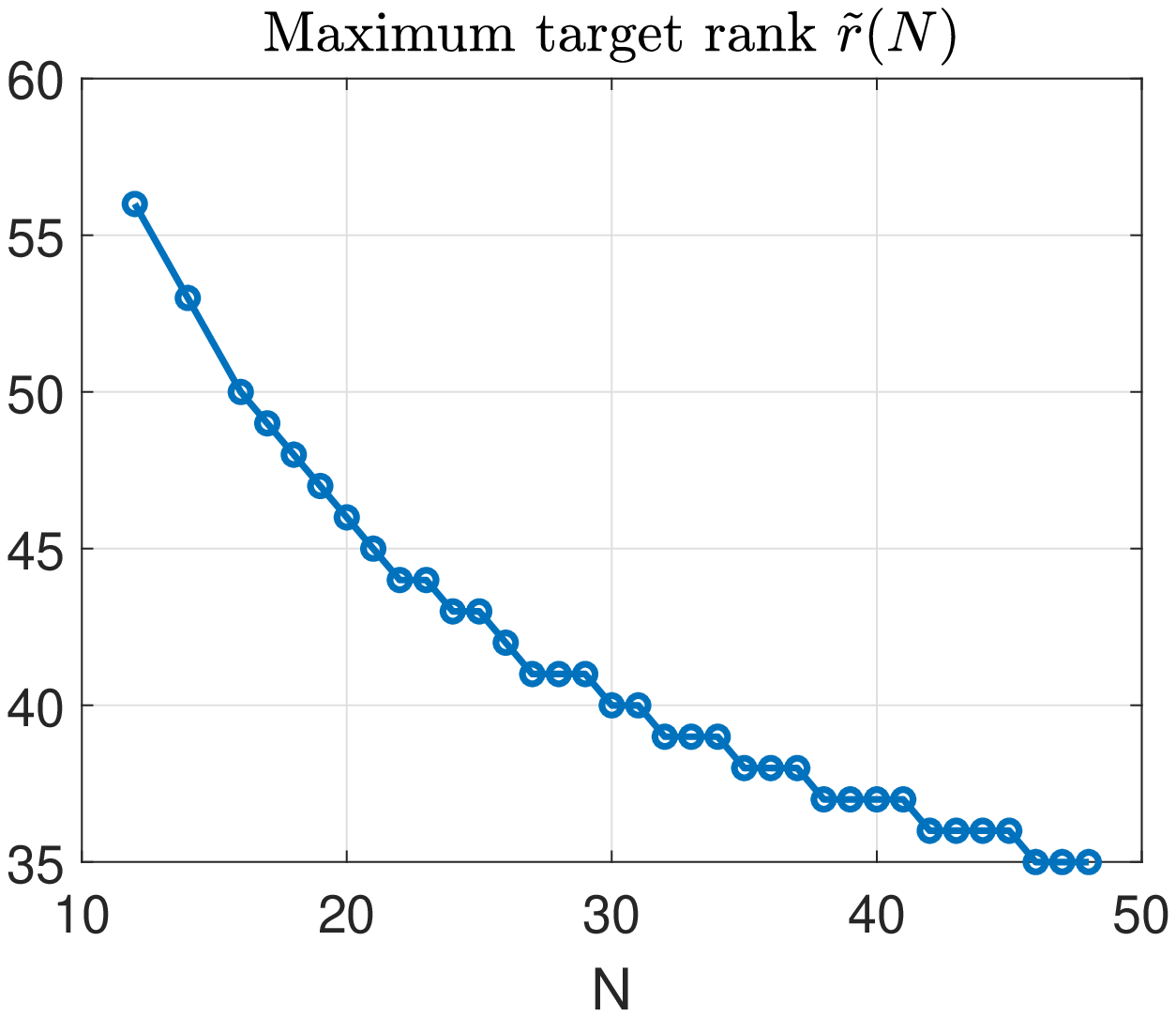}
\captionsetup{justification=justified}
\caption{$\lambda$-$\omega$ system, spiral waves. Left plot: vector $\r=[r_1, \dots, r_N]$ of target ranks used in the subinterval $S_i$, by pDMD for $N=16,48$. Right plot: maximum target rank $\tilde{r}(N)$ in \eqref{rank_N} used by pDMD with respect to $N$.} 
\label{lo_rank_pdmd}
\end{figure}
\newpage
\subsection{DIB model: Turing instability}
%
In this section, we apply the pDMD to reconstruct the Turing instability dynamics of the morphochemical DIB model discussed in Section \ref{sec:turing}. We consider $\overline{tol}=10^{-3}$ and $N=1$ as input of the pDMD Algorithm \ref{Alg:PiecDMD}. 

In the left panel of Figure \ref{dib_err_piecewise} we show the relative error $\mathcal{E}_p(\widetilde{S}^N,{\bf r})$ defined in \eqref{rel_err} for increasing values of $N$ until it is less than $tol = 10^{-6}$ in Algorithm \ref{Alg:PiecDMD_N}, that here happens for $N = 48$.
The first acceptable value is $N=29$, but for $N < 48$ some symbols are missing because for that $N$ the criterion $err(i) < \overline{tol}$ is not satisfied for all $i$. We remark that for $N=29$, we have datasets of dimension $\nu = 345$ whereas with $N=48$, $\nu = 209$ holds. 


In the right panel of Figure \ref{dib_err_piecewise}, we show how the pDMD error \eqref{error_time} evolves in time for $N = 29$ and $N = 48$, corresponding to the maximum and minimum value attained in Figure \ref{dib_err_piecewise}, left plot, given by $\mathcal{E}_p(\widetilde{S}^{29},{\bf r}) = 4.2932 \times 10^{-6}$ and $\mathcal{E}_p(\widetilde{S}^{48},{\bf r}) = 1.4039 \times 10^{-7}$, respectively.
We note that the maximum error is essentially concentrated in the first zone, that is in the reactivity Turing regime, and that for larger $N$ it decreases almost uniformly with respect to time along $[0,T]$.
 
\begin{figure}[htbp]
\centering
\includegraphics[scale=0.45]{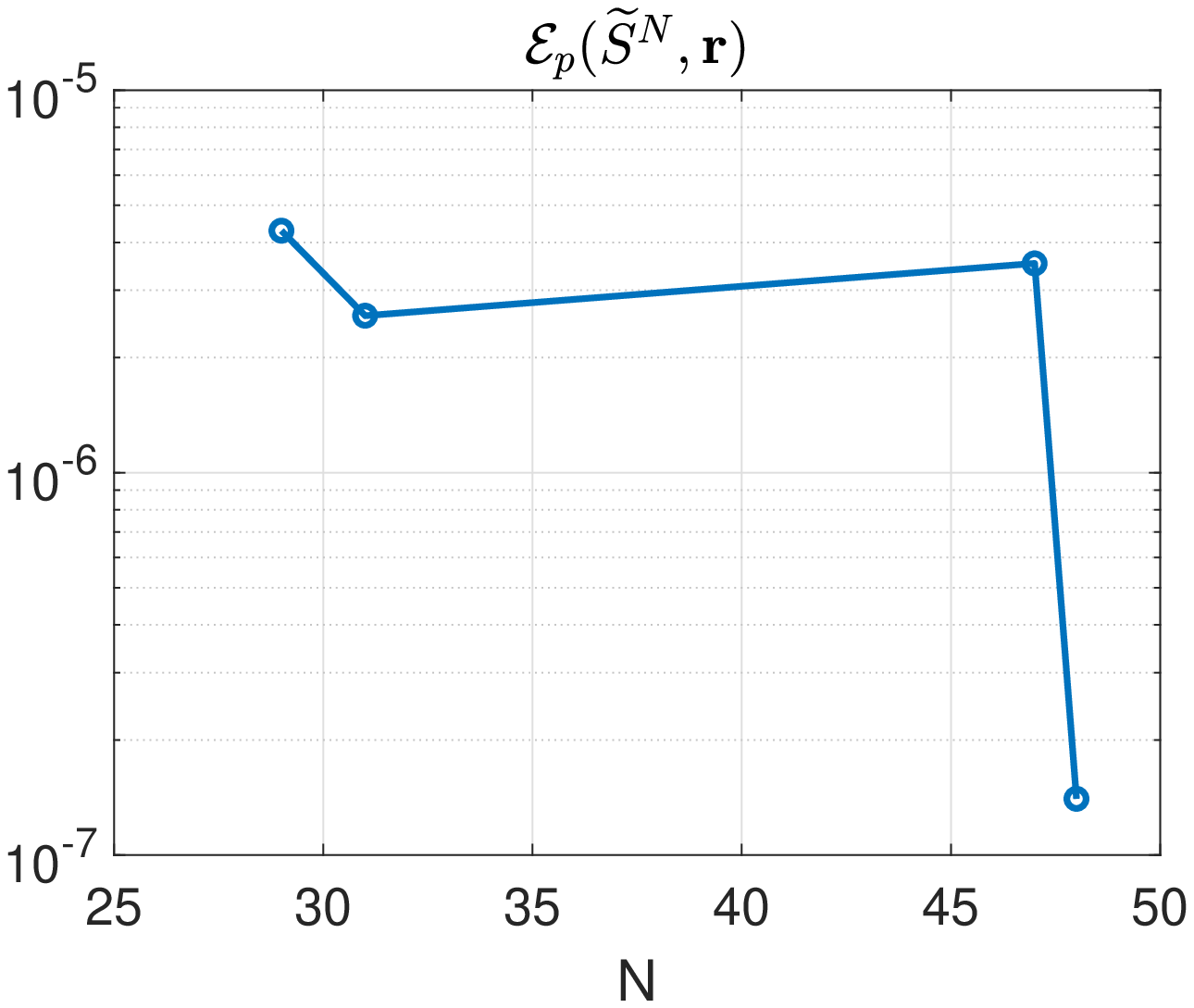}
\includegraphics[scale=0.45]{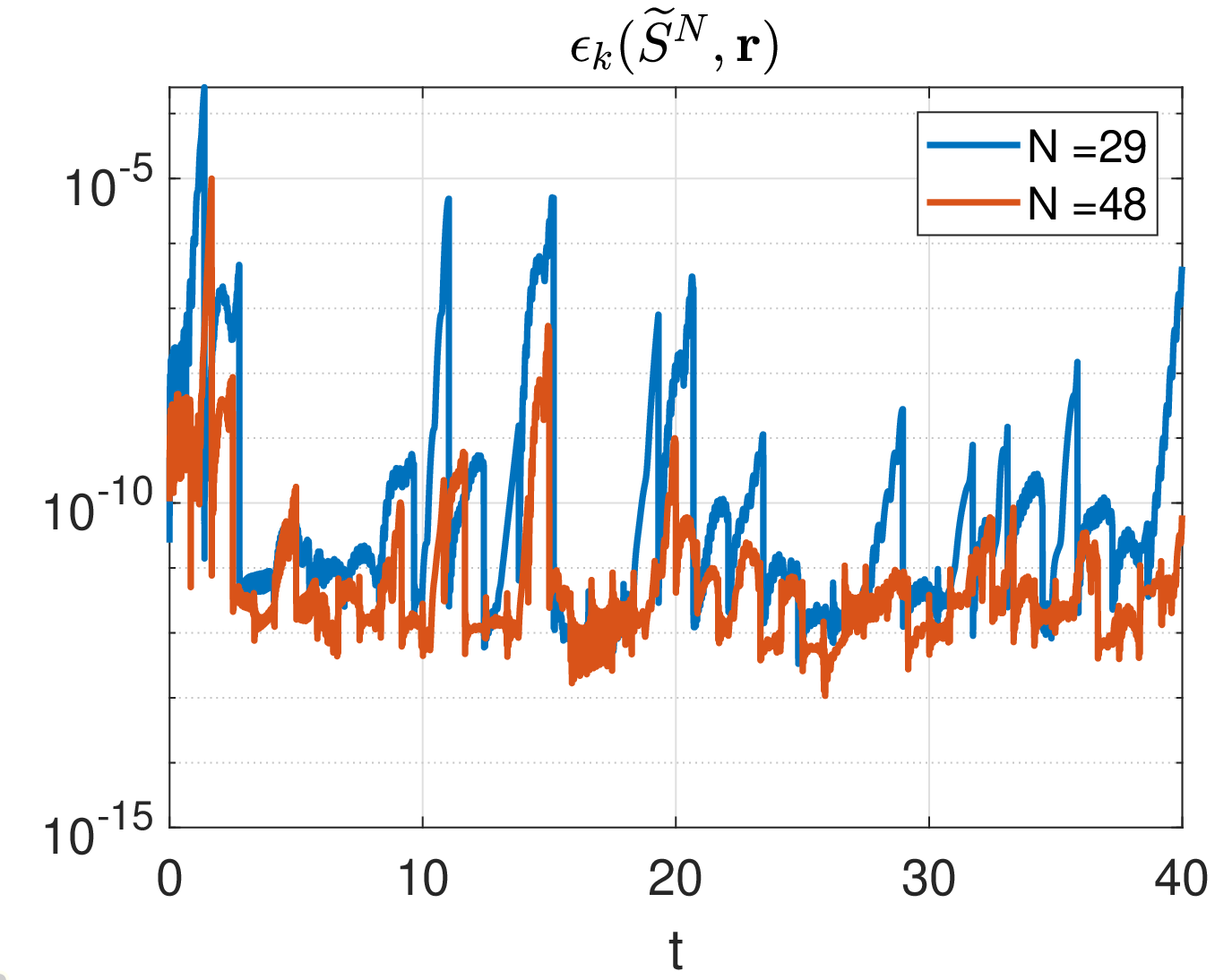}
\captionsetup{justification=justified}
\caption{DIB model: Turing instability. Left plot: relative error $\mathcal{E}_p(\widetilde{S}^N,{\bf r})$ in \eqref{rel_err}. Relative error $\epsilon_k(\widetilde{S}^N,{\bf r})$ over time in \eqref{error_time} for two meaningful values of $N$, as discussed in the main text.}
\label{dib_err_piecewise}
\end{figure}
\noindent
To further confirm this trend, in Figure \ref{dib_mean_piecewise}, we show the comparison of the spatial mean for the variables $u$ (left panel) and $v$ (right panel) obtained by the pDMD reconstruction with $N = 48$ 
with respect to the data. The time dynamics of the pDMD solution for all times matches the spatial means of the data (compare with the right plot in Figure \ref{dib_mean_classic}).

\begin{figure}[htbp]
\centering
\includegraphics[scale=0.45]{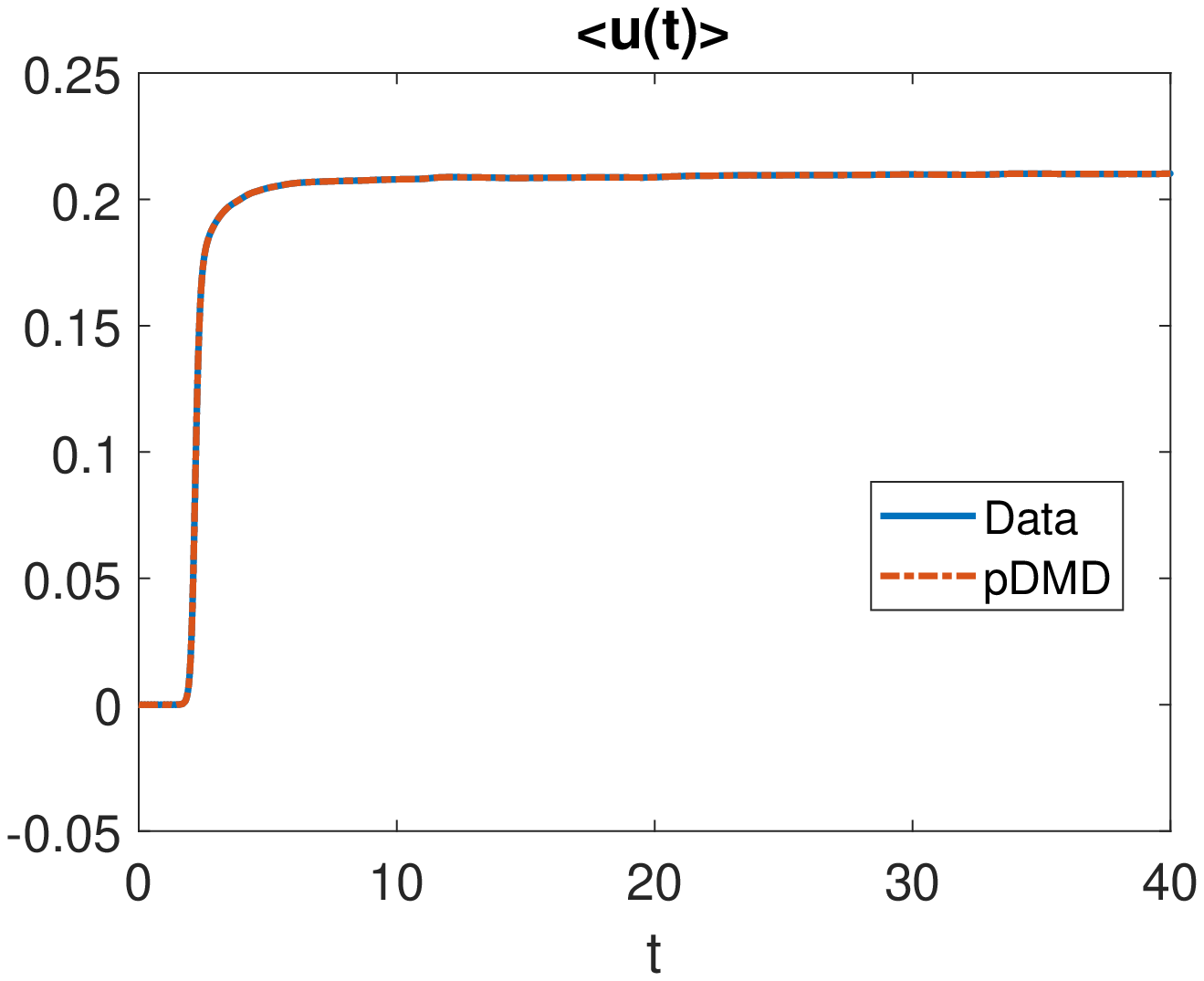}
\includegraphics[scale=0.45]{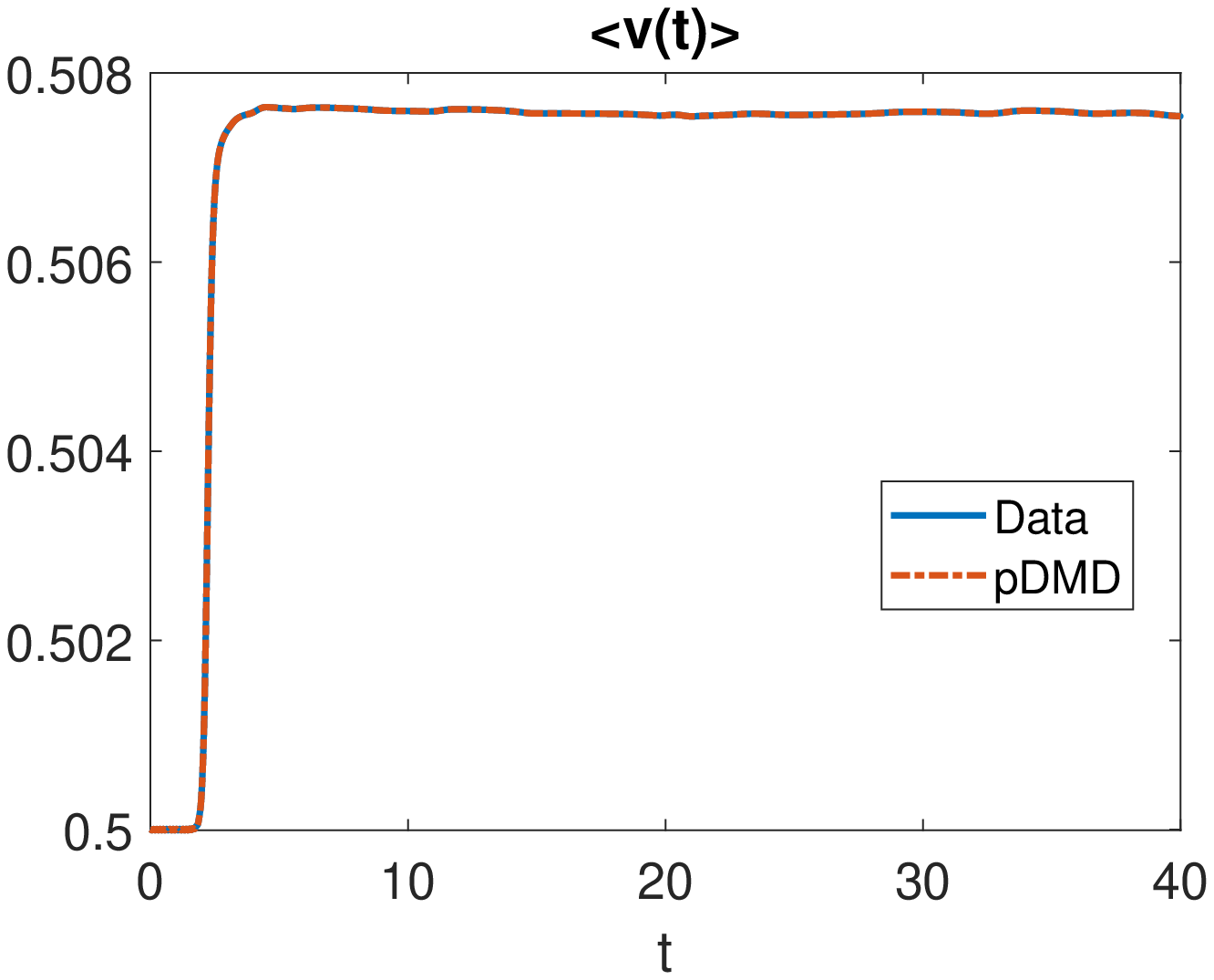}
\captionsetup{justification=justified}
\caption{DIB model: Turing instability. Spatial mean of the variables $u$ (left plot) and $v$ (right plot) for the pDMD reconstruction with $N = 48$ compared with that of the data.}
\label{dib_mean_piecewise}
\end{figure}
\noindent
We emphasize that, for the reconstruction of the Turing instability dynamics the piecewise approach not only is able to remove the ill-conditioning in the ``global'' DMD (see Figure \ref{dib_mean_classic}, left), but also it is able to adapt the choice of the ``local'' target ranks to the peculiar dynamics along time. For this reason, we report how the rank $r_i$ in each subset $S_i$ changes for $i=1, \dots, N$ for the previous partitions with $N = 29$ and $N = 48$.

The left plot in Figure \ref{dib_rank_pdmd} highlights that, for both $N$, in the initial reactivity zone we need to choose higher values of the rank, whereas in the stabilizing zone significantly lower values are required, such that $r_i \leq 20$.\\
In the right plot of Figure \ref{dib_rank_pdmd}, we report the maximum target ranks $\tilde{r}(N) $ for the values of $N$ in Figure \ref{dib_err_piecewise}(left) (that are then attained in the initial part of the time interval). We observe that for all $N$, $\tilde{r}(N) \leq \tilde{r}(29)=43$, that is again much smaller than the rank of the original dataset $S$ given by $303$. This result confirms that the piecewise approach is also able to reduce the overall computational load of the usual DMD implementation.

\begin{figure}[htbp]
\centering
\includegraphics[scale=0.45]{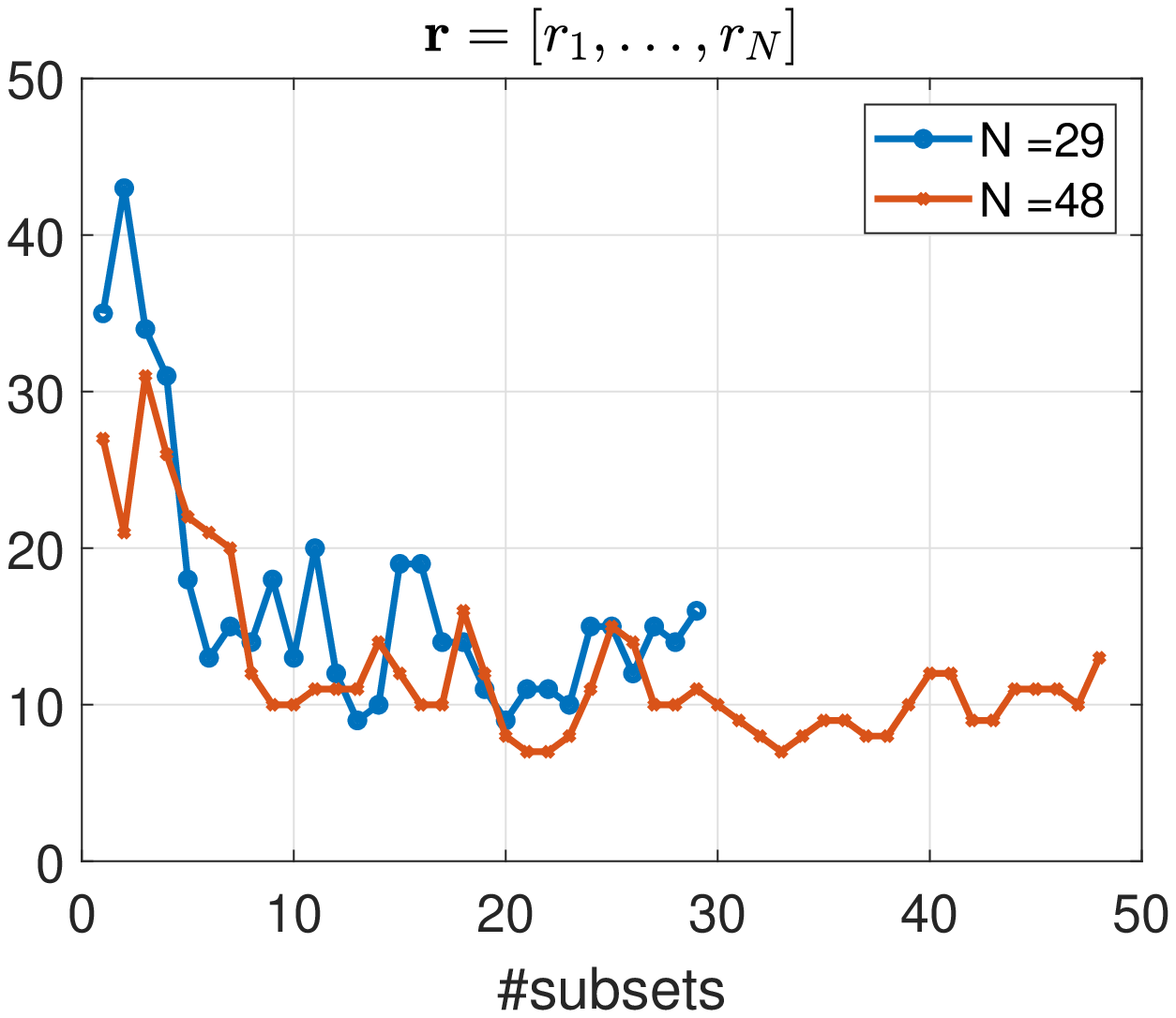}
\includegraphics[scale=0.45]{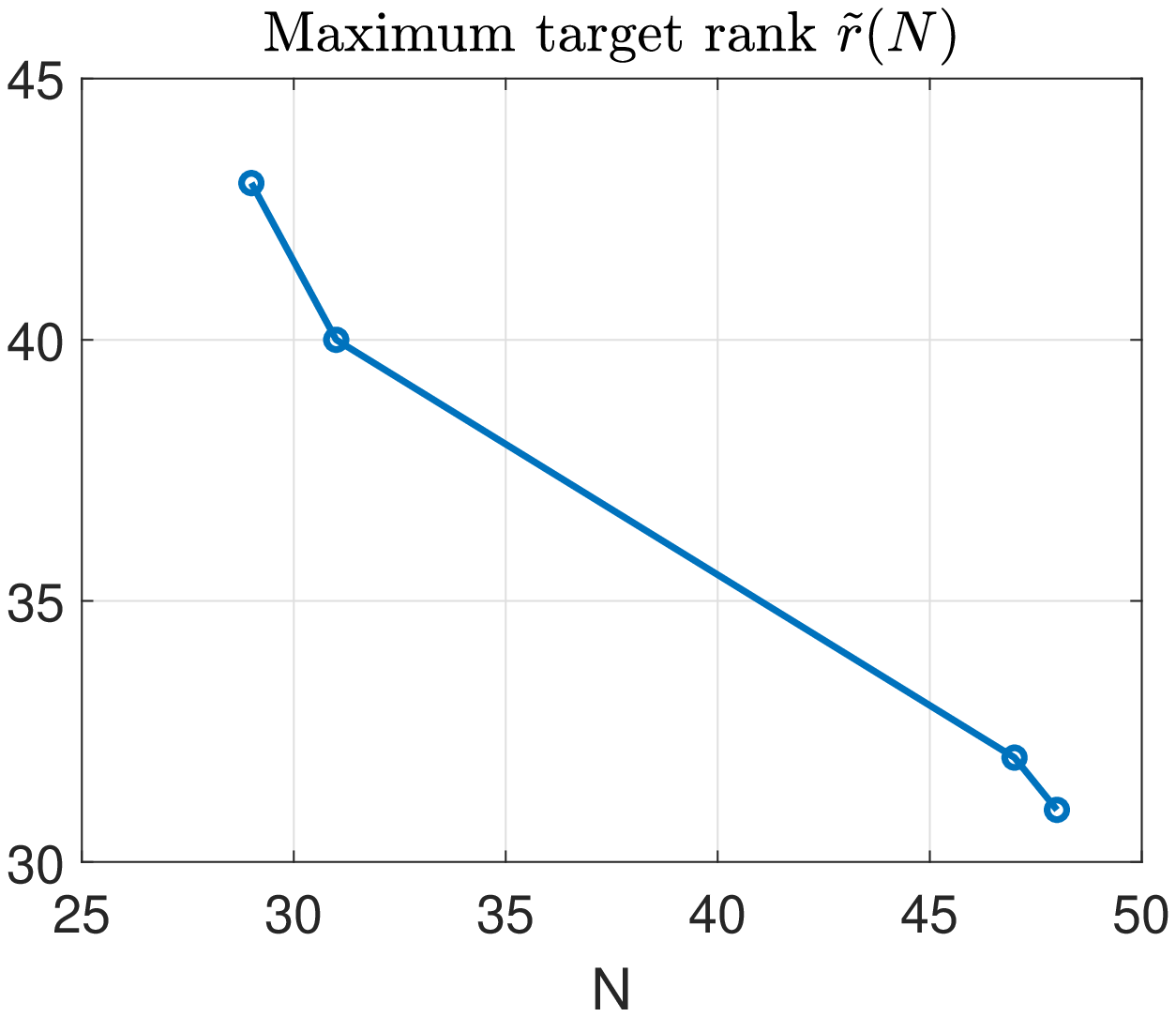}
\captionsetup{justification=justified}
\caption{DIB model: Turing instability. Left plot: vector $\r=[r_1, \dots, r_N]$ of target ranks used in the subsets $S_i$, for $N=29, 48$. Right plot: maximum target rank $\tilde{r}(N)$ in \eqref{rank_N} used by pDMD with respect to $N$.}
\label{dib_rank_pdmd}
\end{figure}
%
%
\subsection{DIB model: Turing-Hopf instability}
The last experiment concerns again the DIB morphochemical RD system but with a different choice of the model parameters that gives rise to the so-called Turing-Hopf instability. We recall from Section \ref{sec_fail_dmd}, that this is the most complicated dynamics considered in this paper, because the PDE solutions exhibit an initial instability and then the formation of a pattern oscillating both in space and time. We apply the pDMD to the same dataset $S$ generated in Section \ref{sec_fail_dmd} starting the Algorithm \ref{Alg:PiecDMD} with $N=5$ and $\overline{tol}=10^{-3}$. We obtain the first useful partition for $N = 225$, that is $\nu =50$. Then, we increment $N$ by $5$ and check the error \eqref{rel_err} until $tol =10^{-5}$ in Algorithm \ref{Alg:PiecDMD_N}. The final value is $N=535$ with $\nu = 21$ where $\mathcal{E}_p(\widetilde{S}^{535},{\bf r}) = 8.5031 \times 10^{-6}$}. In Figure \ref{dib_hopf_err_piecewise} (left plot) we show that the error \eqref{rel_err} slowly decreases for increasing $N$, even though there are small oscillations in a neighborhood of $N = 500$. In this case, it is clear that the dynamics is very complex to catch and many iterations are needed which means that submatrices $S_i$ of low dimensions are required.

As for the previous numerical experiments, we consider the values of $N$ for which the error $\mathcal{E}_p$ has its maximum and minimum, that is $N=230$ where $\mathcal{E}_p(\widetilde{S}^{230}, {\bf r}) = 2.6059 \times 10^{-4}$ and $N=535$ where the pDMD stopped. In the right panel of Figure \ref{dib_hopf_err_piecewise}, we show the behaviour of the error $\epsilon_k(\widetilde{S}^N,{\bf r})$ (see \eqref{error_time}) along the integration time interval of the PDE model. For both partition sizes $N$, this error is larger in the second part of the time interval that corresponds to the oscillatory behaviour of the spatial mean (as shown in the bottom left plot of Figure \ref{dib_hopf_mean_classic}), but it uniformly decreases on the whole time interval for the larger value $N=535$.

\begin{figure}[htbp]
\centering
\includegraphics[scale=0.45]{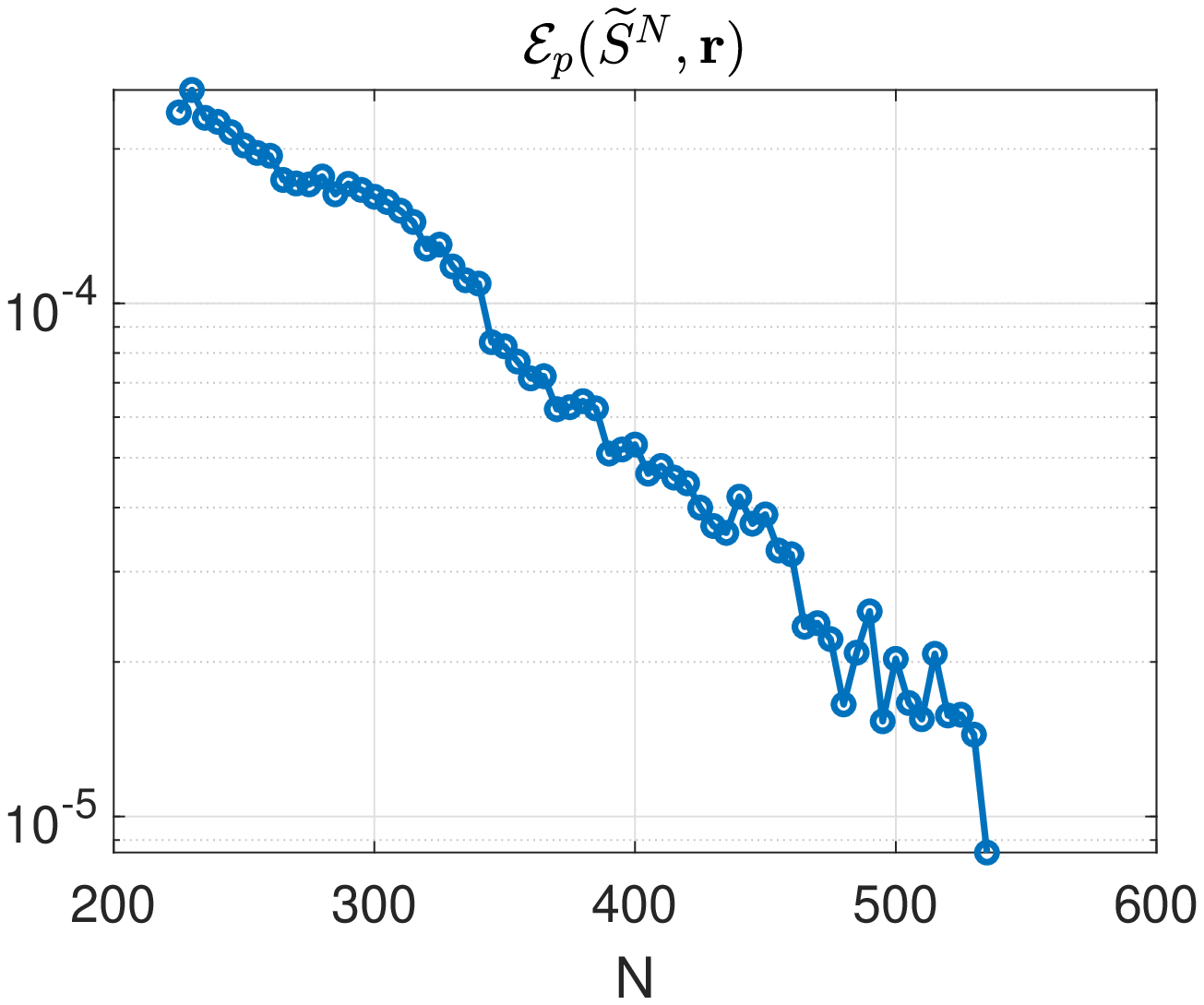}
\includegraphics[scale=0.45]{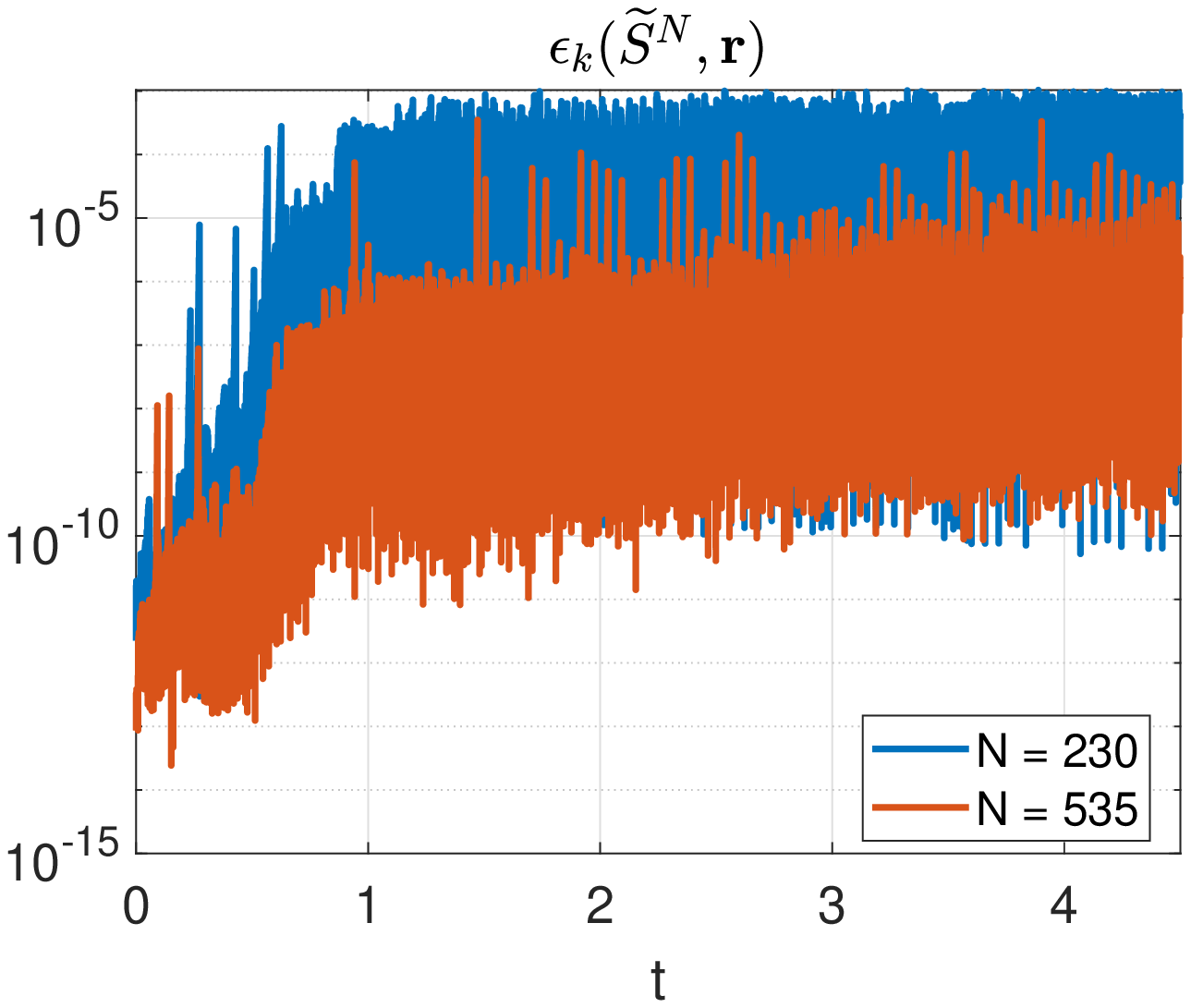}
\captionsetup{justification=justified}
\caption{DIB model: Turing-Hopf instability. Left plot: relative error $\mathcal{E}_p(\widetilde{S}^N, {\bf r})$ in \eqref{rel_err}. Relative error $\epsilon_k(\widetilde{S}^N, {\bf r})$ over time in \eqref{error_time} for two meaningful values of $N$.}
\label{dib_hopf_err_piecewise}
\end{figure}
\begin{figure}[htbp]
\centering
\includegraphics[scale=0.4]{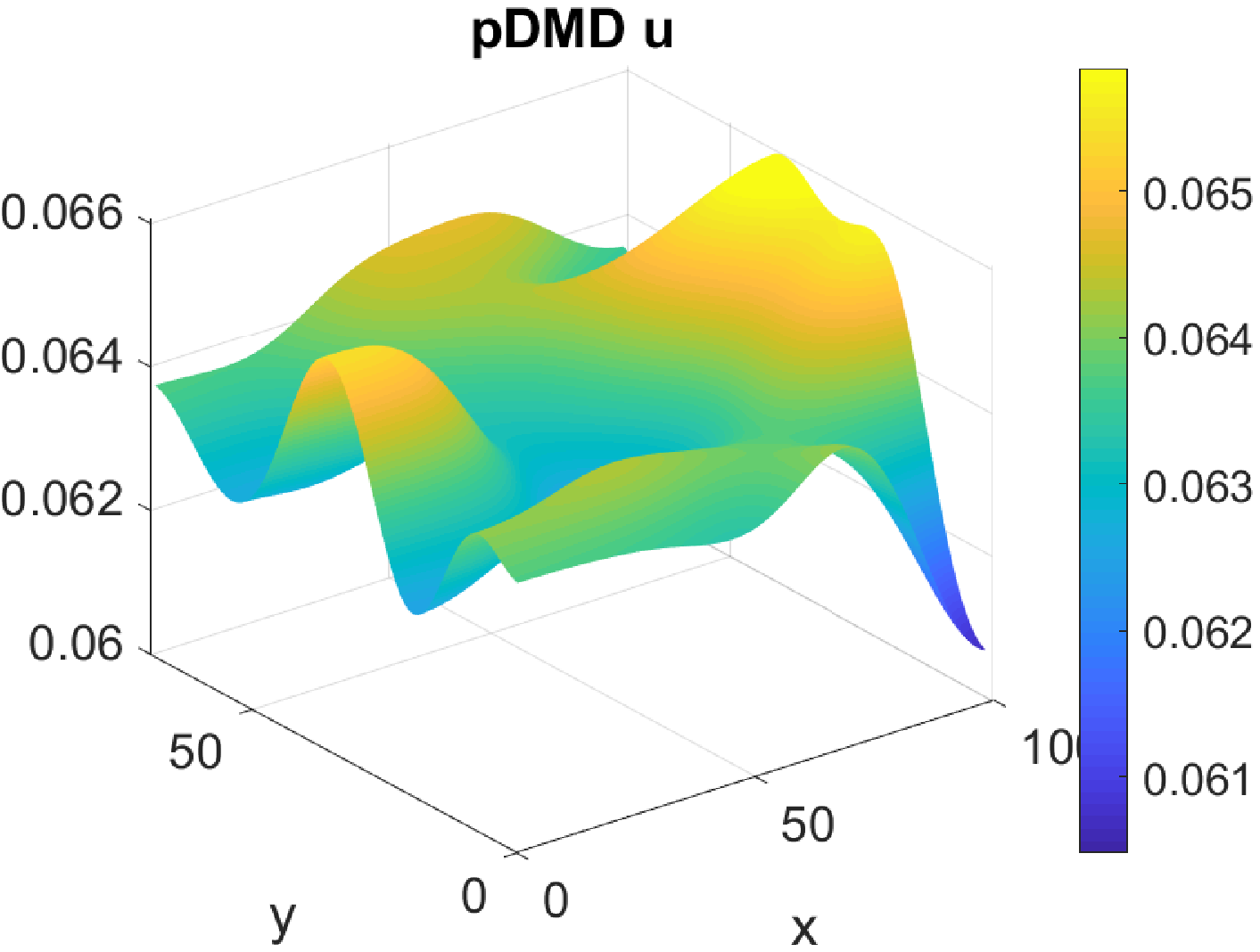}
\includegraphics[scale=0.4]{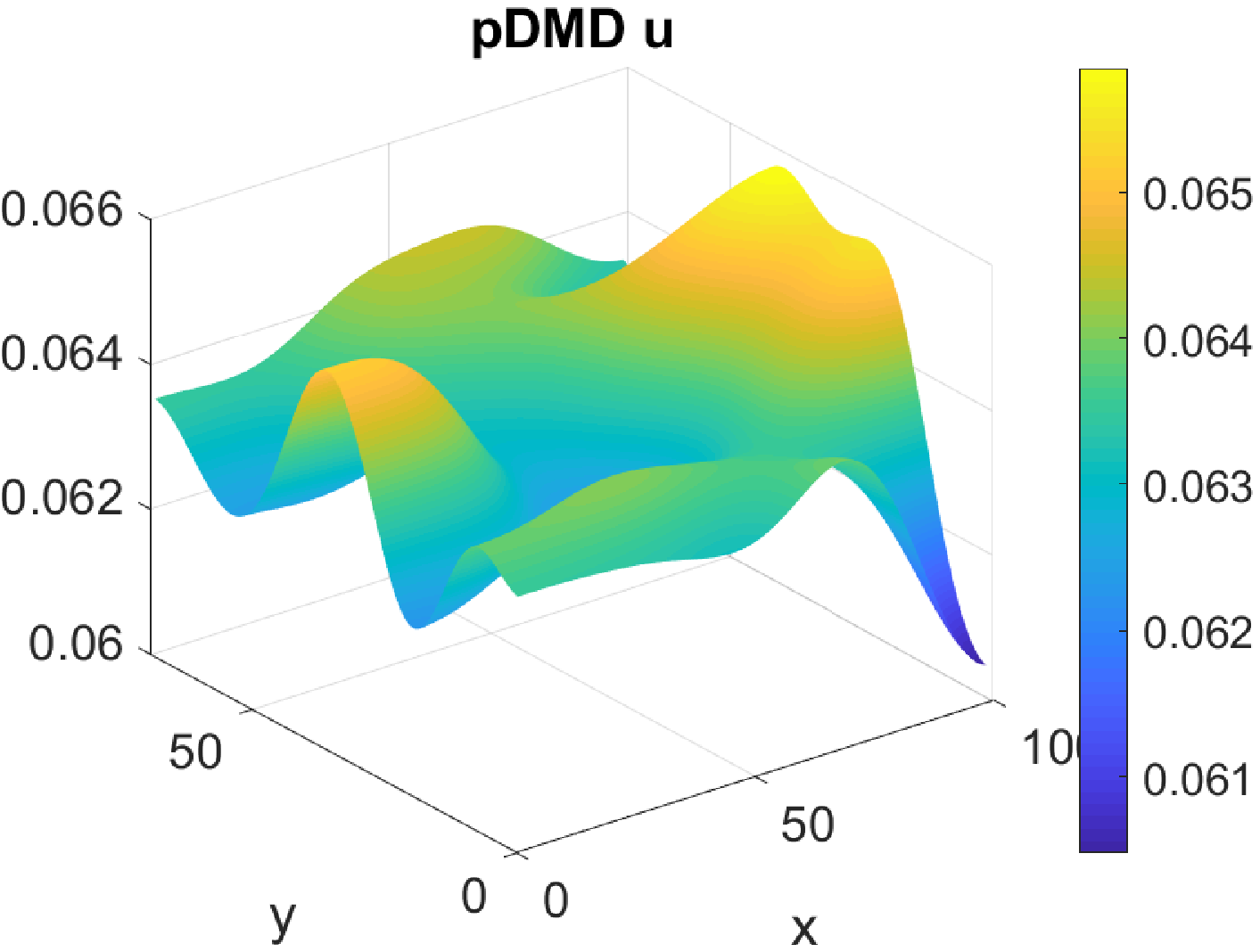}
 \includegraphics[scale=0.4]{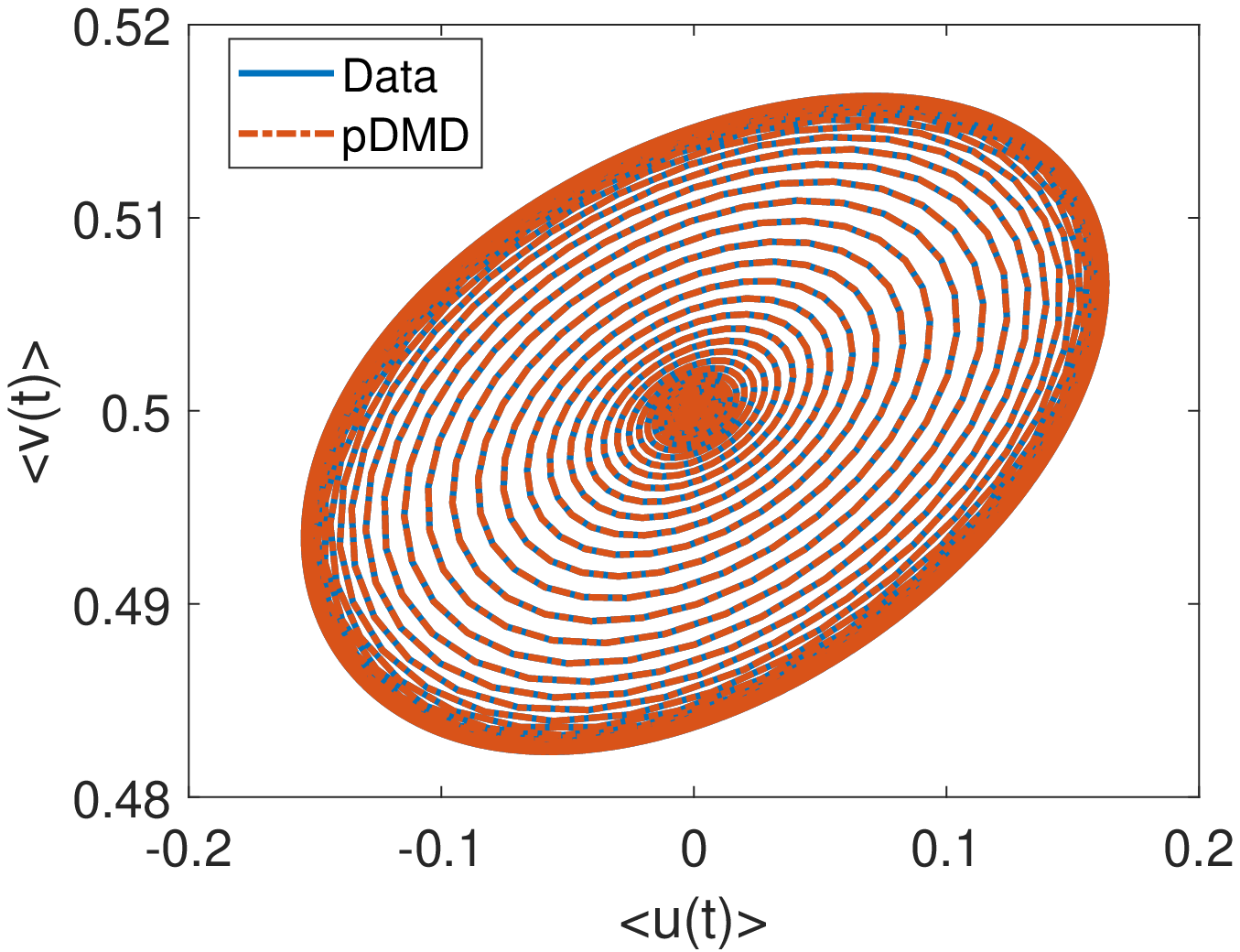}
\captionsetup{justification=justified}
\caption{DIB model: Turing-Hopf instability. pDMD reconstruction of the variable $u$ at the final time $T=4.5$ for $N = 230$ (left plot) and $N =535$ (center plot). Right plot: limit cycle for the pDMD solution with $N = 535$ in the phase plane of the spatial means, for which the error has its minimum that is $\mathcal{E}_p(\widetilde{S}^N,{\bf r}) = 8.5031 \times 10^{-6}$.}
\label{dib_hopf_sol_phase_piecewise}
\end{figure}

\noindent
In Figure \ref{dib_hopf_sol_phase_piecewise}, we compare the solutions obtained at the final time $T$ for the $u$ variable by applying the pDMD algorithm with $N = 230$ (left plot) and $N = 535$ (middle plot). 
The reconstruction with $N = 535$ (middle plot) is almost the same as the data (see Figure \ref{dib_hopf_dmd_classic}); for $N = 230$ the almost same shape is reconstructed by pDMD (left plot), with very small different amplitude, as evident by a slight different colour distribution.
Moreover, in Figure \ref{dib_hopf_sol_phase_piecewise} (right plot) we compare the time dynamics in the phase plane $(\langle u \rangle, \langle v \rangle)$ of the spatial mean values  obtained by pDMD for $N = 535$ with respect to the data. We can observe that there is no difference between the limit cycles for pDMD and data.
 
\begin{figure}[htbp]
\centering
\includegraphics[scale=0.45]{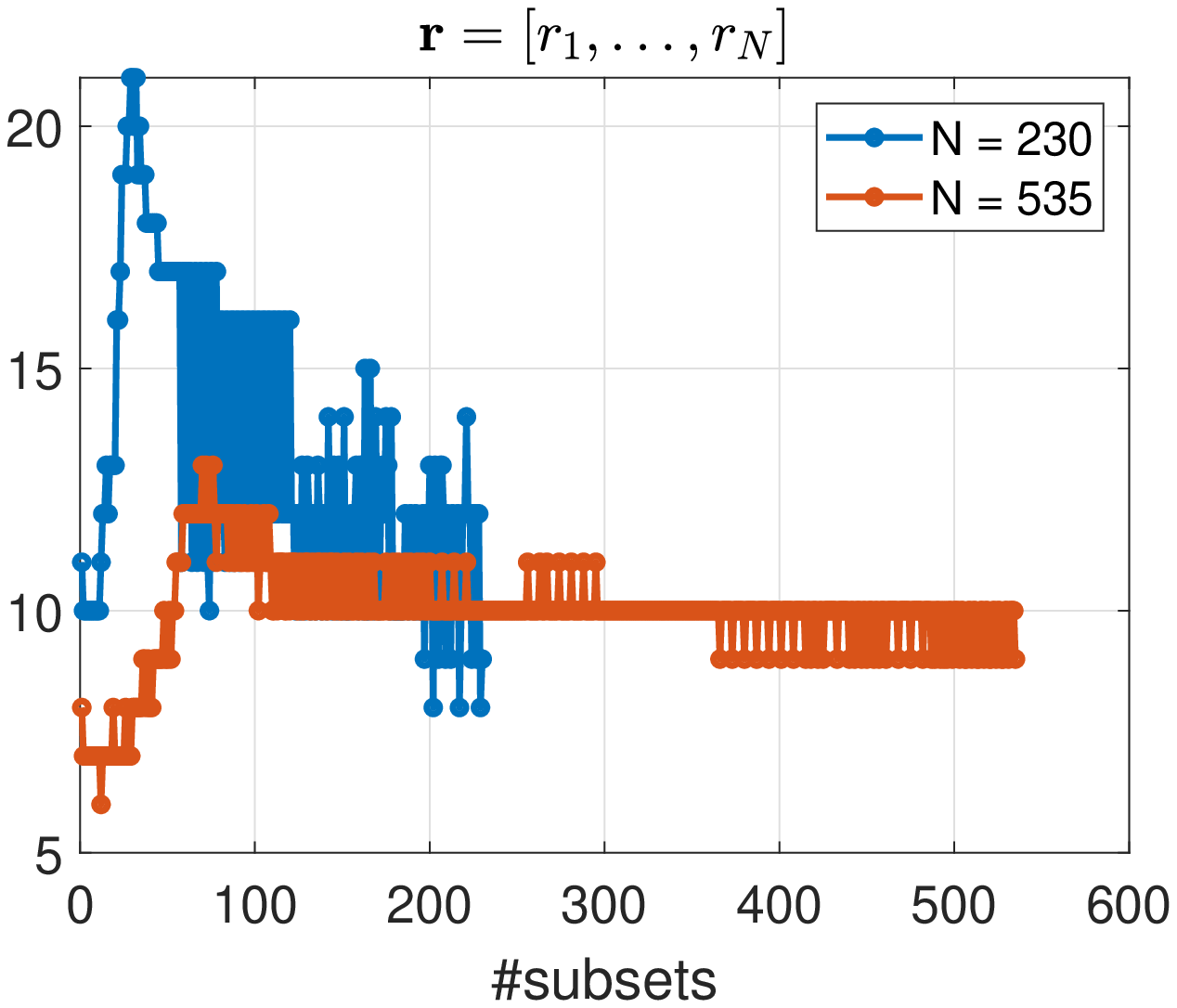}
\includegraphics[scale=0.45]{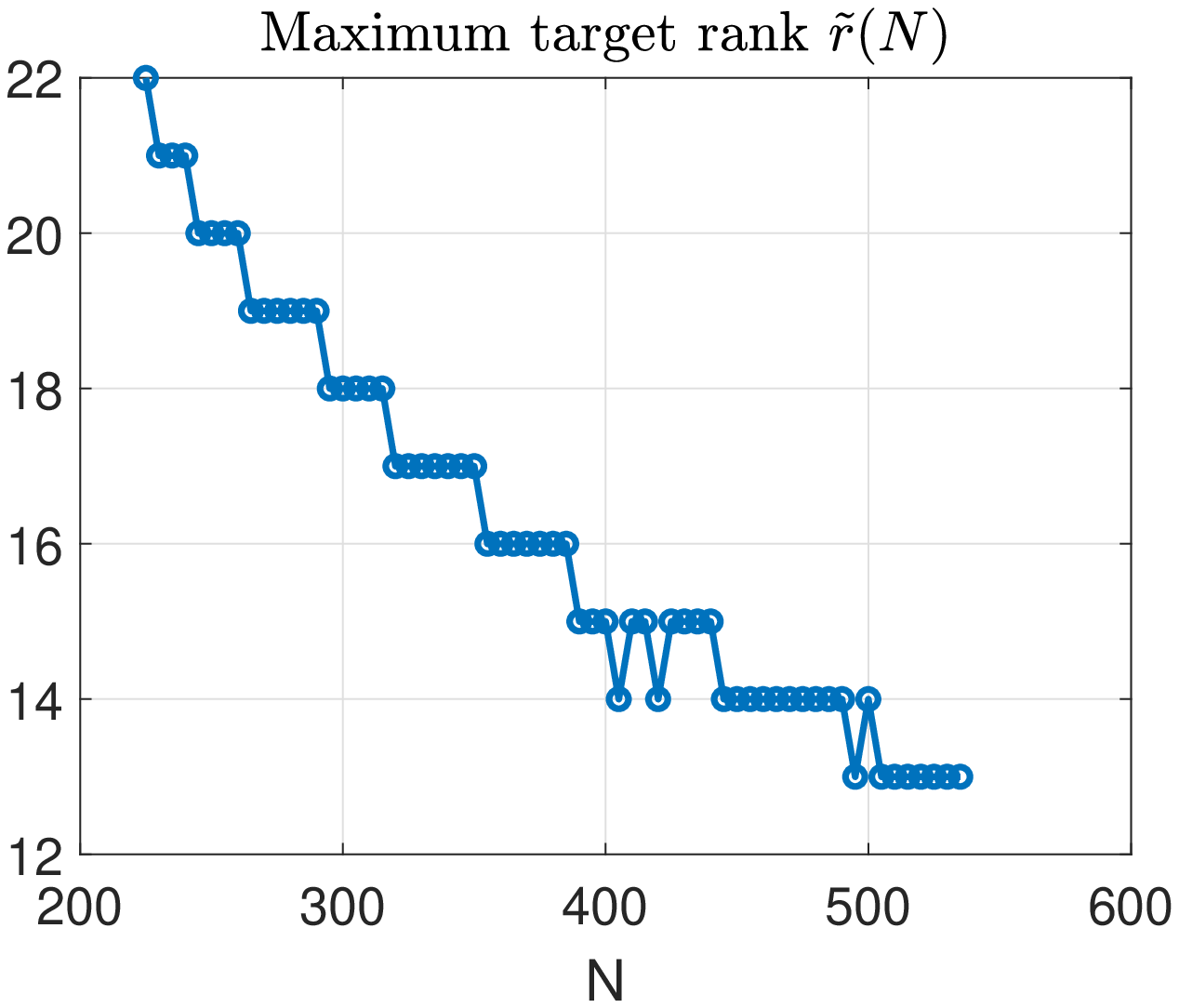}
\captionsetup{justification=justified}
\caption{DIB model: Turing-Hopf instability. Left plot: vector $\r$ of the target ranks $r_i$ used in the subsets $S_i$ for $N=230,535$. Right plot: maximum target rank $\tilde{r}(N)$ with respect to the range of partition sizes $N$ identified in Figure \ref{dib_hopf_err_piecewise}.}
\label{dib_hopf_rank_pdmd}
\end{figure}

\noindent
Finally, in the right panel of Figure \ref{dib_hopf_rank_pdmd}, always for the partition sizes $N = 230$ and $N =535$ we show the ranks $r_i$ chosen in the subsets $S_i$, for $i=1,\dots,N$.  In both cases, in the first part of the time interval, increasing values of the target rank are needed, then after a certain subset (or time, say $\bar{t} \approx 1$ when the limit cycle is reached) smaller and smaller rank values are sufficient to follows the oscillatory spatio-temporal regime. Moreover, for larger $N$, that is $N = 535$, smaller $r_i$ are needed. To confirm this trend, in the right panel of Figure \ref{dib_hopf_rank_pdmd} we show the behaviour of the maximum target rank $\tilde{r}(N)$: it exhibits a monotone decay and in the worst case, the maximum rank needed is $22 \ll 130 = {\tt rank}(S)$, i.e. much smaller than the rank of the original dataset.

\section{Conclusions}\label{sec:end}
Inspired by the classical ``divide and conquer'' principle, in this paper we have introduced a {\it piecewise} version of the exact DMD technique, called pDMD.
Given a temporal dataset and fixed tolerance, the new approach is implemented in Algorithm \ref{Alg:PiecDMD} and Algorithm \ref{Alg:PiecDMD_N}, when the partition size $N$ of the original dataset is increased towards a desired final accuracy of the reconstruction. 
The new approach can be applied to a general dataset, even if our study has been motivated by the failure of the original exact DMD on snapshots describing peculiar spatio--temporal dynamics arising in Reaction-Diffusion (RD) PDE systems. In fact, in Section \ref{sec_fail_dmd} we have shown that DMD exhibits very innacurate reconstructions or ill-conditioning for large target ranks (where a better approximation is expected) for four significant models with: relaxation oscillations (FitzHugh-Nagumo 1D in space), spiral waves ($\lambda$-$\omega$ system), Turing pattern formation and Turing-Hopf patterns oscillating in space and time (DIB morphochemical system for battery modeling).

In Section \ref{sec:testpdmd}, for each kind of the above dynamics, we have shown that pDMD is able to remove all drawbacks previously highlighted.
In some cases, like the FitzHugh-Nagumo and Turing-Hopf dynamics, we have shown that a suitable partition size $N^*$ can be obtained such that for $N \geq N^*$ an error much lower than the best obtained by the ``global'' DMD in Section \ref{sec_fail_dmd} is obtained.
In other cases, like for spiral waves and Turing instability, both the final spiral/pattern and their time histories described by the limit cycle in the phase plane are carefully reconstructed for $N \geq N^*$.
In all cases, pDMD is now able to follow the entire spatio-temporal dynamics, including different regimes (e.g. reactivity-stabilizing for Turing, unstable-oscillating for spiral waves). In particular, in all simulations, for larger $N$ the error uniformly decreases along all the time interval of the entire dataset and a convergence trend can be observed.

We can conclude that, for oscillatory and Turing spatio-temporal dynamics, 
the main assumption underlying the original DMD, that is a ``global'' linear fitting over the full temporal horizon, is not sufficient to recognize different ``phenomena'' arising along the time pathways. Instead, the ``local'' linear fitting by pDMD does it when a sufficient dataset partition size $N$ can be identified.

As a final remark, in all simulation in Section \ref{sec:testpdmd} and for all $N$, we find that the maximum rank considered by the DMDs along the partition is always much smaller than the rank of the original dataset $S$. This result confirms that the piecewise approach is also able to reduce the overall computational load of the usual DMD implementation: several problems of significant small dimensions are solved by the linear fitting behind DMD. We argue that this property is also at the origin of the observed ill-conditioning reduction with respect to the original DMD (see Figure \ref{lo_mean_classic} and \ref{dib_mean_classic}).

\section*{Acknownledgments}
AA, AM, IS  are members of the INdAM-GNCS activity group. The work of IS is
supported by the MIUR through the project PRIN 2020, “Mathematics for Industry 4.0”, project no. 2020F3NCPX and from “National Centre for High Performance Computing, Big Data and Quantum Computing” funded by European Union – NextGenerationEU, PNRR project code CN00000013, CUP F83C22000740001.

\bibliographystyle{plain}
\bibliography{refs_new.bib}

\end{document}